\documentclass[12pt,reqno]{amsart}

\usepackage{amssymb}
\usepackage{enumerate}
\usepackage{tikz}
\usepackage[margin=2.2cm]{geometry}

\usepackage[T1]{fontenc}
\usepackage[english]{babel}
\usepackage{lmodern,microtype,hyperref}

\allowdisplaybreaks
\DeclareMathOperator{\Ad}{Ad}
\DeclareMathOperator{\Area}{Area}
\DeclareMathOperator{\Card}{Card}
\DeclareMathOperator{\Op}{Op}
\DeclareMathOperator{\Rer}{Re}
\DeclareMathOperator{\Vol}{Vol}

\newcommand{\bs}{\backslash}
\newcommand{\C}{\mathbb{C}}
\newcommand{\D}{\mathbb{D}}
\newcommand{\Dc}{\mathcal{D}}
\newcommand{\F}{\mathbb{F}}
\newcommand{\Hb}{\mathbb{H}}
\newcommand{\Hc}{\mathcal{H}}
\newcommand{\Id}{\mathrm{I}}

\newcommand{\Lc}{\mathcal{L}}
\newcommand{\M}{\mathbb{M}}
\newcommand{\N}{\mathbb{N}}
\newcommand{\ol}{\overline}
\newcommand{\PSU}{\mathrm{PSU}}
\newcommand{\R}{\mathbb{R}}
\newcommand{\Sb}{\mathbb{S}}
\newcommand{\Sc}{\mathcal{S}}
\newcommand{\SL}{\mathrm{SL}}
\newcommand{\SO}{\mathrm{SO}}
\newcommand{\SU}{\mathrm{SU}}
\newcommand{\suf}{\mathfrak{su}}
\newcommand{\T}{\mathbb{T}}
\newcommand{\ve}{\varepsilon}
\newcommand{\vp}{\varphi}
\newcommand{\Wc}{\mathcal{W}}
\newcommand{\Wr}{\mathrm{W}}
\newcommand{\Z}{\mathbb{Z}}

\theoremstyle{plain}
\newtheorem{thm}{Theorem}[section]
\newtheorem*{conj}{Conjecture}
\newtheorem{lemma}[thm]{Lemma}
\newtheorem{prop}[thm]{Proposition}
\newtheorem*{thmn}{Theorem}

\theoremstyle{definition}
\newtheorem*{defn}{Definition}
\newtheorem*{rmk}{Remark}

\numberwithin{equation}{section}

\title[Schr\"odinger group along the horocycles and Egorov's theorem]{Schr\"odinger group along the horocycles and Egorov's theorem on compact hyperbolic surfaces}
\author{Xiaolong Han}
\email{xiaolong.han@csun.edu}
\address{Department of Mathematics, California State University, Northridge, CA 91330, USA}

\subjclass[2020]{58J51, 81Q50 (primary); 37A17, 37D40, 58J40, 81Q20 (secondary)}
\keywords{Schr\"odinger group along horocycles, Egorov's theorem, Ehrenfest time, quantum ergodicity, quantum unique ergodicity, hyperbolic surfaces, horocycle flow}
\thanks{}

\begin{document}

\begin{abstract}
In this paper, we study Egorov's theorem for Laplacian eigenfunctions on compact hyperbolic surfaces in Zelditch's quantization. In this quantization, the remainder in Egorov's theorem for eigenfunctions is given explicitly by the second derivatives of the symbols along the weak stable foliation of the geodesic flow. Forward in time, we show that Egorov's theorem for eigenfunctions is valid to all logarithmic times, far beyond the Ehrenfest time $T_E=\frac12|\log h|$, where $h\to0$ is the semiclassical parameter.

Backward in time, we introduce the Schr\"odinger group along the stable horocycles, which absorbs the exponentially growing part of the remainder. We show that Egorov's theorem for eigenfunctions is valid up to $-T_E+T_h$ for any $T_h\to\infty$ with $0\le T_h\le T_E$. Moreover, the limit of the remainder in Egorov's theorem at the Ehrenfest time $-T_E$ is exactly the difference between the semiclassical measure and the Liouville measure, tested against the symbols propagated by the Schr\"odinger group. Consequently, Egorov's theorem for a sequence of eigenfunctions holds at $-T_E$ (and in fact, at $-T_E-T$ for any fixed $T\in\R$) if and only if the sequence tends to be equidistributed in the phase space. Thus, the quantum unique ergodicity (QUE) conjecture is equivalent to this Egorov property holding for every eigenbasis. Furthermore, QUE holds with a polynomial rate if and only if Egorov's theorem holds at $-T_E$ with a polynomial remainder. In this case, Egorov's theorem holds beyond $-T_E$ with a polynomial remainder, and the stronger the rate is in QUE, the longer the time is in Egorov's theorem.
\end{abstract}

\maketitle

\setcounter{tocdepth}{1}
\tableofcontents

\section{Introduction}\label{sec:intro}
Let $\M=\Gamma\bs\D$ be a compact connected hyperbolic surface, where $\D$ is the Poincar\'e disk and $\Gamma\subset\PSU(1,1)$ is a torsion-free cocompact Fuchsian group. Suppose that $\{\vp_j\}_{j=0}^\infty$ is an orthonormal basis of eigenfunctions (i.e., eigenbasis) of the non-negative Laplacian $\Delta_\Gamma$ in $L^2(\M)$, with eigenvalues
\[0=\mu_0<\mu_1\le\mu_2\le\cdots\to\infty.\]
If $\mu_j>\frac54$, then we write $\mu_j=\lambda_j^2+\frac14$ with $\lambda_j>1$, and use $h_j=\lambda_j^{-1}\in(0,1)$ as the semiclassical parameter.

The distribution of the eigenfunction $\vp_j$ in the phase space, i.e., the tangent bundle $T\M=\{(x,\xi):x\in\M,\xi\in T_x\M\}$, is described by the matrix elements
\[\left\langle\Op_{h_j}^\Gamma(f)\vp_j,\vp_j\right\rangle_{L^2(\M)},\]
where $f$ is a symbol on $T\M\setminus\{0\}$ and $\Op_h^\Gamma(f)$ is the semiclassical pseudo-differential operator constructed by Zelditch \cite{Ze1, Ze2}. The quantum unique ergodicity (QUE) conjecture of Rudnick-Sarnak \cite{RuSa} asserts that
\begin{conj}[The quantum unique ergodicity conjecture]
\[\lim_{j\to\infty}\left\langle\Op_{h_j}^\Gamma(f)\vp_j,\vp_j\right\rangle_{L^2(\M)}=\int_{S\M}f\,d\Vol_L\quad\text{for each }f\in S^0(T\M\setminus\{0\})\text{ independent of }h,\]
where $d\Vol_L$ is the (normalized) Liouville measure on the sphere bundle $S\M=\{(x,\xi)\in T\M:|\xi|_x=1\}$ and $S^0(T\M\setminus\{0\})$ is the symbol class of order $0$.
\end{conj}

The QUE conjecture predicts that the eigenfunctions tend to be equidistributed on $S\M$. In particular, if $f\in C^\infty(\M)$, then $\Op_{h_j}^\Gamma(f)=f$, so the QUE conjecture becomes
\[\lim_{j\to\infty}\int_\M f\left|\vp_j\right|^2\,d\Area_L=\int_\M f\,d\Area_L\quad\text{for each }f\in C^\infty(\M),\]
that is, the eigenfunctions tend to be equidistributed in the physical space $\M$ with respect to the (normalized) Riemannian area. The QUE conjecture can be proposed on compact manifolds with variable negative curvature, and it remains open; see Dyatlov \cite{D1, D2} and Nonnenmacher \cite{N} for recent developments related to the QUE conjecture.

The main purpose of this paper is to derive a precise relation between the QUE conjecture on compact hyperbolic surfaces and Egorov's theorem for eigenfunctions in Zelditch's quantization \cite{Ze1, Ze2, Ze3}. We remark that the QUE conjecture does not depend on the choice of quantization, but our discussion on Egorov's theorem and its relation to the QUE conjecture does; see Section \ref{sec:quant} for more discussion.

Egorov's theorem \cite{E} relates the quantum evolution of semiclassical pseudo-differential operators to the classical evolution of their symbols under the geodesic flow. For standard semiclassical quantizations, an energy-localized version holds on suitable logarithmic time intervals; see Zworski \cite[Chapter 11]{Zw}.

Let $a_t$ be the unit-speed geodesic flow on $T\M\setminus\{0\}$. In this paper, we investigate the following form of Egorov's theorem for eigenfunctions $\vp_j$ with $h=h_j$ when $\mu_j>\frac54$, initially for fixed $t\in\R$:
\[\left|\left\langle\Op_{h_j}^\Gamma\left(f\circ a_t\right)\vp_j,\vp_j\right\rangle_{L^2(\M)}-\left\langle\Op_{h_j}^\Gamma(f)\vp_j,\vp_j\right\rangle_{L^2(\M)}\right|=o_{f,t}(1)\quad\text{as }h_j\to0.\]
This formulation concerns the diagonal matrix elements of the eigenfunctions. In Section \ref{sec:SDO}, we show that the matrix elements depend only on the restriction of the symbol $f$ to $S\M$ at $h=h_j$. Conversely, every function in $C^\infty(S\M)$ is the restriction to $S\M$ of a symbol in $S^0(T\M\setminus\{0\})$, e.g., its extension which is homogeneous of degree $0$. Therefore, by this identification, we can write
\[\Wc_j(f)=\left\langle\Op_{h_j}^\Gamma(f)\vp_j,\vp_j\right\rangle_{L^2(\M)}\quad\text{for }f\in C^\infty(S\M).\]
The symbol $f$ may depend on $h$, in which case it is evaluated at $h=h_j$. The map $f\mapsto\Wc_j(f)$ defines a distribution on $S\M$. In particular, $\Wc_j(f\circ a_t)$ is defined for each $t\in\R$, even when $t$ depends on $h_j$. Therefore, for eigenfunctions, Egorov's theorem concerns the difference
\[\Wc_j\left(f\circ a_t\right)-\Wc_j(f).\]
We investigate the time ranges in which this difference tends to $0$ as $j\to\infty$.

Since the geodesic flow $a_t$ is Anosov with expansion rate $1$, derivatives of order $k$ of $f\circ a_t$ may grow like $e^{k|t|}$. Thus, the standard Egorov theorem is usually valid up to the Ehrenfest time
\[T_E=\frac12|\log h|.\]
Here and below, $T_E=\frac12|\log h_j|$ for the eigenfunction $\vp_j$. Our results show that the Ehrenfest time $T_E$ is exactly the threshold at which Egorov's theorem for eigenfunctions detects QUE:
\begin{itemize}
\item Egorov's theorem for eigenfunctions is valid forward in time to all logarithmic times, and backward in time up to $-T_E+T_h$ for any $T_h\to\infty$ with $0\le T_h\le T_E$ as $h\to0$ (Theorem \ref{thm:Egorov}).
\item At the Ehrenfest time $-T_E$, the limit of the remainder in Egorov's theorem is exactly the difference between the semiclassical measure and the Liouville measure, tested against the symbols propagated by the Schr\"odinger group along the stable horocycles defined below. So Egorov's theorem for a sequence of eigenfunctions holds at $-T_E$ if and only if the sequence tends to be equidistributed in the phase space, in which case it also holds beyond $-T_E$ (Theorem \ref{thm:equiv}).
\item Egorov's theorem for a sequence of eigenfunctions holds at $-T_E$ with a polynomial remainder if and only if the sequence tends to be equidistributed in the phase space with a polynomial rate. In this case, it also holds beyond $-T_E$ with a polynomial remainder, and a stronger rate of equidistribution permits a longer time range (Theorem \ref{thm:rate}).
\end{itemize}
In particular, the QUE conjecture on compact hyperbolic surfaces is equivalent to Egorov's theorem for eigenfunctions at the Ehrenfest time.

Our starting point is that, in Zelditch's quantization, the remainder in Egorov's theorem for eigenfunctions is given explicitly (Theorem \ref{thm:EefnM}): for $\mu_j>\frac54$, $f\in C^\infty(S\M)$, and $t\in\R$,
\begin{equation}\label{eq:exact}
\Wc_j\left(f\circ a_t\right)-\Wc_j(f)=\frac{ih_j}{2}\int_0^t\Wc_j\left(\left(Q^2f\right)\circ a_s+e^{-2s}\left(U_+^2f\right)\circ a_s\right)\,ds,
\end{equation}
where $Q$ generates the geodesic flow and $U_+$ generates the stable horocycle flow. That is, the remainder only involves the second derivatives of $f$ along the leaves of the weak stable foliation. Forward in time, the weight $e^{-2s}$ is bounded by $1$, so \eqref{eq:exact} can be iterated directly. Backward in time, the weight grows and reaches $h_j^{-1}$ at $s=-T_E$, so the term $e^{-2s}(U_+^2f)\circ a_s$ is no longer small. This term is absorbed by the central object of this paper, the Schr\"odinger group along the stable horocycles. Formally, this group is defined by
\begin{equation}\label{eq:Schr}
e^{-\frac{i\sigma}{4}U_+^2}f=u(\sigma)\quad\text{for }\sigma\in\R,f\in C^\infty(S\M),
\end{equation}
where, along each stable horocycle $s\mapsto xn_s$ with $n_s=e^{sU_+}$, $u$ solves the free Schr\"odinger equation
\[\partial_\sigma u(\sigma)=-\frac i4\partial_s^2u(\sigma)\quad\text{and}\quad u(0)=f.\]
We give the precise definition by the spectral theorem in Section \ref{sec:Schrgroup}. The Schr\"odinger group is a unitary group on $L^2(S\M)$, which maps $C^\infty(S\M)$ to itself, commutes with the stable horocycle flow, and preserves $d\Vol_L$ (Lemma \ref{lemma:Schr}). It acts on the symbols, rather than on the eigenfunctions. We explain its role in Section \ref{sec:Schr}.

For $f\in C^\infty(S\M)$, we write $\|f\|_{C^n(S\M)}$ and $\|f\|_{H^n(S\M)}$ for the $C^n$ and Sobolev norms of $f$ with respect to the frame $\{Q,U_+,U_-\}$, respectively; see Section \ref{sec:norms}. Our first main theorem states
\begin{thm}[Egorov's theorem forward and backward in time]\label{thm:Egorov}
Let $f\in C^\infty(S\M)$. Suppose that $j\in\N$ satisfies $\mu_j>\frac54$.
\begin{enumerate}[(i)]
\item Let $\rho>0$ and $k=\lfloor\rho\rfloor+2$. Then there is $C=C(\rho)>0$ such that for $0\le t\le\rho T_E$,
\[\left|\Wc_j\left(f\circ a_t\right)-\Wc_j(f)\right|\le Ch_j\|f\|_{C^{2k+2}(S\M)}.\]
\item There is $C=C(\M)>0$ such that for $-T_E\le t\le0$,
\[\left|\Wc_j\left(f\circ a_t\right)-\Wc_j(f)\right|\le Ch_je^{2|t|}\|f\|_{H^{20}(S\M)}.\]
\end{enumerate}
\end{thm}

Therefore, (i) of Theorem \ref{thm:Egorov} implies that forward in time, Egorov's theorem for eigenfunctions holds with a remainder $O_{f,\rho}(h)$ up to $\rho T_E$ for every $\rho>0$, far beyond the Ehrenfest time $T_E$. Its proof does not use the Schr\"odinger group: we iterate \eqref{eq:exact} $\lfloor\rho\rfloor+2$ times, and use a uniform bound of $\Wc_j$ involving only two derivatives (Lemma \ref{lemma:apriori}).

Egorov's theorem for all times also holds for the quantum cat maps of Hannay-Berry \cite{HB}, where it is exact in the Weyl quantization because the cat maps are linear. However, QUE fails for the quantum cat maps by Faure-Nonnenmacher-De Bi\`evre \cite{FNDB}, so Egorov's theorem in the Weyl quantization does not detect QUE. An analogue of Theorem \ref{thm:equiv} below holds in a quantization adapted to the stable foliation of the cat map \cite{Han2}; see Appendix \ref{sec:non}.

In our setting, for non-constant $f$ and $t\ge2T_E$ (allowed when $\rho\ge2$), the symbols $f\circ a_t$ do not belong to any reasonable symbol class, including the exotic ones in Zworski \cite{Zw}, Dyatlov-Zahl \cite{DZ}, and Dyatlov-Jin \cite{DyJi}. Thus, $\Op_{h_j}^\Gamma(f\circ a_t)$ is no longer a semiclassical pseudo-differential operator in any known calculus, although its action on eigenfunctions is defined by Theorem \ref{thm:SDOefnM}, and $\Wc_j(f\circ a_t)$ is defined since $\Wc_j$ is a distribution. In this sense, (i) of Theorem \ref{thm:Egorov} is a property of the eigenfunctions rather than of the operators.

Backward in time, (ii) of Theorem \ref{thm:Egorov} gives the remainder
\[O_f\left(e^{-2T}\right)\quad\text{at }t=-T_E+T\text{ with }0\le T\le T_E.\]
So Egorov's theorem with a remainder of order $O_f(e^{-2T_h})=o_f(1)$ for eigenfunctions holds up to $-T_E+T_h$ for any $T_h\to\infty$ with $0\le T_h\le T_E$ as $h\to0$. The difference between the two time directions comes from the Lagrangian foliation associated with Zelditch's quantization; see Section \ref{sec:quant}.

At the Ehrenfest time $t=-T_E$, (ii) of Theorem \ref{thm:Egorov} gives only an $O_f(1)$ bound, which does not imply that the remainder tends to $0$. Our second main theorem computes the limit of the remainder in Egorov's theorem at and beyond $-T_E$ in terms of the semiclassical measures. To describe the distributions of eigenfunctions in the phase space, we recall the notion of a semiclassical measure. Let $J\subset\N$ be an infinite set such that $\mu_j>\frac54$ for $j\in J$. We say that $\{\vp_j\}_{j\in J}$ induces a semiclassical measure $\nu$ if
\[\lim_{j\in J,j\to\infty}\Wc_j(f)=\int_{S\M}f\,d\nu\quad\text{for each }f\in C^\infty(S\M)\text{ independent of }h.\]
It is well known that each semiclassical measure is a probability measure on $S\M$ and is invariant under the geodesic flow $a_t$; see Zworski \cite[Sections 5.1 and 5.2]{Zw} and also Section \ref{sec:scm}. We say that $\{\vp_j\}_{j\in J}$ tends to be equidistributed in the phase space if every semiclassical measure induced by its subsequences is $d\Vol_L$. Hence, QUE is equivalent to the statement that the semiclassical measure is unique and equals $d\Vol_L$.

\begin{thm}[Egorov's theorem at the Ehrenfest time and semiclassical measures]\label{thm:equiv}
Let $J\subset\N$ be an infinite set such that $\mu_j>\frac54$ for $j\in J$. Suppose that $\{\vp_j\}_{j\in J}$ induces a semiclassical measure $\nu$. Then for $T\in\R$ and $f\in C^\infty(S\M)$ independent of $h$,
\begin{equation}\label{eq:defect}
\lim_{j\in J,j\to\infty}\left(\Wc_j\left(f\circ a_{-T_E-T}\right)-\Wc_j(f)\right)=\int_{S\M}\left(e^{-\frac{ie^{2T}}{4}U_+^2}f-f\right)d\nu.
\end{equation}
Moreover, for each $T\in\R$, the limits in \eqref{eq:defect} for $f\in C^\infty(S\M)$ determine $\nu$. In particular, the following are equivalent.
\begin{enumerate}[(i)]
\item $\nu=d\Vol_L$.
\item Egorov's theorem holds at the Ehrenfest time $-T_E$ along $J$, i.e., the limits in \eqref{eq:defect} with $T=0$ vanish for each $f\in C^\infty(S\M)$.
\item Egorov's theorem holds beyond the Ehrenfest time along $J$, i.e., for each $f\in C^\infty(S\M)$ independent of $h$, there are $T_j\ge0$, $j\in J$, with $T_j\to\infty$ such that
\[\lim_{j\in J,j\to\infty}\sup_{-T_E-T_j\le t\le0}\left|\Wc_j\left(f\circ a_t\right)-\Wc_j(f)\right|=0.\]
\end{enumerate}
\end{thm}

Since the Schr\"odinger group preserves $d\Vol_L$, the right-hand side of \eqref{eq:defect} vanishes if $\nu$ is replaced by $d\Vol_L$. That is, the limit of the remainder in Egorov's theorem at $-T_E-T$ is exactly the difference between $\nu$ and $d\Vol_L$, tested against $e^{-\frac{ie^{2T}}{4}U_+^2}f-f$.

Every subsequence of $\{\vp_j\}_{j\in J}$ has a further subsequence which induces a semiclassical measure; see Section \ref{sec:scm}. So by Theorem \ref{thm:equiv}, $\{\vp_j\}_{j\in J}$ tends to be equidistributed in the phase space if and only if Egorov's theorem holds at the Ehrenfest time, i.e.,
\[\lim_{j\in J,j\to\infty}\left(\Wc_j\left(f\circ a_{-T_E}\right)-\Wc_j(f)\right)=0\quad\text{for each }f\in C^\infty(S\M)\text{ independent of }h,\]
in which case Egorov's theorem holds beyond the Ehrenfest time as in (iii). In fact, the time $-T_E$ can be replaced by $-T_E-T$ for any fixed $T\in\R$, i.e., by any fixed time before or beyond $-T_E$ at a bounded distance; see Section \ref{sec:beyond}. In particular, taking $J=\{j\in\N:\mu_j>\frac54\}$, the fixed eigenbasis $\{\vp_j\}_{j=0}^\infty$ tends to be equidistributed in the phase space if and only if Egorov's theorem holds for this eigenbasis at $-T_E$. Hence, the QUE conjecture is equivalent to this Egorov property holding for every Laplacian eigenbasis.

\begin{rmk}[Outline of the proof]
The proof of Theorem \ref{thm:equiv} relies on the following observation: the weight $e^{-2s}$ in \eqref{eq:exact} can be absorbed by the Schr\"odinger group $e^{-\frac{i\sigma}{4}U_+^2}$ with $\sigma=h_j(e^{2|t|}-1)$, which is $O(1)$ at $t=-T_E-T$. So $\Wc_j(f\circ a_{-T_E-T})$ is close to $\Wc_j(e^{-\frac{i\sigma}{4}U_+^2}f)$ with $\sigma\to e^{2T}$, which gives \eqref{eq:defect}. Since $\nu$ is invariant under the geodesic flow, which rescales the time $\sigma$ of the Schr\"odinger group, the limits in \eqref{eq:defect} determine the integrals of $U_+^2f$ against $\nu$. A finite signed measure annihilated by $U_+^2$ is invariant under the stable horocycle flow, so it is a multiple of $d\Vol_L$ by Furstenberg's unique ergodicity of the horocycle flow \cite{Fu}. Hence, these limits determine $\nu$.
\end{rmk}

Our third main theorem is a quantitative version of Theorem \ref{thm:equiv}. We use the effective equidistribution of the stable horocycle flow by Burger \cite{Bur} and Flaminio-Forni \cite{FF} (Theorem \ref{thm:Burger}). That is, the averages of $f\in C^\infty(S\M)$ along the stable horocycles of length $R\ge1$ converge to $\int_{S\M}f\,d\Vol_L$ uniformly with the rate $O_f(R^{-\beta})$. Here, $\beta$ is determined by the first non-zero eigenvalue $\mu_1>0$ of $\Delta_\Gamma$:
\[\beta=\frac12-\Rer\frac{\sqrt{1-4\mu_1}}{2}\quad\text{if }\mu_1\ne\frac14,\]
and if $\mu_1=\frac14$, then $\beta$ can be any number in $(0,\frac12)$.

\begin{thm}[Egorov's theorem with a rate and QUE with a rate]\label{thm:rate}
Let $J\subset\N$ be an infinite set such that $\mu_j>\frac54$ for $j\in J$. The following are equivalent.
\begin{enumerate}[(i)]
\item $\{\vp_j\}_{j\in J}$ tends to be equidistributed in the phase space with a polynomial rate, i.e., there are $\delta_1>0$, $n_1\in\N$, and $C>0$ such that for $j\in J$ and $f\in C^\infty(S\M)$,
\begin{equation}\label{eq:QErate}
\left|\Wc_j(f)-\int_{S\M}f\,d\Vol_L\right|\le Ch_j^{\delta_1}\|f\|_{H^{n_1}(S\M)}.
\end{equation}
\item Egorov's theorem holds at the Ehrenfest time with a polynomial remainder, i.e., there are $\delta_2>0$, $n_2\in\N$, and $C>0$ such that for $j\in J$ and $f\in C^\infty(S\M)$,
\begin{equation}\label{eq:Egorovrate}
\left|\Wc_j\left(f\circ a_{-T_E}\right)-\Wc_j(f)\right|\le Ch_j^{\delta_2}\|f\|_{H^{n_2}(S\M)}.
\end{equation}
\end{enumerate}
More precisely, if (i) holds with $0<\delta_1\le1$ and $n_1\ge1$, then (ii) holds with $\delta_2=\delta_1$ and $n_2=\max\{20,3n_1\}$. Moreover, in this case, Egorov's theorem holds beyond the Ehrenfest time with a polynomial remainder: for each $\ve$ with
\[0<\ve<\min\left\{\frac{1}{13},\frac{\delta_1}{2n_1}\right\},\]
there is $C>0$ such that for $j\in J$ and $f\in C^\infty(S\M)$,
\begin{equation}\label{eq:Egorovbeyond}
\sup_{-(1+\ve)T_E\le t\le0}\left|\Wc_j\left(f\circ a_t\right)-\Wc_j(f)\right|\le C\left(h_j^{1-13\ve}+h_j^{\delta_1-2n_1\ve}\right)\|f\|_{H^{\max\{20,3n_1\}}(S\M)}.
\end{equation}
On the other hand, if (ii) holds with $\delta_2$ and $n_2\in\N$, then (i) holds with $n_1=\max\{n_2,20\}+2$ and any $\delta_1$ with
\[0<\delta_1<\frac{\beta}{2n_1+2+\beta}\cdot\min\left\{\delta_2,\frac12\right\}.\]
\end{thm}

\begin{rmk}
If $\mu_1>\frac14$, then $\beta=\frac12$, so (ii) implies (i), with $n_1=\max\{n_2,20\}+2$, for any
\[0<\delta_1<\frac{1}{4n_1+5}\cdot\min\left\{\delta_2,\frac12\right\}.\]
\end{rmk}

\begin{rmk}
As in Theorem \ref{thm:equiv}, we state the equivalence in Theorem \ref{thm:rate} at the single Ehrenfest time $-T_E$. However, its consequence goes beyond $-T_E$. Indeed, by \eqref{eq:Egorovbeyond}, the stronger the rate in (i) is, the longer the time beyond the Ehrenfest time is in Egorov's theorem: the rate $h_j^{\delta_1}$ gives Egorov's theorem with a polynomial remainder up to $-(1+\ve)T_E$ for every $\ve<\min\{\frac{1}{13},\frac{\delta_1}{2n_1}\}$.
\end{rmk}

\begin{rmk}[Outline of the proof]
The proof that (i) implies (ii) uses the renormalized Egorov estimate in Theorem \ref{thm:renorm} and the invariance of $d\Vol_L$ under the Schr\"odinger group $e^{-\frac{i\sigma}{4}U_+^2}$.

The proof that (ii) implies (i) has three steps. First, the renormalized Egorov estimate and the commutator of $Q$ with $e^{-\frac{i}{4}U_+^2}$ show that
\[\Wc_j\left(U_+^2f\right)=O_f\left(h_j^{\delta_2}+h_j\right),\]
that is, $\Wc_j$ is approximately annihilated by $U_+^2$. Second, $\Wc_j$ is close to a probability measure on $S\M$ up to $O(h_j^{\frac12}|\log h_j|)$; see Lemma \ref{lemma:positive}. Third, a probability measure on $S\M$ which is approximately annihilated by $U_+^2$ is close to $d\Vol_L$ with a polynomial rate, by the effective equidistribution of the stable horocycle flow; see Lemma \ref{lemma:horoRate}.
\end{rmk}

The remainder of the introduction discusses the role of the Schr\"odinger group along the stable horocycles in Section \ref{sec:Schr}, the QUE conjecture in Section \ref{sec:QUE}, quantizations in Section \ref{sec:quant}, and the background on Egorov's theorem in Section \ref{sec:Egorov}.

\subsection{The Schr\"odinger group along the stable horocycles}\label{sec:Schr}
In this subsection, we explain the role of the Schr\"odinger group in \eqref{eq:Schr}. Backward in time, it absorbs the weight $e^{-2s}$ in \eqref{eq:exact}. Precisely, let
\[\sigma(t)=h_j\left(e^{2|t|}-1\right)\quad\text{for }t\le0.\]
Then for $f\in C^\infty(S\M)$ and $t\le0$,
\begin{equation}\label{eq:exactSchr}
\Wc_j\left(f\circ a_t\right)-\Wc_j\left(e^{-\frac{i\sigma(t)}{4}U_+^2}f\right)=\frac{ih_j}{2}\int_0^t\Wc_j\left(\left(Q^2e^{-\frac{i(\sigma(t)-\sigma(s))}{4}U_+^2}f\right)\circ a_s\right)\,ds;
\end{equation}
see Theorem \ref{thm:exactSchr}. Compared with \eqref{eq:exact}, the term $e^{-2s}(U_+^2f)\circ a_s$ disappears. Indeed, the Schr\"odinger group is the evolution generated by $\frac{ih_j}{2}e^{-2s}U_+^2$, with the time $s$ replaced by $\sigma(s)$. Applying \eqref{eq:exactSchr} twice in Section \ref{sec:backward}, we show that $\Wc_j(f\circ a_t)$ and $\Wc_j(e^{-\frac{i\sigma(t)}{4}U_+^2}f)$ differ by $O_f(h_j(1+\sigma(t))^{13})$ (Theorem \ref{thm:renorm}). This difference is small up to $-(1+\ve)T_E$ with $0\le\ve<\frac{1}{13}$. That is, the matrix elements $\Wc_j$ approximately intertwine the geodesic flow and the Schr\"odinger group along the stable horocycles.

Anantharaman-Zelditch \cite{AZ1, AZ2} also obtained exact relations between the quantum evolution and the geodesic flow in Zelditch's quantization, by a different approach. In \cite{AZ1}, they introduced the Patterson-Sullivan distributions associated with the eigenfunctions, which are invariant under the geodesic flow, and showed that they are asymptotically equal to the matrix elements $\Wc_j$, with exact relations between the two families. In \cite{AZ2}, they constructed an explicit intertwining operator on certain Hilbert spaces of symbols, relating the geodesic flow to the conjugation by the Schr\"odinger group of $\Delta_\Gamma$ on $L^2(\M)$. The resulting modified quantization satisfies an exact intertwining relation; see also Gansemer-Hassell \cite{GH} for the traces of the operators in Zelditch's quantization whose symbols are evolved by the geodesic flow.

Our Schr\"odinger group is different from the one in \cite{AZ2}. In \cite{AZ2}, it is generated by the Laplacian on $\M$ and gives the quantum evolution, while in this paper, it is generated by $U_+^2$ on $S\M$ and acts on the symbols. Moreover, in \cite{AZ1, AZ2}, the quantization, or the matrix elements, are modified by an operator independent of time. In contrast, we keep Zelditch's quantization and the matrix elements $\Wc_j$, and modify the symbols by $e^{-\frac{i\sigma(t)}{4}U_+^2}$, which depends on time. This modification is negligible before the Ehrenfest time. At and beyond the Ehrenfest time, it gives exactly the limit of the remainder in Egorov's theorem in \eqref{eq:defect}. This is how Egorov's theorem in Zelditch's quantization detects QUE.

\subsection{The QUE conjecture}\label{sec:QUE}
The quantum ergodicity theorem of Shnirelman-Zelditch-Colin de Verdi\`ere \cite{Sh, Ze2, CdV} asserts that ``almost all'' eigenfunctions tend to be equidistributed in the phase space. The theorem uses the ergodicity of the geodesic flow on the unit sphere bundle with respect to the Liouville measure, and applies to compact manifolds with ergodic geodesic flow. In our setting, it states:

\begin{thmn}[Quantum ergodicity]
There is a full density subsequence $\{\vp_j\}_{j\in J}$ which tends to be equidistributed in the phase space.
\end{thmn}

Here, the density of $J\subset\N$ is
\[\lim_{N\to\infty}\frac{\Card\left(J\cap[0,N]\right)}{N+1},\quad\text{if the limit exists}.\]
When the density is $1$, we call $\{\vp_j\}_{j\in J}$ a full density subsequence. Hence, quantum ergodicity states that ``almost all'' eigenfunctions tend to be equidistributed in the phase space, while the QUE conjecture states that all of them do.

While there is an abundance of invariant measures for the geodesic flow, the QUE conjecture asserts that only the Liouville measure arises as a semiclassical measure. Partial results are known about the semiclassical measures. Anantharaman \cite{A1} and Anantharaman-Nonnenmacher \cite{AN} proved that every semiclassical measure has positive entropy, at least $\frac12$ in the constant curvature case. More recently, Dyatlov-Jin \cite{DyJi} proved that every semiclassical measure has full support, using the fractal uncertainty principle of Bourgain-Dyatlov \cite{BD}; see also Dyatlov-Jin-Nonnenmacher \cite{DJN} for the result on surfaces with variable negative curvature.

In other chaotic systems, QUE may fail. For the quantum cat maps, it fails for the half-scarred eigenstates in \cite{FNDB}, while Kurlberg-Rudnick \cite{KR1} proved that it holds for the joint eigenstates of the quantum cat map and its Hecke operators. The billiard flow on the Bunimovich stadium is ergodic by Bunimovich \cite{Bun}. However, Hassell \cite{Has} proved that for almost every aspect ratio of the rectangle in $[1,2]$, some sequence of the Dirichlet eigenfunctions on the Bunimovich stadium induces a semiclassical measure which is not the Liouville measure. We compare these two models with our results in Appendix \ref{sec:non}.

On compact arithmetic surfaces, there are Hecke operators associated with the symmetry of the structure, and one can find a joint eigenbasis of the Laplacian and these Hecke operators. For such Hecke eigenfunctions, QUE has been verified by Lindenstrauss \cite{L2}; see also Soundararajan \cite{So} for the noncompact modular surface. Therefore, by Theorem \ref{thm:equiv}, Egorov's theorem for a Hecke eigenbasis holds at the Ehrenfest time, and beyond it, i.e., up to $-T_E-T_j$ for some $T_j\to\infty$ depending on the symbol.

A QUE rate of $h_j^{\frac12-\eta}$ for every $\eta>0$ is expected in the arithmetic Hecke setting; see Luo-Sarnak \cite{LS} and Sarnak \cite{Sa2}. The corresponding rate for arbitrary Laplacian eigenbases on compact hyperbolic surfaces remains unclear. For Hecke-Maass eigenforms and symbols independent of the fiber variable, this rate follows from the generalized Lindel\"of hypothesis, by the triple product formulas of Watson \cite{Wa} and Ichino \cite{I}, and a subconvexity bound for the triple product $L$-functions would give a polynomial rate.

To relate this conjectured rate of QUE to Theorem \ref{thm:rate}, suppose that (i) holds with $\delta_1=\frac12-\eta$ for every $\eta>0$, with a Sobolev order $n_1$ independent of $\eta$. Then Egorov's theorem holds at $-T_E$ with a remainder of order $O(h_j^{\frac12-\eta})$ for every $\eta>0$, and up to $-(1+\ve)T_E$ with a polynomial remainder for every
\[0<\ve<\min\left\{\frac{1}{13},\frac{1}{4n_1}\right\}.\]
If the Sobolev order depends on $\eta$, the admissible time range must instead be expressed in terms of $n_1(\eta)$. This direction does not use the effective equidistribution of the stable horocycle flow, so it does not depend on $\mu_1$, the first non-zero eigenvalue of $\Delta_\Gamma$. Hence, even the conjectured rate of QUE only gives Egorov's theorem to $-(1+\ve)T_E$ in our approach. We expect that $-2T_E=-|\log h_j|$ is the natural limit, as in the anisotropic calculus of Dyatlov-Zahl \cite{DZ} and Dyatlov-Jin \cite{DyJi}.

Conversely, the loss from $\delta_2$ to $\delta_1$ in Theorem \ref{thm:rate} comes from the positivity in Lemma \ref{lemma:positive} and the effective equidistribution of the stable horocycle flow in Lemma \ref{lemma:horoRate}. This approach via Egorov's theorem is unlikely to prove QUE with the conjectured rate of $O(h_j^{\frac12-\eta})$, and therefore we do not attempt to optimize these parameters in Theorem \ref{thm:rate}.

We refer to the surveys by Anantharaman \cite{A2, A3}, Sarnak \cite{Sa1, Sa2}, and Zelditch \cite{Ze6, Ze7} for the development related to the QUE conjecture.

\subsection{Quantizations}\label{sec:quant}
There are different frameworks of ``quantization'', which assign semiclassical pseudo-differential operators to symbols. For standard quantizations on a fixed compact energy region, the operators associated with the same symbol independent of $h$ differ by $o_f(1)$ in $L^2(\M)$ operator norm as $h\to0$; see Zworski \cite[Sections 4.1 and 14.2]{Zw} for a systematic treatment. Therefore, the QUE conjecture is independent of the choice of quantization. However, the discussion in this paper on Egorov's theorem is highly dependent on the quantization of Zelditch \cite{Ze1, Ze2, Ze3}, which we now explain.

Zelditch's quantization is associated with a Lagrangian foliation of the phase space. Precisely, for $b\in\Sb=\partial\D$ and $r>0$, let
\[L_{b,r}=\left\{(z,\xi)\in T\D\setminus\{0\}:b_+\left(z,\frac{\xi}{|\xi|_z}\right)=b,|\xi|_z=r\right\},\]
where $b_+$ is the forward endpoint of the geodesic on $\Sb$; see Figure \ref{fig:b}. Using the Riemannian metric to identify $T\D$ with $T^*\D$, each $L_{b,r}$ is a Lagrangian submanifold, namely the graph of the differential of $z\mapsto r\langle z,b\rangle$, where $\langle z,b\rangle$ is the horocyclic distance. The family $\{L_{b,r}\}$ is invariant under $\Gamma$, so it descends to a Lagrangian foliation of $T\M\setminus\{0\}$, whose leaves are the weak stable manifolds of the geodesic flow; see Figure \ref{fig:flows}. The hyperbolic waves
\[w_{b,\frac rh}(z)=e^{(\frac{ir}{h}+\frac12)\langle z,b\rangle}\]
are the Lagrangian states associated with $L_{b,r}$, and $\Op_h(f)$ acts on them by multiplication by the restriction of $f$ to $L_{b,r}$; see Sections \ref{sec:efn} and \ref{sec:SDO}. Consequently, only the derivatives of $f$ along the leaves, i.e., $Qf$ and $U_+f$, appear in the commutator formula in Theorem \ref{thm:commuD} and hence in the remainder in \eqref{eq:exact}. In this sense, Zelditch's quantization has the anisotropic structure along a Lagrangian foliation that is exploited in the calculus of Dyatlov-Zahl \cite{DZ} and Dyatlov-Jin \cite{DyJi}. In these works \cite{DyJi, DZ}, the symbols may oscillate in the direction transversal to the leaves, but vary slowly along the leaves.

However, in our theorems, the symbols $f\circ a_t$ in Egorov's theorem are not necessarily from the anisotropic class associated with this Lagrangian foliation. By Theorem \ref{thm:QU},
\[U_+\left(f\circ a_t\right)=e^{-t}\left(U_+f\right)\circ a_t\quad\text{and}\quad U_-\left(f\circ a_t\right)=e^t\left(U_-f\right)\circ a_t.\]
Forward in time, i.e., $t\ge0$, the derivatives of $f\circ a_t$ along the leaves do not grow, while those in the transversal direction $U_-$ may grow like $e^t$. Thus, forward propagation has the derivative structure associated with the foliation $\{L_{b,r}\}$. The anisotropic symbol classes of \cite{DyJi, DZ} allow transverse oscillations at scales $h^\rho$ with $\rho<1$, giving the corresponding propagation time $\rho|\log h|$. In contrast, (i) of Theorem \ref{thm:Egorov} reaches $\rho T_E$ for every fixed $\rho>0$, including $\rho\ge2$, where $f\circ a_t$ oscillates transversally to the leaves at the scale $e^{-t}\le h_j$, by iterating the exact identity \eqref{eq:exact} for diagonal matrix elements. We remark that Egorov's theorem in \cite{DyJi, DZ} concerns operator norms, whereas our theorem concerns only diagonal matrix elements of eigenfunctions.

Backward in time, i.e., $t\le0$, $f\circ a_t$ can develop oscillations at the scale $e^{-|t|}$ along the stable horocycles, which lie in the leaves $L_{b,r}$. Thus, backward propagation has the derivative structure associated with the other Lagrangian foliation, given by the backward endpoint $b_-$, i.e., by the weak unstable manifolds. In particular, around $t=-T_E$, our quantization is adapted to one foliation while the symbol $f\circ a_t$ is adapted to the other foliation, and oscillates at the scale $h^{\frac12}$. At this scale, the term $e^{-2s}(U_+^2f)\circ a_s$ in \eqref{eq:exact} is of order $h^{-1}$ near $s=-T_E$. So its contribution is of order $1$. So at and beyond the Ehrenfest time, Egorov's theorem for eigenfunctions does not follow from the symbol calculus of \cite{DyJi, DZ}, and it requires a cancellation in the matrix elements of the eigenfunctions. By Theorem \ref{thm:equiv}, this cancellation holds if and only if the eigenfunctions tend to be equidistributed in the phase space. Roughly speaking, our results backward in time correspond to the setting of \cite{DyJi, DZ} in which the quantization is chosen in the anisotropic calculus associated with one Lagrangian foliation, while the symbols are chosen in the anisotropic class associated with the transversal one. Forward in time, the quantization and the symbols are adapted to the same foliation.

If $b_+$ is replaced by the backward endpoint $b_-$ in Zelditch's quantization, then the roles of the two foliations, and the roles of the two time directions, are interchanged. In particular, (i) of Theorem \ref{thm:Egorov} becomes Egorov's theorem backward in time up to $-\rho T_E$ for every $\rho>0$. So beyond the Ehrenfest time, Egorov's theorem depends on the quantization. In this paper, we use Zelditch's quantization through $b_+$ throughout.

\subsection{Egorov's theorem}\label{sec:Egorov}
For standard semiclassical quantizations on $\M$ with symbols localized in a fixed compact energy region, long-time Egorov estimates hold uniformly on logarithmic time intervals strictly within the Ehrenfest range; see Bambusi-Graffi-Paul \cite{BGP}, Bouzouina-Robert \cite{BR}, and Zworski \cite[Sections 11.4 and 14.2]{Zw}.

In the anisotropic calculus of Dyatlov-Zahl \cite{DZ} and Dyatlov-Jin \cite{DyJi} associated with a Lagrangian foliation, Egorov's theorem reaches $2(1-\ve)T_E$ for the symbols adapted to the foliation. Since Zelditch's quantization and the symbols $f\circ a_t$ with $t\le0$ are associated with different foliations (Section \ref{sec:quant}), our diagonal matrix-element estimates cannot be directly compared with the operator-norm estimates in the anisotropic calculus of \cite{DyJi, DZ} beyond the Ehrenfest time $T_E$.

Egorov's theorem at logarithmic times is a key ingredient in the entropy bounds of Anantharaman \cite{A1} and Anantharaman-Nonnenmacher \cite{AN}, in the full support theorems of Dyatlov-Jin \cite{DyJi} and Dyatlov-Jin-Nonnenmacher \cite{DJN}, in the rate of quantum ergodicity by Zelditch \cite{Ze4} and Schubert \cite{Schu}, and in the small scale quantum ergodicity by Han \cite{Han1} and Hezari-Rivi\`ere \cite{HR}.

\subsection*{Organization of the paper}
In Section \ref{sec:pre}, we recall the geometry of $\D$ and $\M$ and the Laplacian eigenfunctions; in Section \ref{sec:SDO}, we recall Zelditch's semiclassical pseudo-differential operators and prove Egorov's theorem for eigenfunctions on $\M$ with an explicit remainder (Theorem \ref{thm:EefnM}); in Section \ref{sec:rep}, we recall the representation theory of $\PSU(1,1)$ on $L^2(S\M)$, and prove a uniform bound and the positivity of the matrix elements $\Wc_j$; in Section \ref{sec:forward}, we prove (i) of Theorem \ref{thm:Egorov}; in Section \ref{sec:backward}, we give the precise definition and establish the properties of the Schr\"odinger group along the stable horocycles, prove Theorem \ref{thm:exactSchr}, and prove (ii) of Theorem \ref{thm:Egorov}; in Section \ref{sec:beyond}, we prove Theorem \ref{thm:equiv}; in Section \ref{sec:rate}, we recall the effective equidistribution of the stable horocycle flow and prove Theorem \ref{thm:rate}; in Appendix \ref{sec:non}, we compare our results with the quantum cat maps and the Bunimovich stadium, in which QUE fails.

\section{Preliminaries}\label{sec:pre}
In this section, we introduce the basic terminology and facts used throughout the paper. We recall the geometry of the Poincar\'e disk model $\D$, the relevant Lie theory of $\SU(1,1)$ and $\PSU(1,1)$, the geodesic and horocycle flows on $S\D$, and the corresponding structures on the compact quotient $\M$. Throughout the paper, $\N=\{0,1,2,\dots\}$.

We fix a compact connected hyperbolic surface $\M=\Gamma\bs\D$, where $\Gamma\subset\PSU(1,1)$ is a torsion-free cocompact Fuchsian group. The geometric structures, flows, and differential operators on $\M$ and $S\M$ are obtained by descending the corresponding $\Gamma$-invariant objects from $\D$ and $S\D$.

Our setup follows standard treatments of hyperbolic geometry, Lie theory, and Laplacian eigenfunctions, such as Helgason \cite[Chapter 4]{He} and Zelditch \cite{Ze1, Ze2, Ze3}. This preliminary section is not intended to be self-contained; we refer to these works for details. We adopt different normalizations, however, which must be borne in mind when comparing formulas. In particular:
\begin{itemize}
\item The Gaussian curvature of $\D$ and $\M$ is normalized to be $-1$, which fixes the metric and Laplacian used below. This differs from the convention of \cite[Chapter 4]{He} and \cite{Ze1, Ze2, Ze3} by a constant rescaling of the metric and a corresponding rescaling of the Laplacian.
\item The Laplacians $\Delta$ on $\D$ and $\Delta_\Gamma$ on $\M$ are taken to be non-negative. This is opposite in sign to the convention of \cite[Chapter 4]{He} and \cite{Ze1, Ze2, Ze3}.
\end{itemize}

\subsection{The Poincar\'e disk $\D$}\label{sec:D}
The Poincar\'e disk $\D$ is the unit disk $\{z\in\C:|z|<1\}$ equipped with the Riemannian metric
\[ds^2=\frac{4|dz|^2}{\left(1-|z|^2\right)^2}=\frac{4\left(dx^2+dy^2\right)}{\left(1-|z|^2\right)^2},\quad\text{where }z=x+iy.\]
Hence, $\D$ has constant Gaussian curvature $-1$. Denote by $d_\D(z,z')$ the Riemannian distance between $z$ and $z'$ in $\D$, and by $o$ the origin in $\D$.

\begin{prop}[The isometry group of $\D$]
Let
\[\SU(1,1)=\left\{\begin{pmatrix}
\alpha & \beta\\
\ol\beta & \ol\alpha\end{pmatrix}:\alpha,\beta\in\C\text{ and }|\alpha|^2-|\beta|^2=1\right\}\]
be the indefinite special unitary group of the Hermitian form of signature $(1,1)$. The elements of $\SU(1,1)$ act on $\D$ by linear fractional transformations:
\[\gamma=
\begin{pmatrix}\alpha & \beta\\
\ol\beta & \ol\alpha\end{pmatrix}\in\SU(1,1),\quad z\mapsto\gamma z=\frac{\alpha z+\beta}{\ol\beta z+\ol\alpha}\quad\text{for }z\in\D.\]
The kernel of this action is $\{\pm\Id\}$. Hence, the action descends to an effective action of
\[\PSU(1,1)=\SU(1,1)/\{\pm\Id\}.\]
Moreover, $\PSU(1,1)$ is the group of orientation-preserving isometries of $\D$, and its action on $\D$ is transitive.
\end{prop}

The Riemannian area form on $\D$ is
\[d\Area=\frac{4\,dxdy}{\left(1-|z|^2\right)^2},\]
and it is invariant under $\PSU(1,1)$.

\subsection{The Lie theory of $\SU(1,1)$ and $\PSU(1,1)$}\label{sec:Lie}
We record the Lie-theoretic facts used below. The group $\SU(1,1)$ is a $3$-dimensional Lie group, and its Lie algebra $\suf(1,1)$ has a basis $\{W,V,Q\}$, where
\[W=\frac12\begin{pmatrix}i & 0\\
0 & -i\end{pmatrix},\quad V=\frac12\begin{pmatrix}0 & 1\\
1 & 0\end{pmatrix},\quad Q=\frac12\begin{pmatrix}0 & i\\
-i & 0\end{pmatrix}.\]
These elements satisfy the commutation relations
\[[W,V]=Q,\quad[Q,V]=W,\quad[Q,W]=V.\]
Define
\[U_+=V+W=\frac12\begin{pmatrix}i & 1\\
1 & -i\end{pmatrix}\quad\text{and}\quad U_-=V-W=\frac12\begin{pmatrix}-i & 1\\
1 & i\end{pmatrix}.\]
Then
\begin{equation}\label{eq:WUQ}
\left[W,U_\pm\right]=Q,\quad\left[V,U_\pm\right]=\mp Q,\quad\left[Q,U_\pm\right]=\pm U_\pm,\quad\left[U_+,U_-\right]=2Q.
\end{equation}
The elements $W$, $U_+$, $Q$ are the infinitesimal generators of the one-parameter subgroups $K$, $N$, $A$ of $\SU(1,1)$, respectively:
\begin{enumerate}[(1)]
\item
\[K=\left\{k_\theta=e^{\theta W}=\begin{pmatrix}e^{\frac{i\theta}{2}} & 0\\
0 & e^{-\frac{i\theta}{2}}\end{pmatrix}:\theta\in\R\right\}\cong\SO(2).\]
\item
\[N=\left\{n_s=e^{sU_+}=\frac12\begin{pmatrix}2+is & s\\
s & 2-is\end{pmatrix}:s\in\R\right\}.\]
\item
\[A=\left\{a_t=e^{tQ}=\begin{pmatrix}\cosh\frac t2 & i\sinh\frac t2\\
-i\sinh\frac t2 & \cosh\frac t2\end{pmatrix}:t\in\R\right\}.\]
\end{enumerate}
By $[Q,U_+]=U_+$ in \eqref{eq:WUQ}, we have $\Ad_{a_t}U_+=e^tU_+$, and therefore
\[n_sa_t=a_tn_{se^{-t}}\quad\text{for }s,t\in\R.\]
Thus, for fixed $s\in\R$, the trajectories $\gamma n_sa_t$ and $\gamma a_t$ approach each other as $t\to+\infty$. We therefore call $N$ the stable horocyclic subgroup.

We denote the image in $\PSU(1,1)$ of an element or a subgroup of $\SU(1,1)$ by the same symbol. In $\PSU(1,1)$, the map $\theta\mapsto k_\theta$ has period $2\pi$, and $k_\theta$ acts on $\D$ as the rotation $z\mapsto e^{i\theta}z$. Thus, $K$ is the stabilizer of the origin $o$ in $\PSU(1,1)$.

Write $\Sb=\partial\D$, the boundary of the Poincar\'e disk $\D$ in $\C$. The action of $\PSU(1,1)$ by linear fractional transformations extends to $\Sb$.

\begin{defn}[Lie derivatives]
Let $Y\in\suf(1,1)$. The Lie derivative $\Lc_Y$ on $\SU(1,1)$ is defined by
\[\left(\Lc_Yu\right)(\gamma)=\partial_t\left[u\left(\gamma e^{tY}\right)\right]_{t=0}\quad\text{for }u\in C^\infty(\SU(1,1))\text{ and }\gamma\in\SU(1,1).\]
\end{defn}

The Lie derivative $\Lc_Y$ is a left-invariant vector field on $\SU(1,1)$. For notational simplicity, we write $\Lc_Y=Y$.

\begin{defn}[Casimir operator]
The Casimir operator $\Omega$ on $C^\infty(\SU(1,1))$ is defined by
\[\Omega=-Q^2-V^2+W^2.\]
\end{defn}

Using the commutation relations, one checks that $\Omega$ commutes with $Q,V,W$. See Knapp \cite[Chapter III]{Kn} for further background.

\subsection{The sphere bundle of the Poincar\'e disk}\label{sec:bundle}
Denote by $S\D$ the sphere bundle of $\D$, i.e., the unit tangent bundle. We identify each tangent space $T_z\D$ with $\C$. An element $\gamma\in\SU(1,1)$ acts on $S\D$ by the differential of its action on $\D$. Explicitly,
\[\gamma=\begin{pmatrix}\alpha & \beta\\\ol\beta & \ol\alpha\end{pmatrix}:(z,\xi)\mapsto(\gamma z,\gamma\xi)=\left(\frac{\alpha z+\beta}{\ol\beta z+\ol\alpha},\frac{\xi}{\left(\ol\beta z+\ol\alpha\right)^2}\right)\quad\text{for }(z,\xi)\in S\D.\]
Both elements of $\{\pm\Id\}$ act trivially, so this action factors through $\PSU(1,1)$. For each $(z,\xi)\in S\D$, there is a unique element $\gamma\in\PSU(1,1)$ such that
\[(z,\xi)=\left(\gamma o,\gamma\frac i2\right).\]
Thus, $S\D$ is identified with $\PSU(1,1)$. Under this identification, the geodesic flow and the stable horocycle flow on $S\D$ are defined by right translations in $\PSU(1,1)$.

\begin{defn}[The geodesic flow and the horocycle flows]
Let $\gamma\in S\D\simeq\PSU(1,1)$.
\begin{itemize}
\item The geodesic flow $a_t$ on $S\D$ is defined by $\gamma\mapsto\gamma e^{tQ}$, $t\in\R$.
\item The stable horocycle flow on $S\D$ is defined by $\gamma\mapsto\gamma e^{sU_+}$, $s\in\R$.
\item The unstable horocycle flow on $S\D$ is defined by $\gamma\mapsto\gamma e^{sU_-}$, $s\in\R$.
\end{itemize}
\end{defn}

We use the same symbol for an element of $\PSU(1,1)$ and for the diffeomorphism of $S\D$ that it induces by right translation. Thus, $a_t=e^{tQ}$ denotes both the element of $A$ defined in Section \ref{sec:Lie} and the geodesic flow on $S\D$, and likewise for $n_s=e^{sU_+}$.

\begin{figure}[htb]
\centering
\begin{tikzpicture}
\draw[thick] (0,0) circle (4);
\fill (0,0) circle (.05) node[anchor=north west] {$o$};

\draw[thick] (0,2) circle (2);
\draw[thick] (0,-2) circle (2);

\draw[->] (0,0) -- (0,4);
\draw[->] (2,2) arc (270:180:2);
\draw[->] (-2,2) arc (270:360:2);
\draw[->] (1.414,0.586) arc (225:180:4.828);
\draw[->] (-1.414,0.586) arc (315:360:4.828);

\draw[->,dashed] (0,-4) -- (0,0);
\draw[->,dashed] (0,-4) arc (180:90:2);
\draw[->,dashed] (0,-4) arc (0:90:2);
\draw[->,dashed] (0,-4) arc (180:135:4.828);
\draw[->,dashed] (0,-4) arc (0:45:4.828);

\draw (2.4,3) node[below] {$e^{sU_+}$};
\draw (2.4,-3) node[above] {$e^{sU_-}$};
\draw (0,.7) node[left] {$e^{tQ}$};
\draw (0,.7) node[right] {$t>0$};
\draw (0,-.7) node[left] {$e^{tQ}$};
\draw (0,-.7) node[right] {$t<0$};

\draw[->] (0,0) -- (0,2) node[right] {$\frac i2$};
\draw[->] (2,2) -- (1,2);
\draw[->] (-2,2) -- (-1,2);
\draw[->] (2,-2) -- (3,-2);
\draw[->] (-2,-2) -- (-3,-2);

\draw[->] (1.414,0.586) -- (0.207,1.793);
\draw[->] (-1.414,0.586) -- (-0.207,1.793);
\draw[->] (1.414,-0.586) -- (2.621,0.621);
\draw[->] (-1.414,-0.586) -- (-2.621,0.621);
\end{tikzpicture}
\caption{The stable and unstable horocycle flows at $(o,\frac i2)\in S\D$}\label{fig:flows}
\end{figure}
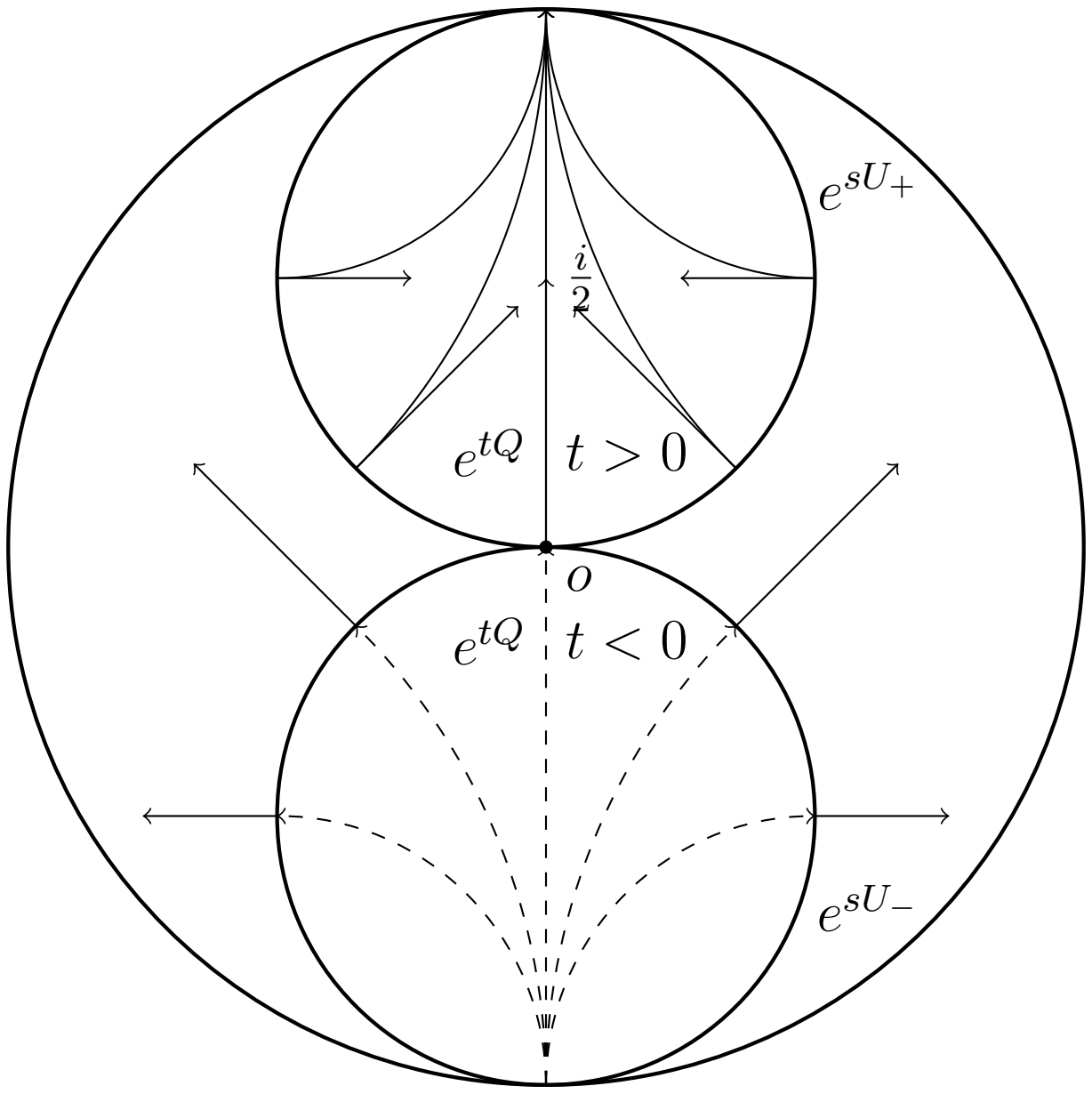

In Figure \ref{fig:flows}, the geodesic trajectory starting from $(o,\frac i2)\in S\D$ at $t=0$ is shown as a solid line for $t>0$ and as a dashed line for $t<0$. The stable and unstable horocycles through $o$ are labeled by $e^{sU_+}$ and $e^{sU_-}$, respectively. Along the stable horocycle flow $e^{sU_+}$, the geodesics shown as solid curves converge to the geodesic starting from $(o,\frac i2)$ as $t\to\infty$; along the unstable horocycle flow $e^{sU_-}$, the geodesics shown as dashed curves converge to it as $t\to-\infty$.

\begin{figure}[htb]
\centering
\begin{tikzpicture}
\draw[thick] (0,0) circle (4);
\draw[thick,->] (-3.464,2) arc (240:300:6.928);
\draw[thick,->] (0,1.072) -- (4/3,1.072) node[below] {$\xi$};
\draw (0,1.072) node[above] {$(z,\xi)$};

\fill (0,0) circle (.05) node[below] {$o$};
\fill (0,1.072) circle (.05) node[below] {$z$};
\fill (3.464,2) circle (.05) node[anchor=south west] {$b_+(z,\xi)$};
\fill (-3.464,2) circle (.05) node[anchor=south east] {$b_-(z,\xi)$};
\end{tikzpicture}
\caption{The forward and backward endpoints $b_\pm(z,\xi)\in\Sb$}\label{fig:b}
\end{figure}
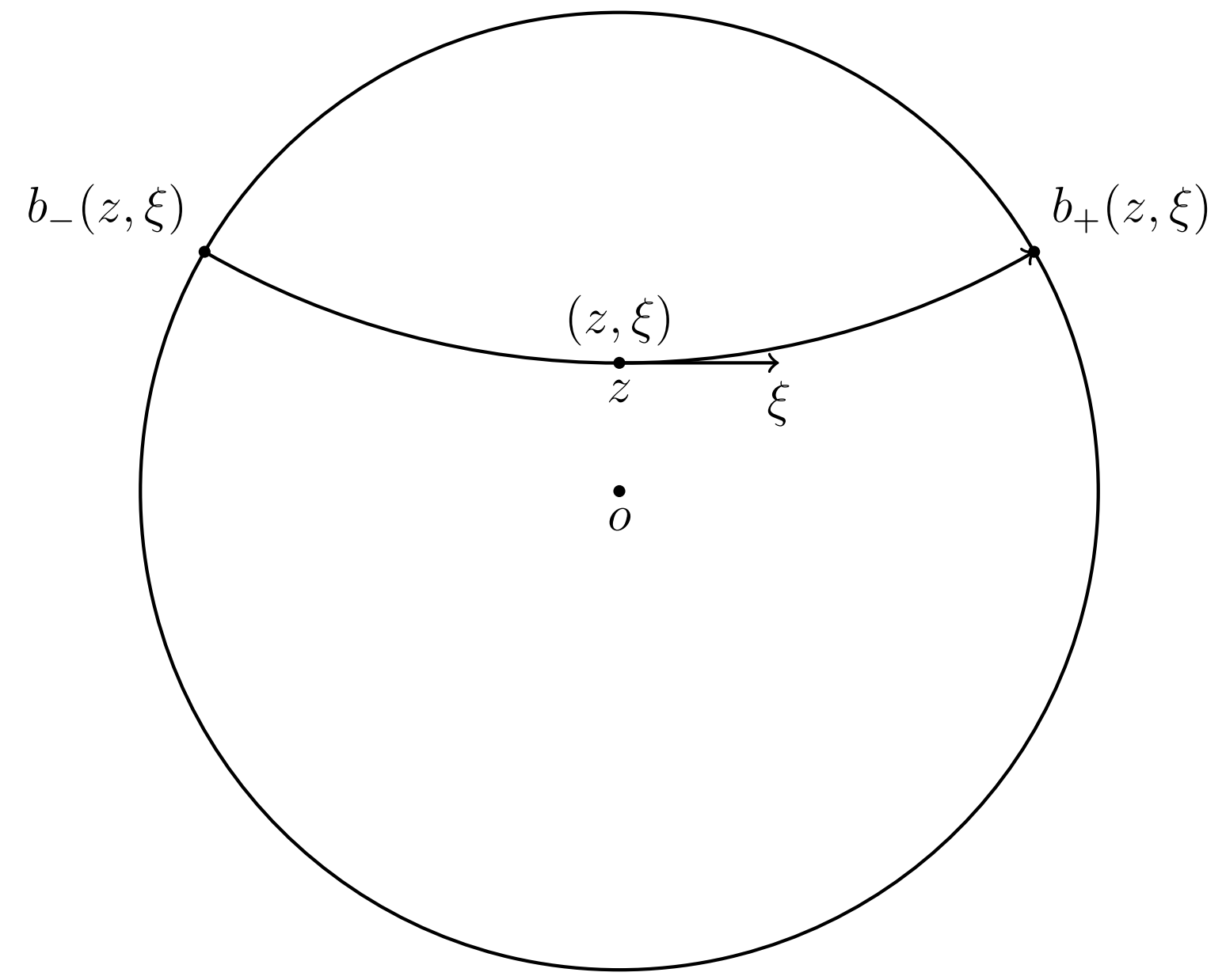

Define the natural projection $\Pi:S\D\to\D$ by
\[\Pi(z,\xi)=z\quad\text{for }(z,\xi)\in S\D.\]
For each $(z,\xi)\in S\D$, define, with the limits taken in $\ol\D\subset\C$,
\[b_+(z,\xi)=\lim_{t\to\infty}\Pi\left(a_t(z,\xi)\right)\in\Sb\quad\text{and}\quad b_-(z,\xi)=\lim_{t\to-\infty}\Pi\left(a_t(z,\xi)\right)\in\Sb;\]
see Figure \ref{fig:b}.

Denote by $T\D$ the tangent bundle of $\D$. Each vector field $X\in\{Q,U_+,U_-\}$ on $S\D$ has a unique extension to $T\D\setminus\{0\}$ which commutes with the fiber dilations and restricts to $X$ on $S\D$. We denote this extension by $X$ as well. Thus, the geodesic flow $a_t=e^{tQ}$ and the stable horocycle flow $e^{sU_+}$ extend to $T\D\setminus\{0\}$. The extended flow $a_t$ preserves $|\xi|_z$ and moves the base point along the geodesic determined by $\xi$ at unit speed, independently of $|\xi|_z$.

\begin{defn}\leavevmode
\begin{itemize}
\item For $(z,\xi)\in S\D$, define
\[B_+(z,\xi)=\left(z,b_+(z,\xi)\right)\in\D\times\Sb.\]
Thus, the map $B_+$ identifies $S\D$ with $\D\times\Sb$.
\item For $(z,\xi)\in T\D\setminus\{0\}$, define
\[B_+(z,\xi)=\left(z,b_+\left(z,\frac{\xi}{|\xi|_z}\right),|\xi|_z\right)\in\D\times\Sb\times\R^+,\]
where $\R^+=(0,\infty)$. Thus, the map $B_+$ identifies $T\D\setminus\{0\}$ with $\D\times\Sb\times\R^+$.
\end{itemize}
\end{defn}

Define the vector field $D=\xi\cdot\nabla_\xi$, which is the generator of dilations in the fiber direction of $T\D\setminus\{0\}$. Then $\{Q,U_+,U_-,D\}$ forms a frame on $T\D\setminus\{0\}$. Since $Q,U_+,U_-$ are homogeneous of degree $0$, they commute with $D$.

\begin{thm}\label{thm:QU}
Let $f\in C^\infty(T\D\setminus\{0\})$. For $t\in\R$ and $(\alpha_1,\alpha_2,\alpha_3,\alpha_4)\in\N^4$,
\[Q^{\alpha_1}U_+^{\alpha_2}U_-^{\alpha_3}D^{\alpha_4}\left(f\circ a_t\right)=e^{\left(\alpha_3-\alpha_2\right)t}\left(Q^{\alpha_1}U_+^{\alpha_2}U_-^{\alpha_3}D^{\alpha_4}f\right)\circ a_t.\]
The same holds for the products of $Q,U_+,U_-,D$ in any order, with $\alpha_2$ and $\alpha_3$ the numbers of the factors $U_+$ and $U_-$, respectively.
\end{thm}

\begin{proof}
For $Y\in\suf(1,1)$, we have
\[Y\left(f\circ a_t\right)=\left(\left(\Ad_{a_{-t}}Y\right)f\right)\circ a_t.\]
By \eqref{eq:WUQ}, $\Ad_{a_{-t}}Q=Q$ and $\Ad_{a_{-t}}U_\pm=e^{\mp t}U_\pm$. Moreover, $D$ commutes with $a_t$. The claim follows by induction on the number of factors.
\end{proof}

\begin{figure}[htb]
\centering
\begin{tikzpicture}
\draw[thick] (0,0) circle (4);
\fill (0,0) circle (.05) node[anchor=north west] {$o$};

\draw[thick] (0,3) circle (1);
\draw[thick] (0,-1) circle (3);

\fill (0,2) circle (.05) node[anchor=north west] {$z$};

\fill (0,4) circle (.05) node[above] {$b_+(z,\xi)$};
\fill (0,-4) circle (.05) node[below] {$b_-(z,\xi)$};
\draw[->,thick] (0,0) -- (0,4);
\draw[-,thick] (0,-4) -- (0,0);
\draw[-,thick] (0,0) -- (0,2) node[midway,left] {$(z,\xi)$};
\draw[->,thick] (0,2) -- (0,3) node[midway,left] {$\xi$};
\end{tikzpicture}
\caption{The horocycles through $z$ tangent to $\Sb$ at $b_\pm(z,\xi)$}\label{fig:horo}
\end{figure}
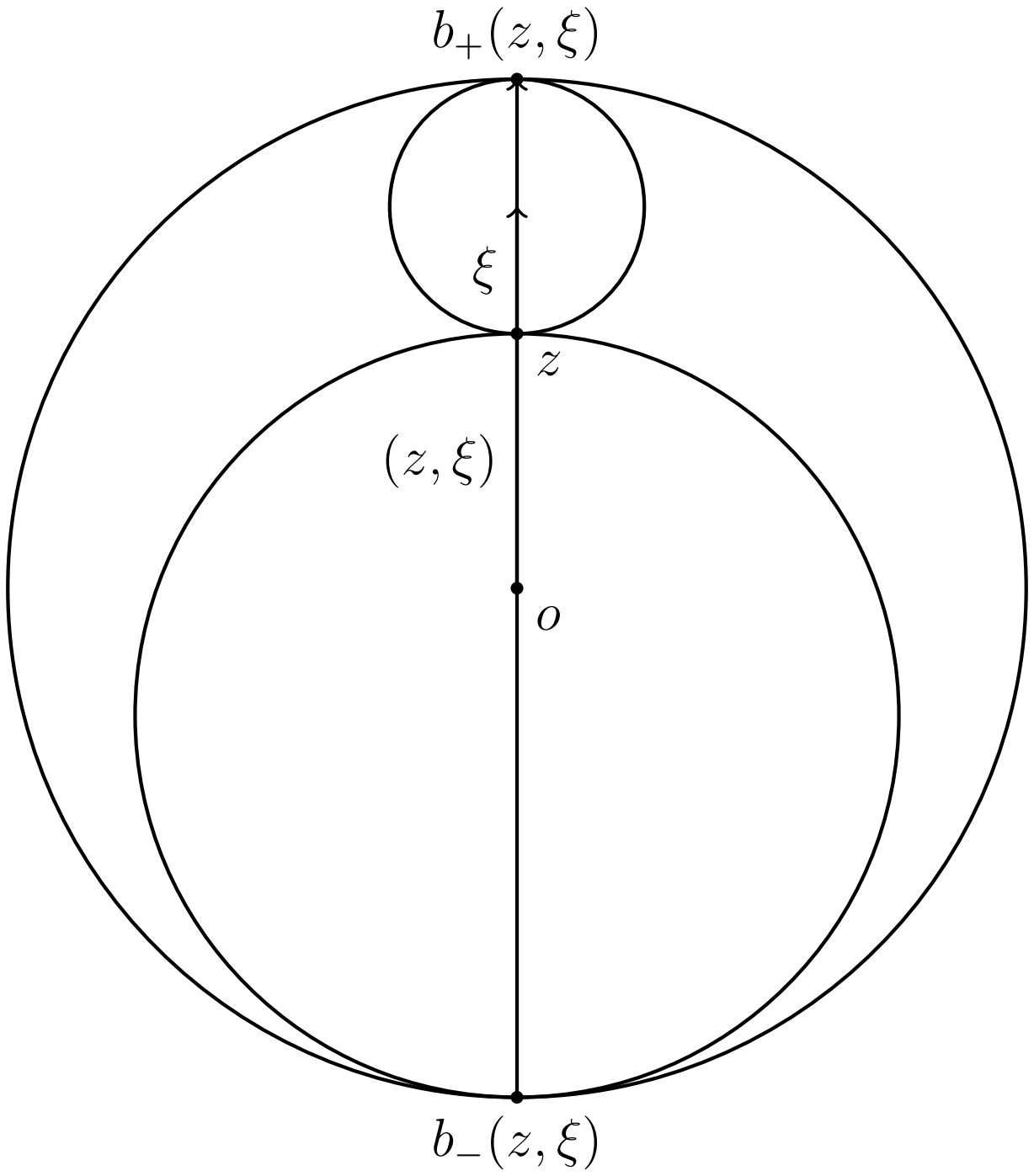

For $z\in\D$ and $b\in\Sb$, define $\langle z,b\rangle$ to be the signed distance from $o$ to the horocycle passing through $z$ and tangent to $\Sb$ at $b$. The sign is chosen so that
\[\langle z,b\rangle\begin{cases}
>0,&\text{if $o$ lies in the exterior of the horocycle},\\
=0,&\text{if $o$ lies on the horocycle},\\<0,&\text{if $o$ lies in the interior of the horocycle}.
\end{cases}\]
In Figure \ref{fig:horo}, $o$ lies in the exterior of the horocycle through $z$ tangent to $\Sb$ at $b_+(z,\xi)$, and in the interior of the horocycle through $z$ tangent to $\Sb$ at $b_-(z,\xi)$. So $\langle z,b_+(z,\xi)\rangle>0$ and $\langle z,b_-(z,\xi)\rangle<0$.

\begin{prop}\label{prop:horo}
Let $\gamma\in\PSU(1,1)$ and $(z,b)\in\D\times\Sb$.
\begin{enumerate}[(i)]
\item
\[\langle\gamma z,\gamma b\rangle=\langle z,b\rangle+\langle\gamma o,\gamma b\rangle.\]
\item
\[\frac{d(\gamma b)}{db}=e^{-\langle\gamma o,\gamma b\rangle},\]
where the derivative is taken with respect to the positively oriented angular coordinate on $\Sb$.
\item The function $(z,b)\mapsto\langle z,b\rangle$ is smooth on $\D\times\Sb$, and satisfies
\[Q\langle z,b\rangle=1\quad\text{and}\quad U_+\langle z,b\rangle=0.\]
\end{enumerate}
\end{prop}

\begin{proof}
See Helgason \cite[Chapter 4]{He} for proofs of (i) and (ii). To show (iii), identify $S\D$ with $\D\times\Sb$ using $B_+$. Under the geodesic flow $a_t=e^{tQ}$, the base point moves toward $b$ along a geodesic at unit speed, so the signed distance from $o$ to the horocycle determined by $a_t(z,b)$ increases by $t$. Hence, $Q\langle z,b\rangle=1$. Under the stable horocycle flow $e^{sU_+}$, the base point moves along this horocycle, which therefore remains fixed. Hence, $U_+\langle z,b\rangle=0$; see Figure \ref{fig:flows}.
\end{proof}

\begin{prop}[The Liouville measure]\label{prop:Liouville}
Identify $(z,\xi)\in S\D$ with $(z,b)=B_+(z,\xi)\in\D\times\Sb$. Then the Liouville measure on $S\D$ is given by
\[d\Vol=\frac{1}{2\pi}e^{\langle z,b\rangle}\,dzdb,\]
where $dz=d\Area(z)$ is the Riemannian area form on $\D$, and $db$ is the Lebesgue measure on $\Sb$. Moreover, $d\Vol$ is invariant under $\PSU(1,1)$, and under the geodesic flow $e^{tQ}$ and the stable horocycle flow $e^{sU_+}$.
\end{prop}

See Zelditch \cite[Proposition 1.3]{Ze1} and Anantharaman-Zelditch \cite[Lemma 3.1]{AZ2} for a proof. The invariance under $\PSU(1,1)$ follows from (i) and (ii) of Proposition \ref{prop:horo}, and the invariance of $d\Area$. Under the identification $S\D\simeq\PSU(1,1)$, $d\Vol$ is a Haar measure on $\PSU(1,1)$. Since $\PSU(1,1)$ is unimodular, $d\Vol$ is also invariant under right translations.

The fiber of $S\D$ over $z=\gamma o$, $\gamma\in\PSU(1,1)$, is $\{\gamma k_\theta:\theta\in[0,2\pi)\}$, and $b_+(\gamma k_\theta)=\gamma b'$ with $b'=ie^{i\theta}$. By (ii) of Proposition \ref{prop:horo}, the substitution $b=\gamma b'$ gives
\[d\theta=db'=e^{\langle\gamma o,b\rangle}db=e^{\langle z,b\rangle}db.\]
Hence, for $u\in C^\infty(S\D)$,
\begin{equation}\label{eq:fiber}
\frac{1}{2\pi}\int_0^{2\pi}u\left(\gamma k_\theta\right)\,d\theta=\frac{1}{2\pi}\int_\Sb u(z,b)e^{\langle z,b\rangle}\,db,
\end{equation}
where the left-hand side is the average of $u$ over the fiber over $z$.

\subsection{Laplacian eigenfunctions on the Poincar\'e disk}\label{sec:efn}
\begin{defn}[Laplacian on $\D$]
Denote by $\Delta$ the non-negative Laplacian on $\D$. In the coordinates $z=x+iy\in\D$,
\[\Delta=-\frac{\left(1-x^2-y^2\right)^2}{4}\left(\partial_x^2+\partial_y^2\right).\]
\end{defn}

A function $\vp\in C^\infty(\D)$ is a Laplacian eigenfunction on $\D$ with eigenvalue $\mu\in\C$ if $\Delta\vp=\mu\vp$. We identify $C^\infty(\PSU(1,1))$ with the space of $\{\pm\Id\}$-invariant functions in $C^\infty(\SU(1,1))$. The Lie derivatives $Q,V,W,U_\pm$ and the Casimir operator $\Omega$ of Section \ref{sec:Lie} act on this space, hence on $C^\infty(S\D)$ under the identification $S\D\simeq\PSU(1,1)$.

\begin{thm}\label{thm:LonD}
Under the identification $S\D\simeq\PSU(1,1)$ and $(z,b)=B_+(z,\xi)$ of Section \ref{sec:bundle}, let $u\in C^\infty(S\D)$. Then
\[\Delta u(z,b)=\left(-Q^2-U_+^2+Q\right)u(z,b),\]
where $\Delta$ is the Laplacian acting in the variable $z$ on $\D$ with $b\in\Sb$ held fixed. Moreover, if $u$ is independent of the $b$-variable, then $u$ defines a function on $\D$, and
\[\Omega u(z)=\Delta u(z)\quad\text{for }z\in\D.\]
\end{thm}

\begin{proof}
See Zelditch \cite[Proposition 2.1]{Ze3} for the proof of the first identity. We provide a short proof here for completeness. Identify $\D$ with the upper half-plane
\[\Hb=\{x+iy\in\C:y>0\}\quad\text{ equipped with the metric }ds^2=\frac{dx^2+dy^2}{y^2}.\]
By the invariance under the orientation-preserving isometry group, it suffices to consider the boundary point $b=\infty$. Under the identification of the isometry group with $\mathrm{PSL}(2,\R)$, write
\[n_x=\begin{pmatrix}
1 & x\\
0 & 1
\end{pmatrix}\quad\text{and}\quad
a_t=\begin{pmatrix}
e^{\frac t2} & 0\\
0& e^{-\frac t2}
\end{pmatrix}.\]
Here, $n_x$ and $a_t$ denote elements of $\mathrm{PSL}(2,\R)$. They correspond to the elements $n_s$ and $a_t$ of Section \ref{sec:Lie} under this identification, so right multiplication by $a_t$ and $n_s$ gives the geodesic flow and the stable horocycle flow in Section \ref{sec:bundle}. The element $n_xa_{\log y}$ corresponds to the unit tangent vector at $x+iy$ pointing toward $\infty$. Right multiplication gives
\[n_xa_{\log y}a_t=n_xa_{\log y+t}\quad\text{and}\quad n_xa_{\log y}n_s=n_{x+ys}a_{\log y}.\]
Consequently, on functions $u(x+iy,\infty)$, the vector fields $Q$ and $U_+$ act as
\[Q=y\partial_y\quad\text{and}\quad U_+=y\partial_x.\]
Therefore,
\begin{align*}
\left(-Q^2-U_+^2+Q\right)u
&=-\left(y\partial_y\right)^2u-\left(y\partial_x\right)^2u+y\partial_yu\\
&=-y^2\left(\partial_x^2+\partial_y^2\right)u\\
&=\Delta u.
\end{align*}
By equivariance, the identity holds for every $b\in\Sb$.

For the second identity, suppose that $u$ is independent of the $b$-variable. Since $W$ generates the fiber rotation $K$, this is equivalent to $Wu=0$. As $V=U_+-W$, we have
\begin{align*}
\Omega&=-Q^2-V^2+W^2\\
&=-Q^2-(U_+-W)^2+W^2\\
&=-Q^2-U_+^2+U_+W+WU_+\\
&=-Q^2-U_+^2+Q+2U_+W.
\end{align*}
Thus, $\Omega u=(-Q^2-U_+^2+Q)u=\Delta u$.
\end{proof}

We follow Helgason \cite[Chapter 4]{He} to characterize Laplacian eigenfunctions in terms of hyperbolic waves on $\D$.

\begin{defn}[Hyperbolic waves]
Let $b\in\Sb$ and $\lambda\in\C$. The hyperbolic wave $w_{b,\lambda}$ on $\D$ is defined by
\[w_{b,\lambda}(z)=e^{\left(i\lambda+\frac12\right)\langle z,b\rangle}\quad\text{for }z\in\D.\]
\end{defn}

\begin{prop}\label{prop:wave}
Let $b\in\Sb$ and $\lambda\in\C$. Then $w_{b,\lambda}$ is a Laplacian eigenfunction on $\D$ with eigenvalue
\[\mu=\lambda^2+\frac14.\]
\end{prop}

\begin{proof}
By (iii) of Proposition \ref{prop:horo}, we have $U_+\langle z,b\rangle=0$ and $Q\langle z,b\rangle=1$. Hence, by Theorem \ref{thm:LonD},
\begin{align*}
\Delta w_{b,\lambda}&=\left(-Q^2-U_+^2+Q\right)w_{b,\lambda}\\
&=\left[-\left(i\lambda+\frac12\right)^2+\left(i\lambda+\frac12\right)\right]w_{b,\lambda}\\
&=\left(\lambda^2+\frac14\right)w_{b,\lambda}.
\end{align*}
\end{proof}

For $\tau\in\Dc'(\Sb)$ and $v\in C^\infty(\Sb)$, we write the pairing between $\tau$ and $v$ formally as
\[\langle v,\tau\rangle=\int_\Sb v(b)\tau(b)\,db.\]
Let $z\in\D$ and $\lambda\in\C$. The function $b\mapsto w_{b,\lambda}(z)$ belongs to $C^\infty(\Sb)$. Hence, the integral in the following theorem is understood as such a pairing.

\begin{thm}\label{thm:efnD}
Let $\lambda\in\C$ with $i\lambda+\frac12\notin\{0,-1,-2,\dots\}$. Then the map
\[\tau\mapsto\vp_{\tau,\lambda}(z)=\int_\Sb e^{\left(i\lambda+\frac12\right)\langle z,b\rangle}\tau(b)\,db\]
is a bijection between $\Dc'(\Sb)$ and the space of Laplacian eigenfunctions $\vp$ on $\D$ with eigenvalue $\lambda^2+\frac14$ that satisfy, for some constant $C>0$,
\[\left|\vp(z)\right|\le Ce^{Cd_\D(o,z)}\quad\text{for }z\in\D.\]
\end{thm}

See \cite[Chapter 4]{He}.

\subsection{Compact hyperbolic surfaces}\label{sec:compact}
Recall the compact connected hyperbolic surface $\M=\Gamma\bs\D$ fixed in Section \ref{sec:pre}. Each $u\in C^\infty(\M)$ is identified with an element $\tilde u$ of $C^\infty_\Gamma(\D)$, the space of smooth $\Gamma$-invariant functions on $\D$, i.e.,
\[\tilde u(\gamma z)=\tilde u(z)\quad\text{for }\gamma\in\Gamma,z\in\D.\]
The Riemannian area form $d\Area$ descends to $\M$, and
\[\int_\M u\,d\Area=\int_\F \tilde u\,d\Area\quad\text{for }u\in C^\infty(\M),\]
where $\F\subset\D$ is a fundamental domain for $\Gamma$.

Let $T\M$ denote the tangent bundle of $\M$ and $S\M$ denote the sphere bundle of $\M$, realized as the quotient
\[S\M=\Gamma\bs S\D=\Gamma\bs\PSU(1,1).\]
Each $u\in C^\infty(T\M\setminus\{0\})$ is identified with an element $\tilde u$ of $C^\infty_\Gamma(T\D\setminus\{0\})$, the space of smooth $\Gamma$-invariant functions on $T\D\setminus\{0\}$, i.e.,
\[\tilde u(\gamma z,\gamma\xi)=\tilde u(z,\xi)\quad\text{for }\gamma\in\Gamma,(z,\xi)\in T\D\setminus\{0\},\]
and the same identification, restricted to $S\D$, identifies $C^\infty(S\M)$ with $C^\infty_\Gamma(S\D)$. We write $\Dc'_\Gamma(S\D)$ for the space of $\Gamma$-invariant distributions on $S\D$, which is identified with $\Dc'(S\M)$ in the same way.

The Liouville measure $d\Vol$ on $S\M$ is obtained from $d\Vol$ on $S\D$, and
\[\int_{S\M}u\,d\Vol=\int_{S\F}\tilde u\,d\Vol\quad\text{for }u\in C^\infty(S\M),\]
where $S\F$ is the sphere bundle over $\F$, and where we use the same notation for a measure on $S\M$ and for its lift to $S\D$. We normalize by
\[d\Area_L=\frac{d\Area}{\Area(\M)}\quad\text{and}\quad d\Vol_L=\frac{d\Vol}{\Vol(S\M)},\]
so that $d\Area_L$ and $d\Vol_L$ are probability measures. By \eqref{eq:fiber}, $\Vol(S\M)=\Area(\M)$, and the pushforward of $d\Vol_L$ under the natural projection $S\M\to\M$ is $d\Area_L$.

For $1\le p<\infty$, the spaces $L^p(\M)$ and $L^p(S\M)$ are the completions of $C^\infty(\M)$ and of $C^\infty(S\M)$ with respect to the norms $\|\cdot\|_{L^p(\M)}$ and $\|\cdot\|_{L^p(S\M)}$, taken with respect to $d\Area_L$ and $d\Vol_L$, respectively.

The right translations by elements of $\PSU(1,1)$ on $S\D$ and $T\D\setminus\{0\}$ commute with the action of $\Gamma$, and therefore descend to $S\M$ and $T\M\setminus\{0\}$. In particular, the geodesic flow $a_t$, the stable horocycle flow $e^{sU_+}$, and the vector fields $Q,U_+,U_-,W$ descend to $S\M$, and $d\Vol_L$ is invariant under these flows.

\subsection{Laplacian eigenbasis on the compact hyperbolic surface}\label{sec:ebasis}
Using the identification $C^\infty(\M)=C^\infty_\Gamma(\D)$, the Laplacian $\Delta$ on $\D$ descends to the Laplacian $\Delta_\Gamma$ on $\M=\Gamma\bs\D$. Precisely, if $u\in C^\infty(\M)$ is identified with $\tilde u\in C^\infty_\Gamma(\D)$, then $\Delta\tilde u\in C^\infty_\Gamma(\D)$ is the lift of $\Delta_\Gamma u\in C^\infty(\M)$.

Thus, if $\vp$ is a Laplacian eigenfunction on $\M$ with eigenvalue $\mu$, then its lift $\tilde\vp$ satisfies
\[\begin{cases}\Delta\tilde\vp(z)=\mu\tilde\vp(z)&\text{for }z\in\D,\\\tilde\vp(\gamma z)=\tilde\vp(z)&\text{for }\gamma\in\Gamma,z\in\D.\end{cases}\]
It is standard that $\Delta_\Gamma$ is self-adjoint and has a discrete spectrum of non-negative eigenvalues; see, e.g., Zworski \cite[Section 14.3]{Zw}. Moreover, there is an orthonormal basis $\{\vp_j\}_{j=0}^\infty$ of Laplacian eigenfunctions in $L^2(\M)$, with eigenvalues
\[0=\mu_0<\mu_1\le\mu_2\le\mu_3\le\cdots\to\infty\quad\text{as }j\to\infty,\]
counted with multiplicity. Throughout the paper, we fix such a Laplacian eigenbasis $\{\vp_j\}_{j=0}^\infty$ on $\M$.

Finally, we characterize the Laplacian eigenfunctions $\vp_j$ on $\M$ using Helgason's theory in Section \ref{sec:efn}. For $\mu_j\ge\frac14$, we fix the spectral parameter
\begin{equation}\label{eq:lambda}
\lambda_j=\sqrt{\mu_j-\frac14}\ge0,
\end{equation}
so that $\mu_j=\lambda_j^2+\frac14$ and $i\lambda_j+\frac12\notin\{0,-1,-2,\dots\}$. If $\mu_j>\frac54$, then $\lambda_j>1$, and $h_j=\lambda_j^{-1}\in(0,1)$ is the semiclassical parameter in Section \ref{sec:intro}. Since $\M$ is compact, $\tilde\vp_j\in C^\infty_\Gamma(\D)$ is bounded, hence of uniform exponential growth in the sense of Theorem \ref{thm:efnD}. Therefore, there is a unique $\tau_j\in\Dc'(\Sb)$ such that
\begin{equation}\label{eq:tildevp}
\tilde\vp_j(z)=\int_\Sb w_{b,\lambda_j}(z)\tau_j(b)\,db=\int_\Sb e^{\left(i\lambda_j+\frac12\right)\langle z,b\rangle}\tau_j(b)\,db.
\end{equation}
Due to the invariance of $\tilde\vp_j$ under $\Gamma$, Zelditch \cite[Propositions 2.1 and 2.2]{Ze2} proved the following result.

\begin{thm}\label{thm:Phi}
Let $\mu_j\ge\frac14$. The distribution
\[\tilde\Phi_j(z,b)=e^{\left(i\lambda_j-\frac12\right)\langle z,b\rangle}\tau_j(b)\]
belongs to $\Dc'_\Gamma(S\D)$ and therefore descends to a distribution $\Phi_j\in\Dc'(S\M)$. Moreover,
\[Q\Phi_j=\left(i\lambda_j-\frac12\right)\Phi_j\quad\text{and}\quad U_+\Phi_j=0.\]
\end{thm}

\section{Semiclassical pseudo-differential operators}\label{sec:SDO}
In this section, we recall Zelditch's construction of the semiclassical pseudo-differential operators on the Poincar\'e disk and on compact hyperbolic surfaces, and prove Egorov's theorem for eigenfunctions on $\M$ with an explicit remainder.

Throughout, let $h\in(0,1)$ be the semiclassical parameter. If $f=f(z,\xi;h)$ is defined on $(T\D\setminus\{0\})\times(0,1)$, we write $f\in C^\infty(T\D\setminus\{0\})$ to mean that $f(\cdot,\cdot;h)$ is smooth on $T\D\setminus\{0\}$ for each $h\in(0,1)$, and we suppress the dependence on $h$ from the notation.

If $f=f(z,\xi)\in C^\infty(T\D\setminus\{0\})$ satisfies appropriate growth conditions, then $f$ is called a symbol, and one assigns to it a semiclassical pseudo-differential operator $\Op_h(f)$. This construction is global in the sense that it uses the global coordinates $(z,\xi)\in T\D\setminus\{0\}$ and the Helgason-Fourier transform on $\D$. It follows Zelditch's hyperbolic pseudo-differential calculus \cite{Ze1, Ze2}, which is based on Helgason's non-Euclidean Fourier analysis on $\D$ \cite[Chapter 4]{He}. See also Anantharaman-Zelditch \cite{AZ1, AZ2}.

In this global theory of pseudo-differential operators, Zelditch proved that if a symbol $f$ is invariant under $\Gamma$, i.e., $f(\gamma z,\gamma\xi)=f(z,\xi)$ for $(z,\xi)\in T\D\setminus\{0\}$ and $\gamma\in\Gamma$, then $\Op_h(f)$ at least \textit{formally} maps $\Gamma$-invariant functions to $\Gamma$-invariant functions, and hence induces an operator on $\M$, denoted by $\Op_h^\Gamma(f)$. In fact, $\Op_{h_j}^\Gamma(f)\vp_j$ for eigenfunctions $\vp_j$ can be defined via the boundary values $\tau_j$ of $\vp_j$ when the eigenvalues $\mu_j>\frac54$, so that $h_j\in(0,1)$.

\subsection{Semiclassical pseudo-differential operators on $\D$}
We define symbol classes for semiclassical pseudo-differential operators on $\D$. These are functions on $T\D\setminus\{0\}$ satisfying growth conditions with respect to the frame $\{Q,U_+,U_-,D\}$ of Section \ref{sec:bundle}.

\begin{defn}[Symbols]
Let $m\in\R$ and $\rho\ge0$. We say that $f\in C^\infty(T\D\setminus\{0\})$ belongs to $S^m_\rho(T\D\setminus\{0\})$, the symbol class of order $m$ with semiclassical derivative growth $\rho$, if for each multi-index $\alpha=(\alpha_1,\alpha_2,\alpha_3,\alpha_4)\in\N^4$ there is $C=C(\alpha)>0$ such that
\[\left|Q^{\alpha_1}U_+^{\alpha_2}U_-^{\alpha_3}D^{\alpha_4}f(z,\xi;h)\right|\le Ch^{-\rho|\alpha|}(1+r)^m\quad\text{for }(z,\xi)\in T\D\setminus\{0\}\text{ and }h\in(0,1),\]
where $r=|\xi|_z$. The topology of $S^m_\rho(T\D\setminus\{0\})$ is defined by the seminorms $\|f\|_{\rho,m,\alpha}$, where $\|f\|_{\rho,m,\alpha}$ is the infimum of the admissible constants $C(\alpha)$.

We omit the superscript $m$ in $S_\rho^m(T\D\setminus\{0\})$ when $m=0$, and the subscript $\rho$ when $\rho=0$. Moreover, define
\[S_\rho^\infty(T\D\setminus\{0\})=\bigcup_{m\in\R}S_\rho^m(T\D\setminus\{0\}).\]
\end{defn}

At logarithmic times, the geodesic flow maps symbols in $S^m$ to symbols with positive semiclassical derivative growth $\rho$.

\begin{prop}\label{prop:seminorm}
Let $f\in S^m(T\D\setminus\{0\})$, $\rho\ge0$, and $|t|\le\rho|\log h|$. Then $f\circ a_t\in S^m_\rho(T\D\setminus\{0\})$, and for each $\alpha\in\N^4$,
\[\left\|f\circ a_t\right\|_{\rho,m,\alpha}\le\|f\|_{0,m,\alpha}.\]
\end{prop}

\begin{proof}
By Theorem \ref{thm:QU} and the invariance of $r=|\xi|_z$ under $a_t$,
\begin{align*}
\left|Q^{\alpha_1}U_+^{\alpha_2}U_-^{\alpha_3}D^{\alpha_4}\left(f\circ a_t\right)(z,\xi;h)\right|
&=\left|e^{(\alpha_3-\alpha_2)t}\left(Q^{\alpha_1}U_+^{\alpha_2}U_-^{\alpha_3}D^{\alpha_4}f\right)\circ a_t(z,\xi;h)\right|\\
&\le e^{(\alpha_2+\alpha_3)|t|}\|f\|_{0,m,\alpha}(1+r)^m\\
&\le h^{-\rho|\alpha|}\|f\|_{0,m,\alpha}(1+r)^m,
\end{align*}
since $|\alpha_3-\alpha_2|\le\alpha_2+\alpha_3\le|\alpha|$ and $e^{|t|}\le h^{-\rho}$.
\end{proof}

Recall from Section \ref{sec:bundle} the identification $B_+:T\D\setminus\{0\}\to\D\times\Sb\times\R^+$. In the coordinates $(z,b,r)=B_+(z,\xi)$, the fiber radius is $r=|\xi|_z$, and we identify it with the semiclassical frequency $r=h\lambda$ of the hyperbolic wave $w_{b,\lambda}$. This identification is the point at which the quantization is made.

\begin{defn}[Semiclassical pseudo-differential operators on $\D$]
Let $f\in S_\rho^\infty(T\D\setminus\{0\})$. We prescribe the action of the semiclassical pseudo-differential operator $\Op_h(f)$ on the hyperbolic waves by
\[\Op_h(f)w_{b,\frac rh}(z)=f\circ B_+^{-1}(z,b,r;h)w_{b,\frac rh}(z)\quad\text{for }z\in\D,b\in\Sb,r\in\R^+.\]
We say that $f$ is the full symbol of $\Op_h(f)$. Here, $f\circ B_+^{-1}(z,b,r;h)$ means $f(B_+^{-1}(z,b,r);h)$; the map $B_+$ itself is independent of $h$.
\end{defn}

The same prescription extends by linearity to symbols of the form $h^Ng$, where $N\in\R$ and $g\in S_\rho^\infty(T\D\setminus\{0\})$, by setting $\Op_h(h^Ng)=h^N\Op_h(g)$.

By the Helgason-Fourier inversion formula \cite[Chapter 4]{He}, every $u\in C_0^\infty(\D)$ is a superposition of the hyperbolic waves $w_{b,\frac rh}$, $b\in\Sb$, $r\in\R^+$. Hence, the above prescription defines $\Op_h(f):C_0^\infty(\D)\to C^\infty(\D)$; see Zelditch \cite{Ze1}.

We next extend the action of $\Op_h(f)$ to the Laplacian eigenfunctions on $\D$, which are not compactly supported. We do this by applying the prescribed action on the hyperbolic waves to their boundary representations. Recall from Theorem \ref{thm:efnD} that a bounded eigenfunction $\vp_{\tau,\lambda}$ with eigenvalue $\lambda^2+\frac14$, $\lambda\in\R^+$, can be written as
\[\vp_{\tau,\lambda}(z)=\int_\Sb e^{\left(i\lambda+\frac12\right)\langle z,b\rangle}\tau(b)\,db,\]
where $\tau\in\Dc'(\Sb)$ is the boundary value of $\vp_{\tau,\lambda}$. For $f\in S_\rho^\infty(T\D\setminus\{0\})$, we define
\begin{equation}\label{eq:SDOefnD}
\Op_h(f)\vp_{\tau,\lambda}(z)=\int_\Sb f\circ B^{-1}_+(z,b,h\lambda;h)e^{\left(i\lambda+\frac12\right)\langle z,b\rangle}\tau(b)\,db\quad\text{for }z\in\D.
\end{equation}
Since $b\mapsto f\circ B^{-1}_+(z,b,h\lambda;h)e^{(i\lambda+\frac12)\langle z,b\rangle}$ is smooth on $\Sb$ and depends smoothly on $z$, the function $\Op_h(f)\vp_{\tau,\lambda}$ is smooth on $\D$.

\begin{thm}\label{thm:commuD}
Let $f\in S_\rho^\infty(T\D\setminus\{0\})$. Then for $b\in\Sb$ and $r\in\R^+$,
\[\left[h^2\left(\Delta-\frac14\right),\Op_h(f)\right]w_{b,\frac rh}=\Op_h(q)w_{b,\frac rh},\quad\text{where }q=-h^2\left(Q^2+U_+^2\right)f-2ihrQf,\]
with $r=|\xi|_z$ at $(z,\xi)\in T\D\setminus\{0\}$. The quantization of $q$ is understood in the extended sense above.
\end{thm}

\begin{proof}
Identify $(z,\xi)\in T\D\setminus\{0\}$ with $(z,b,r)=B_+(z,\xi)\in\D\times\Sb\times\R^+$, and write $w=w_{b,\frac rh}$. Throughout the proof, we write $f$ for $f\circ B_+^{-1}(z,b,r;h)$, viewed as a function of $z$ at fixed $(b,r)$.

By Theorem \ref{thm:LonD},
\[\Delta=-Q^2-U_+^2+Q.\]
Since the coordinates $b=b_+$ and $r=|\xi|_z$ are preserved by the flows $e^{tQ}$ and $e^{sU_+}$, the vector fields $Q$ and $U_+$ are tangent to the leaf $\{(z,b,r):z\in\D\}$ of fixed $(b,r)$, which we identify with $\D$. By (iii) of Proposition \ref{prop:horo},
\[Qw=\left(i\frac rh+\frac12\right)w\quad\text{and}\quad U_+w=0.\]
Then, by Proposition \ref{prop:wave},
\begin{align*}
\Delta(\Op_h(f)w)(z)
&=\Delta(f\cdot w)(z)\\
&=w(z)\Delta f-2Qf\cdot Qw-2U_+f\cdot U_+w+f\cdot\Delta w\\
&=\left(\Delta f-2\left(i\frac rh+\frac12\right)Qf+\left(\frac{r^2}{h^2}+\frac14\right)f\right)w(z),
\end{align*}
and
\[\Op_h(f)(\Delta w)(z)=\left(\frac{r^2}{h^2}+\frac14\right)fw(z).\]
Hence,
\begin{align*}
\left[h^2\left(\Delta-\frac14\right),\Op_h(f)\right]w(z)
&=h^2\left(\Delta f-2\left(i\frac rh+\frac12\right)Qf\right)w(z)\\
&=\left(-h^2(Q^2+U_+^2)f-2ihrQf\right)w(z).
\end{align*}
\end{proof}

\subsection{Semiclassical pseudo-differential operators on $\M$}\label{sec:SDOonM}
Recall that $\M=\Gamma\bs\D$ is the fixed compact hyperbolic surface. Let $m\in\R$ and $\rho\ge0$. We say that $f\in C^\infty(T\M\setminus\{0\})$ belongs to the symbol class $S^m_\rho(T\M\setminus\{0\})$ if its $\Gamma$-invariant lift $\tilde f\in C^\infty_\Gamma(T\D\setminus\{0\})$ belongs to $S^m_\rho(T\D\setminus\{0\})$, and we define $S^\infty_\rho(T\M\setminus\{0\})$ accordingly. Since $Q,U_+,U_-,D$ commute with the action of $\Gamma$, the seminorms $\|\cdot\|_{\rho,m,\alpha}$ in the definition of $S^m_\rho(T\D\setminus\{0\})$ descend to $S^m_\rho(T\M\setminus\{0\})$, and Proposition \ref{prop:seminorm} holds on $T\M\setminus\{0\}$.

In this paper, we define the action of $\Op_{h_j}^\Gamma(f)$ on each eigenfunction $\vp_j$ with $\mu_j>\frac54$, rather than construct $\Op_h^\Gamma(f)$ as an operator on all of $C^\infty(\M)$. This is sufficient for our purpose. To consider semiclassical pseudo-differential operators acting on $C^\infty(\M)$, one would usually need to impose additional conditions on the symbol $f$, such as that $f(z,\xi)\in S^\infty_\rho(T\M\setminus\{0\})$ can be extended holomorphically to a strip of appropriate width in the $r=|\xi|_z$ variable. See Anantharaman-Zelditch \cite[Section 4.2]{AZ2} for a detailed discussion.

Continue with the Laplacian eigenbasis $\{\vp_j\}_{j=0}^\infty$ on $\M$, and suppose that $\mu_j\ge\frac14$. Recall the spectral parameter $\lambda_j$ in \eqref{eq:lambda}. Then there is a unique $\tau_j\in\Dc'(\Sb)$ such that
\[\tilde\vp_j(z)=\int_\Sb e^{\left(i\lambda_j+\frac12\right)\langle z,b\rangle}\tau_j(b)\,db.\]
By Theorem \ref{thm:Phi},
\[\tilde\Phi_j(z,b)=e^{\left(i\lambda_j-\frac12\right)\langle z,b\rangle}\tau_j(b)\in\Dc'_\Gamma(S\D).\]
Suppose that $\mu_j>\frac54$. Let $f\in S^\infty_\rho(T\M\setminus\{0\})$. By \eqref{eq:SDOefnD},
\begin{align}
\Op_h(\tilde f)\tilde\vp_j(z)
&=\int_\Sb\tilde f\circ B^{-1}_+\left(z,b,h\lambda_j;h\right)e^{\left(i\lambda_j+\frac12\right)\langle z,b\rangle}\tau_j(b)\,db\nonumber\\
&=\int_\Sb\tilde f\circ B^{-1}_+\left(z,b,h\lambda_j;h\right)\tilde\Phi_j(z,b)e^{\langle z,b\rangle}\,db.\label{eq:fvptilde}
\end{align}
We verify that $\Op_h(\tilde f)\tilde\vp_j$ is $\Gamma$-invariant. Let $\gamma\in\Gamma$. By Proposition \ref{prop:horo},
\begin{align*}
\Op_h(\tilde f)\tilde\vp_j(\gamma z)
&=\int_\Sb\tilde f\circ B^{-1}_+\left(\gamma z,b,h\lambda_j;h\right)\tilde\Phi_j(\gamma z,b)e^{\langle\gamma z,b\rangle}\,db\\
&=\int_\Sb\tilde f\circ B^{-1}_+\left(\gamma z,\gamma b',h\lambda_j;h\right)\tilde\Phi_j(\gamma z,\gamma b')e^{\langle\gamma z,\gamma b'\rangle}\,d(\gamma b')\\
&=\int_\Sb\tilde f\circ B^{-1}_+\left(\gamma z,\gamma b',h\lambda_j;h\right)\tilde\Phi_j(\gamma z,\gamma b')e^{\langle z,b'\rangle+\langle\gamma o,\gamma b'\rangle}e^{-\langle\gamma o,\gamma b'\rangle}\,db'\\
&=\int_\Sb\tilde f\circ B^{-1}_+\left(z,b',h\lambda_j;h\right)\tilde\Phi_j(z,b')e^{\langle z,b'\rangle}\,db'\\
&=\Op_h(\tilde f)\tilde\vp_j(z).
\end{align*}
Here, we used the $\Gamma$-invariance of $\tilde f$ and $\tilde\Phi_j$, together with the change of variables $b=\gamma b'$ on $\Sb$.

\begin{thm}\label{thm:SDOefnM}
Let $\mu_j>\frac54$ and $f\in S^\infty_\rho(T\M\setminus\{0\})$. Then $\Op_{h_j}(\tilde f)\tilde\vp_j\in C^\infty_\Gamma(\D)$, and therefore we define it as $\Op_{h_j}^\Gamma(f)\vp_j$. Moreover,
\[\left\langle\Op_{h_j}^\Gamma(f)\vp_j,\vp_j\right\rangle_{L^2(\M)}=2\pi\int_{S\M}f\circ B^{-1}_+\left(z,b,1;h_j\right)\Phi_j(z,b)\ol{\vp_j(z)}\,d\Vol_L.\]
\end{thm}

\begin{proof}
The $\Gamma$-invariance is verified above. Since $h_j\lambda_j=1$, by \eqref{eq:fvptilde} and Proposition \ref{prop:Liouville},
\begin{align*}
\left\langle\Op_{h_j}^\Gamma(f)\vp_j,\vp_j\right\rangle_{L^2(\M)}
&=\frac{1}{\Area(\M)}\int_\F\int_\Sb\tilde f\circ B^{-1}_+\left(z,b,1;h_j\right)\tilde\Phi_j(z,b)\ol{\tilde\vp_j(z)}e^{\langle z,b\rangle}\,dbdz\\
&=\frac{2\pi}{\Area(\M)}\int_{S\F}\tilde f\circ B^{-1}_+\left(z,b,1;h_j\right)\tilde\Phi_j(z,b)\ol{\tilde\vp_j(z)}\,d\Vol,
\end{align*}
in the sense of distributions. Since $\Vol(S\M)=\Area(\M)$, this is the claim.
\end{proof}

In particular, the matrix element $\langle\Op_{h_j}^\Gamma(f)\vp_j,\vp_j\rangle_{L^2(\M)}$ depends only on the restriction of $f$ to $S\M$, i.e., on the values of $f$ at $r=|\xi|_z=1$. This justifies the identification in Section \ref{sec:intro} of a symbol with its restriction to $S\M$; if the symbol depends on $h$, then it is evaluated at $h=h_j$. Then, by Theorem \ref{thm:SDOefnM},
\[\Wc_j(f)=\left\langle\Op_{h_j}^\Gamma(f)\vp_j,\vp_j\right\rangle_{L^2(\M)}=2\pi\int_{S\M}f\,\Phi_j\ol{\vp_j}\,d\Vol_L\quad\text{for }f\in C^\infty(S\M),\]
so $\Wc_j$ is a distribution on $S\M$. Since the flows $a_t$ and $e^{sU_+}$ preserve $r=|\xi|_z$, the restriction to $S\M$ commutes with $a_t$ and $e^{sU_+}$, and with the vector fields $Q,U_+,U_-$. Notice that if $f\in S^0(T\M\setminus\{0\})$ is independent of $h$ and $|t|\le\rho|\log h|$, then $f\circ a_t\in S^0_\rho(T\M\setminus\{0\})$ by Proposition \ref{prop:seminorm}. So $\Wc_j(f\circ a_t)$ in Section \ref{sec:intro} is the matrix element of $\Op_{h_j}^\Gamma(f\circ a_t)$ in Theorem \ref{thm:SDOefnM}, also when $t$ depends on $h_j$. We remark that the definition of $S^m_\rho(T\D\setminus\{0\})$ imposes no upper bound on $\rho$, and for $\rho\ge\frac12$, these classes do not admit a symbol calculus. We do not use any symbol calculus in this paper: we only use the action \eqref{eq:SDOefnD} of $\Op_h(f)$ on eigenfunctions, which is defined for every $f\in S^\infty_\rho(T\D\setminus\{0\})$.

We next establish the commutation formula for semiclassical pseudo-differential operators and the semiclassical Laplacian on $\M$. This follows directly from Theorem \ref{thm:commuD}.

\begin{thm}\label{thm:commuM}
Let $\mu_j>\frac54$ and $f\in S^\infty_\rho(T\M\setminus\{0\})$. Then
\[\left[h_j^2\left(\Delta_\Gamma-\frac14\right),\Op_{h_j}^\Gamma(f)\right]\vp_j=\Op_{h_j}^\Gamma(q)\vp_j,\quad\text{where }q=-h_j^2\left(Q^2+U_+^2\right)f-2ih_jrQf.\]
Here, $r=|\xi|_z$ at $(z,\xi)\in T\M\setminus\{0\}$, and the quantization of $q$ is understood in the extended sense of Theorem \ref{thm:commuD}.
\end{thm}

We prove an important characterization of the distribution $\Phi_j\ol{\vp_j}\in\Dc'(S\M)$. This result is due to Zelditch \cite[Proposition 2.1]{Ze3}; see also Wolpert \cite[Lemma 4.7]{Wo1}.

\begin{thm}\label{thm:Phiode}
Let $\mu_j>\frac54$. Then
\[\left(Q^2+U_+^2-2i\lambda_jQ\right)\Phi_j\ol{\vp_j}=0.\]
\end{thm}

\begin{proof}
Let $f\in S^0(T\M\setminus\{0\})$ be independent of $h$. Since $(\Delta_\Gamma-\frac14)\vp_j=\lambda_j^2\vp_j$ and $\Delta_\Gamma$ is self-adjoint,
\begin{align*}
\left\langle\left(\Delta_\Gamma-\frac14\right)\Op_{h_j}^\Gamma(f)\vp_j,\vp_j\right\rangle_{L^2(\M)}&=\left\langle\Op_{h_j}^\Gamma(f)\vp_j,\left(\Delta_\Gamma-\frac14\right)\vp_j\right\rangle_{L^2(\M)}\\
&=\lambda_j^2\left\langle\Op_{h_j}^\Gamma(f)\vp_j,\vp_j\right\rangle_{L^2(\M)},
\end{align*}
which is also $\langle\Op_{h_j}^\Gamma(f)(\Delta_\Gamma-\frac14)\vp_j,\vp_j\rangle_{L^2(\M)}$. Therefore, by Theorem \ref{thm:commuM}, with $r=1$ on $S\M$,
\begin{align*}
0&=\left\langle\left[h_j^2\left(\Delta_\Gamma-\frac14\right),\Op_{h_j}^\Gamma(f)\right]\vp_j,\vp_j\right\rangle_{L^2(\M)}\\
&=\left\langle\Op_{h_j}^\Gamma(q)\vp_j,\vp_j\right\rangle_{L^2(\M)}\\
&=-h_j^2\Wc_j\left(\left(Q^2+U_+^2+2i\lambda_jQ\right)f\right).
\end{align*}
Since $d\Vol_L$ is invariant under the flows generated by $Q$ and $U_+$, the distributional adjoints of $Q$ and $U_+$ are $-Q$ and $-U_+$, respectively. Since the restriction of $f$ to $S\M$ is an arbitrary function in $C^\infty(S\M)$, this implies the claim.
\end{proof}

We now derive Egorov's theorem for eigenfunctions on $\M$ with an explicit remainder from Theorem \ref{thm:Phiode}.

\begin{thm}[Egorov's theorem for eigenfunctions on $\M$]\label{thm:EefnM}
Let $\mu_j>\frac54$. Then
\begin{equation}\label{eq:transport}
\Wc_j(Qf)=\frac{ih_j}{2}\Wc_j\left(\left(Q^2+U_+^2\right)f\right)\quad\text{for }f\in C^\infty(S\M).
\end{equation}
Consequently, for $f\in C^\infty(S\M)$ and $t_0,t_1\in\R$,
\begin{equation}\label{eq:Duhamel}
\Wc_j\left(f\circ a_{t_1}\right)-\Wc_j\left(f\circ a_{t_0}\right)=\frac{ih_j}{2}\int_{t_0}^{t_1}\Wc_j\left(\left(Q^2+U_+^2\right)\left(f\circ a_s\right)\right)\,ds.
\end{equation}
\end{thm}

\begin{proof}
Let $f\in C^\infty(S\M)$. Pairing the identity in Theorem \ref{thm:Phiode} with $f$, and using that the distributional adjoints of $Q$ and $U_+$ are $-Q$ and $-U_+$, we obtain
\[\Wc_j\left(\left(Q^2+U_+^2\right)f\right)+2i\lambda_j\Wc_j(Qf)=0.\]
This is \eqref{eq:transport}, since $-(2i\lambda_j)^{-1}=\frac{ih_j}{2}$.

The map $s\mapsto f\circ a_s$ is $C^1$ from $\R$ to $C^\infty(S\M)$, with derivative $Q(f\circ a_s)$, and $\Wc_j$ is a distribution on $S\M$. Hence, by \eqref{eq:transport},
\[\frac{d}{ds}\Wc_j\left(f\circ a_s\right)=\Wc_j\left(Q\left(f\circ a_s\right)\right)=\frac{ih_j}{2}\Wc_j\left(\left(Q^2+U_+^2\right)\left(f\circ a_s\right)\right).\]
Integrating from $t_0$ to $t_1$ gives \eqref{eq:Duhamel}.
\end{proof}

By Theorem \ref{thm:QU},
\begin{equation}\label{eq:QUa}
\left(Q^2+U_+^2\right)\left(f\circ a_s\right)=\left(Q^2f\right)\circ a_s+e^{-2s}\left(U_+^2f\right)\circ a_s\quad\text{for }f\in C^\infty(S\M),s\in\R.
\end{equation}
Together with Theorem \ref{thm:EefnM} for $t_0=0$, this gives \eqref{eq:exact}.

\section{Representation theory}\label{sec:rep}
In this section, we recall the representation theory of $\PSU(1,1)$ on $L^2(S\M)$ and use it to prove a uniform bound and the positivity of the matrix elements $\Wc_j$. Recall that $\M=\Gamma\bs\D=\Gamma\bs\PSU(1,1)/\SO(2)$ and $S\M=\Gamma\bs\PSU(1,1)$. The right regular representation of $\PSU(1,1)$ on $L^2(S\M)$ decomposes as a Hilbert direct sum of irreducible unitary representations. We only need the principal series representations generated by the Laplacian eigenfunctions; see Knapp \cite[Chapter II]{Kn} for the general representation theory.

\subsection{Raising and lowering operators}
For $\gamma\in\PSU(1,1)$, denote by $R_\gamma$ the right translation
\[\left(R_\gamma f\right)(x)=f(x\gamma)\quad\text{for }f\in L^2(S\M),x\in S\M.\]
Since the Liouville measure $d\Vol_L$ on $S\M$ is invariant under right translations, $R_\gamma$ is unitary on $L^2(S\M)$. The derived action of $Y\in\suf(1,1)$ is the Lie derivative $Y$ in Section \ref{sec:Lie}. Since the representation is unitary, $Y$ is skew-symmetric on $C^\infty(S\M)$. We now introduce the raising and lowering operators.

\begin{defn}[The raising and lowering operators]
In the complexified Lie algebra of $\suf(1,1)$, define the raising and lowering operators $E_+$ and $E_-$ by
\[E_+=Q+iV=\begin{pmatrix}0 & i\\
0 & 0\end{pmatrix}\quad\text{and}\quad E_-=Q-iV=\begin{pmatrix}0 & 0\\
-i & 0\end{pmatrix}.\]
\end{defn}

Using the commutation relations $[W,V]=Q$ and $[Q,W]=V$, we obtain
\[\left[W,E_+\right]=[W,Q]+i[W,V]=-V+iQ=i(Q+iV)=iE_+,\]
and similarly $[W,E_-]=-iE_-$. Hence, if $f\in C^\infty(S\M)$ has weight $m\in\Z$, i.e., $Wf=imf$, then
\[W\left(E_+f\right)=i(m+1)\left(E_+f\right)\quad\text{and}\quad W\left(E_-f\right)=i(m-1)\left(E_-f\right).\]
Thus, $E_+$ raises the weight by one, while $E_-$ lowers the weight by one.

\begin{prop}\label{prop:RL}\leavevmode
\begin{enumerate}[(i)]
\item With respect to $L^2(S\M)$,
\[E_+^\star=-E_-\quad\text{and}\quad E_-^\star=-E_+.\]
\item
\[E_+E_-=Q^2+V^2-iW\quad\text{and}\quad E_-E_+=Q^2+V^2+iW.\]
\item Let $f_1,f_2\in C^\infty(S\M)$ and $m_1,m_2\in\Z$. Suppose that $Wf_1=im_1f_1$ and $Wf_2=im_2f_2$ with $m_1\ne m_2$. Then
\[\left\langle f_1,f_2\right\rangle_{L^2(S\M)}=0.\]
In particular, if $m_1\ne0$, then
\[\int_{S\M}f_1\,d\Vol_L=0.\]
\item Let $f\in C^\infty(S\M)$, $\mu\in\R$, and $m\in\Z$. Suppose that $\Omega f=\mu f$ and $Wf=imf$. Let $\lambda^2=\mu-\frac14$. Then
\[\left\|E_+f\right\|_{L^2(S\M)}^2=\left(\mu+m^2+m\right)\|f\|_{L^2(S\M)}^2=\left[\lambda^2+\left(m+\frac12\right)^2\right]\|f\|_{L^2(S\M)}^2\]
and
\[\left\|E_-f\right\|_{L^2(S\M)}^2=\left(\mu+m^2-m\right)\|f\|_{L^2(S\M)}^2=\left[\lambda^2+\left(m-\frac12\right)^2\right]\|f\|_{L^2(S\M)}^2.\]
\end{enumerate}
\end{prop}

\begin{proof}\leavevmode
\begin{enumerate}[(i)]
\item Since $Q,V,W$ act as skew-symmetric operators, we have $Q^\star=-Q$ and $V^\star=-V$. Hence,
\[E_+^\star=(Q+iV)^\star=-Q+iV=-(Q-iV)=-E_-,\]
and similarly $E_-^\star=-E_+$.

\item Using $[Q,V]=W$, we have
\[E_+E_-=(Q+iV)(Q-iV)=Q^2+V^2+i[V,Q]=Q^2+V^2-iW.\]
The second identity is similar.

\item Since the right regular representation is unitary,
\[\left\langle f_1,f_2\right\rangle_{L^2(S\M)}=\left\langle e^{\theta W}f_1,e^{\theta W}f_2\right\rangle_{L^2(S\M)}=e^{i\theta\left(m_1-m_2\right)}\left\langle f_1,f_2\right\rangle_{L^2(S\M)}.\]
Since $m_1\ne m_2$, this implies $\langle f_1,f_2\rangle_{L^2(S\M)}=0$. Taking $f_2=1$ gives the final statement.

\item Using (i) and (ii), we compute
\begin{align*}
\left\|E_+f\right\|_{L^2(S\M)}^2
&=-\left\langle E_-E_+f,f\right\rangle_{L^2(S\M)}\\
&=\left\langle\left(-Q^2-V^2-iW\right)f,f\right\rangle_{L^2(S\M)}\\
&=\left\langle\left(\Omega-W^2-iW\right)f,f\right\rangle_{L^2(S\M)}\\
&=\left(\mu+m^2+m\right)\|f\|_{L^2(S\M)}^2\\
&=\left[\lambda^2+\left(m+\frac12\right)^2\right]\|f\|_{L^2(S\M)}^2.
\end{align*}
The formula for $E_-f$ follows similarly.
\end{enumerate}
\end{proof}

\subsection{The principal series}
Let
\[\mu_j\ge\frac14\quad\text{and}\quad\lambda_j=\sqrt{\mu_j-\frac14}\ge0\]
as in \eqref{eq:lambda}. Denote by $v_{j,0}$ the lift of $\vp_j$ to $S\M$, which has $L^2(S\M)$ norm $1$. By Theorem \ref{thm:LonD},
\[\Omega v_{j,0}=\mu_jv_{j,0}\quad\text{and}\quad Wv_{j,0}=0.\]
Starting from $v_{j,0}$, define $v_{j,m}$ inductively by
\begin{equation}\label{eq:RLP}
\begin{cases}
v_{j,m+1}=\frac{1}{i\lambda_j+m+\frac12}E_+v_{j,m}, & \text{if }m\ge0,\\
v_{j,m-1}=\frac{1}{i\lambda_j-m+\frac12}E_-v_{j,m}, & \text{if }m\le0.
\end{cases}
\end{equation}
Since $\Omega$ commutes with $Q,V,W$, we have
\[\Omega v_{j,m}=\mu_jv_{j,m}\quad\text{and}\quad Wv_{j,m}=imv_{j,m}.\]
In fact, the recursions in \eqref{eq:RLP} hold for all $m\in\Z$ in both directions, i.e.,
\begin{equation}\label{eq:RLall}
E_+v_{j,m}=\left(i\lambda_j+m+\frac12\right)v_{j,m+1}\quad\text{and}\quad E_-v_{j,m+1}=\left(i\lambda_j-m-\frac12\right)v_{j,m}\quad\text{for }m\in\Z.
\end{equation}
Indeed, by (ii) of Proposition \ref{prop:RL} and $\Omega=-Q^2-V^2+W^2$,
\[E_-E_+=-\Omega+W^2+iW\quad\text{and}\quad E_+E_-=-\Omega+W^2-iW.\]
Notice that
\[\mu_j+m^2+m=\lambda_j^2+\left(m+\frac12\right)^2=\left(i\lambda_j+m+\frac12\right)\left(-i\lambda_j+m+\frac12\right).\]
If $m\ge0$, then the first identity in \eqref{eq:RLall} is the definition of $v_{j,m+1}$, and
\[E_-v_{j,m+1}=\frac{E_-E_+v_{j,m}}{i\lambda_j+m+\frac12}=-\frac{\mu_j+m^2+m}{i\lambda_j+m+\frac12}v_{j,m}=\left(i\lambda_j-m-\frac12\right)v_{j,m}.\]
If $m\le-1$, then the second identity in \eqref{eq:RLall} is the definition of $v_{j,m}$, and the first identity follows in the same way from $E_+E_-$. By \eqref{eq:RLall} and (iv) of Proposition \ref{prop:RL},
\[\left\|v_{j,m+1}\right\|^2_{L^2(S\M)}=\frac{\left\|E_+v_{j,m}\right\|^2_{L^2(S\M)}}{\lambda_j^2+\left(m+\frac12\right)^2}=\left\|v_{j,m}\right\|^2_{L^2(S\M)}\quad\text{for }m\in\Z.\]
Since different weights are orthogonal by (iii) of Proposition \ref{prop:RL}, the set $\{v_{j,m}\}_{m=-\infty}^\infty$ is orthonormal in $L^2(S\M)$. The closure $\Hc_j$ of the span of $\{v_{j,m}\}_{m=-\infty}^\infty$ in $L^2(S\M)$ is an invariant subspace of the right regular representation. It realizes the spherical unitary principal series representation of $\PSU(1,1)$ with spectral parameter $\lambda_j$, generated by $\vp_j$; see Knapp \cite[Chapter II]{Kn}.

\subsection{Semiclassical pseudo-differential operators by representation theory}\label{sec:SDOrep}
We provide another view of the action of semiclassical pseudo-differential operators on Laplacian eigenfunctions, using the principal series $\Hc_j$.

Let $\mu_j>\frac54$. Since $v_{j,m}$ has weight $m$, i.e.,
\[v_{j,m}\left(xk_\theta\right)=e^{im\theta}v_{j,m}(x),\]
its fiber average vanishes when $m\ne0$. Thus, only the weight-zero component projects to the Laplacian eigenfunction $\vp_j$ on $\M$. The following $K$-type expansion of $\Phi_j$ is due to Zelditch \cite[Proposition 2.2]{Ze2}. We include a proof and track the constant in the proof, as our normalization of the measures is different from \cite{Ze2}.

\begin{thm}\label{thm:SDOrep}
Let $\mu_j>\frac54$. Then
\[\Phi_j=\frac{1}{2\pi}\sum_{m=-\infty}^\infty v_{j,m}\quad\text{in }\Dc'(S\M).\]
Consequently,
\[\Wc_j(f)=\sum_{m=-\infty}^\infty\left\langle fv_{j,m},v_{j,0}\right\rangle_{L^2(S\M)}\quad\text{for }f\in C^\infty(S\M).\]
\end{thm}

\begin{proof}
For $m\in\Z$, let $\Phi_{j,m}$ be the component of $\Phi_j$ of weight $m$, that is,
\[\Phi_{j,m}(x)=\frac{1}{2\pi}\int_0^{2\pi}\Phi_j\left(xk_\theta\right)e^{-im\theta}\,d\theta\]
in the sense of distributions. Then
\[W\Phi_{j,m}=im\Phi_{j,m}\quad\text{and}\quad\Phi_j=\sum_{m=-\infty}^\infty\Phi_{j,m}\text{ in }\Dc'(S\M).\]
By \eqref{eq:fiber}, the weight-zero component of $\tilde\Phi_j$ at $z\in\D$ is
\[\frac{1}{2\pi}\int_\Sb e^{\left(i\lambda_j-\frac12\right)\langle z,b\rangle}\tau_j(b)e^{\langle z,b\rangle}\,db.\]
By \eqref{eq:tildevp}, this is $\frac{1}{2\pi}\tilde\vp_j(z)$. Thus, $\Phi_{j,0}=\frac{1}{2\pi}v_{j,0}$.

Since $V=U_+-W$, Theorem \ref{thm:Phi} gives
\[E_\pm\Phi_j=Q\Phi_j\pm iV\Phi_j=\left(i\lambda_j-\frac12\right)\Phi_j\mp iW\Phi_j.\]
As $E_\pm$ shifts the weights by $\pm1$, comparing the components of weight $m\pm1$ gives
\[E_+\Phi_{j,m}=\left(i\lambda_j+m+\frac12\right)\Phi_{j,m+1}\quad\text{and}\quad E_-\Phi_{j,m}=\left(i\lambda_j-m+\frac12\right)\Phi_{j,m-1}.\]
These are the recursions in \eqref{eq:RLP}, whose coefficients do not vanish. Starting from $\Phi_{j,0}=\frac{1}{2\pi}v_{j,0}$, induction gives $\Phi_{j,m}=\frac{1}{2\pi}v_{j,m}$ for $m\in\Z$.

Finally, let $f\in C^\infty(S\M)$. Since $f\ol{v_{j,0}}$ is smooth on $S\M$, the $K$-type expansion of $\Phi_j$ gives
\[\Wc_j(f)=2\pi\int_{S\M}f\Phi_j\ol{\vp_j}\,d\Vol_L=\sum_{m=-\infty}^\infty\int_{S\M}fv_{j,m}\ol{v_{j,0}}\,d\Vol_L=\sum_{m=-\infty}^\infty\left\langle fv_{j,m},v_{j,0}\right\rangle_{L^2(S\M)}.\]
\end{proof}

\subsection{A uniform bound}\label{sec:norms}
For $n\in\N$ and $f\in C^\infty(S\M)$, define
\[\|f\|_{C^n(S\M)}=\max_{0\le l\le n}\max_{X_1,\dots,X_l\in\left\{Q,U_+,U_-\right\}}\left\|X_1\cdots X_lf\right\|_{L^\infty(S\M)}\]
and
\[\|f\|_{H^n(S\M)}^2=\sum_{l=0}^n\sum_{X_1,\dots,X_l\in\left\{Q,U_+,U_-\right\}}\left\|X_1\cdots X_lf\right\|_{L^2(S\M)}^2,\]
where all orderings of the vector fields are allowed, and the terms with $l=0$ are $\|f\|_{L^\infty(S\M)}$ and $\|f\|_{L^2(S\M)}^2$, respectively. Since $\{Q,U_+,U_-\}$ is a global frame on $S\M$, these norms are equivalent to the standard $C^n$ and Sobolev norms on $S\M$. Denote by $H^n(S\M)$ the completion of $C^\infty(S\M)$ with respect to $\|\cdot\|_{H^n(S\M)}$. By the Sobolev embedding on the $3$-dimensional compact manifold $S\M$,
\begin{equation}\label{eq:Sobolev}
\|f\|_{C^2(S\M)}\le C\|f\|_{H^4(S\M)}\quad\text{for }f\in C^\infty(S\M),
\end{equation}
where $C=C(\M)>0$; see Zworski \cite[Section 14.2]{Zw}. Moreover, by Theorem \ref{thm:QU},
\begin{equation}\label{eq:Cnflow}
\left\|f\circ a_t\right\|_{C^n(S\M)}\le e^{n|t|}\|f\|_{C^n(S\M)}\quad\text{for }f\in C^\infty(S\M),t\in\R.
\end{equation}
\begin{lemma}[Uniform bound]\label{lemma:apriori}
Let $\mu_j>\frac54$. Then for $f\in C^\infty(S\M)$,
\[\left|\Wc_j(f)\right|\le\|f\|_{L^\infty(S\M)}+\frac{\pi^2}{3}\left\|W^2f\right\|_{L^\infty(S\M)}\le5\|f\|_{C^2(S\M)}.\]
\end{lemma}

\begin{proof}
Let $f\in C^\infty(S\M)$. For $m\in\Z$, let $f_m$ be the component of $f$ of weight $m$, i.e.,
\[f_m(x)=\frac{1}{2\pi}\int_0^{2\pi}f\left(xk_\theta\right)e^{-im\theta}\,d\theta.\]
Then $Wf_m=imf_m$. Since $W^2f_m=-m^2f_m$ is the component of $W^2f$ of weight $m$, we have that if $m\ne0$, then
\[\left\|f_m\right\|_{L^\infty(S\M)}\le\frac{1}{m^2}\left\|W^2f\right\|_{L^\infty(S\M)}.\]
Hence, $f=\sum_mf_m$ with uniform convergence on $S\M$. Now $f_{m_1}v_{j,m_2}\ol{v_{j,0}}$ has weight $m_1+m_2$, so its integral vanishes unless $m_1=-m_2$ by (iii) of Proposition \ref{prop:RL}. Notice that
\[\left\|f_0\right\|_{L^\infty(S\M)}\le\|f\|_{L^\infty(S\M)}\quad\text{and}\quad\sum_{m\ne0}\frac{1}{m^2}=\frac{\pi^2}{3}.\]
Therefore, by Theorem \ref{thm:SDOrep} and the orthonormality of $\{v_{j,m}\}_{m=-\infty}^\infty$,
\begin{align*}
\left|\Wc_j(f)\right|&=\left|\sum_{m=-\infty}^\infty\left\langle f_{-m}v_{j,m},v_{j,0}\right\rangle_{L^2(S\M)}\right|\\
&\le\sum_{m=-\infty}^\infty\left\|f_m\right\|_{L^\infty(S\M)}\\
&=\left\|f_0\right\|_{L^\infty(S\M)}+\sum_{m\ne0}\left\|f_m\right\|_{L^\infty(S\M)}\\
&\le\left\|f_0\right\|_{L^\infty(S\M)}+\sum_{m\ne0}\frac{1}{m^2}\left\|W^2f\right\|_{L^\infty(S\M)}\\
&\le\|f\|_{L^\infty(S\M)}+\frac{\pi^2}{3}\left\|W^2f\right\|_{L^\infty(S\M)}.
\end{align*}
Finally, $W=\frac12(U_+-U_-)$, so
\[\left\|W^2f\right\|_{L^\infty(S\M)}\le\|f\|_{C^2(S\M)}.\]
Since $1+\frac{\pi^2}{3}<5$, this completes the proof.
\end{proof}

The first estimate in Lemma \ref{lemma:apriori} depends only on $f$ and its second fiber derivative $W^2f$, and its constant is independent of $j$. For functions transported by the geodesic flow, \eqref{eq:Cnflow} then gives
\begin{equation}\label{eq:aprioriflow}
\left|\Wc_j\left(f\circ a_t\right)\right|\le5e^{2|t|}\|f\|_{C^2(S\M)}\quad\text{for }f\in C^\infty(S\M),t\in\R.
\end{equation}
\begin{rmk}
If $\Op_{h_j}^\Gamma(f)$ extends to a bounded operator on $L^2(\M)$, then $\Wc_j(f)$ can also be bounded by its $L^2$ operator norm. However, standard $L^2$ boundedness estimates for semiclassical pseudo-differential operators, such as the Calder\'on-Vaillancourt theorem in \cite[Theorem 4.23]{Zw}, involve more derivatives of the symbols than the one in Lemma \ref{lemma:apriori}, with the required number depending on the dimension. This matters in Sections \ref{sec:forward} and \ref{sec:backward}, where we apply Lemma \ref{lemma:apriori} to $f\circ a_t$: each derivative of $f\circ a_t$ costs a factor $e^{|t|}$, and the factor $e^{2|t|}$ in \eqref{eq:aprioriflow} is exactly compensated by $h$ at $|t|=T_E$, i.e., $he^{2T_E}=1$. Moreover, Lemma \ref{lemma:apriori} only concerns the diagonal matrix elements of the eigenfunctions, which are defined by Theorem \ref{thm:SDOefnM} without extending $\Op_h^\Gamma$ to $L^2(\M)$.
\end{rmk}

\subsection{Positivity}
The matrix elements $\Wc_j$ are distributions on $S\M$, but are not positive, i.e., $\Wc_j(f)$ need not be non-negative for non-negative $f$. However, $\Wc_j$ is close to a probability measure on $S\M$, as we now show. The positivity is used in Sections \ref{sec:beyond} and \ref{sec:rate}: there are many distributions on $S\M$ which are invariant under the stable horocycle flow by Flaminio-Forni \cite{FF}, but $d\Vol_L$ is the only invariant probability measure by Furstenberg \cite{Fu}.

The positivity result here is not new in the literature. For example, in the standard quantization, the positivity of the semiclassical measures follows from the sharp G\r{a}rding inequality; see Zworski \cite[Section 5.1]{Zw}. In addition, Zelditch \cite{Ze2} obtained the positivity of the semiclassical measures by the Friedrichs quantization, which is positive and differs from $\Op^\Gamma_h(f)$ by $O_f(h)$ for each fixed symbol $f$. Here, the implicit constant in $O_f(h)$ depends on finitely many derivatives of $f$.

We show that $\Wc_j$ is close to a probability measure with a polynomial rate in $h_j$, uniformly in $f$, with a bound involving only two derivatives of $f$, as in Lemma \ref{lemma:apriori}. The proof follows an approach similar to the construction of microlocal lifts using raising and lowering operators by Wolpert \cite{Wo2} and Lindenstrauss \cite{L1}; see also Silberman-Venkatesh \cite{SV}.

\begin{lemma}[Positivity]\label{lemma:positive}
Let $\mu_j>\frac54$ and $m_0\in\N$ with $m_0\ge1$. Let
\[\Psi_{j,m_0}=\frac{1}{\sqrt{2m_0+1}}\sum_{m=-m_0}^{m_0}v_{j,m}\quad\text{and}\quad d\nu_{j,m_0}=\left|\Psi_{j,m_0}\right|^2d\Vol_L.\]
Then $\nu_{j,m_0}$ is a probability measure on $S\M$, and for $f\in C^\infty(S\M)$,
\[\left|\Wc_j(f)-\int_{S\M}f\,d\nu_{j,m_0}\right|\le C\log\left(2+m_0\right)\left(m_0^{-1}+h_jm_0\right)\|f\|_{C^2(S\M)},\]
where $C>0$ is an absolute constant. In particular, taking $m_0=\lfloor h_j^{-\frac12}\rfloor$, there is a probability measure $\nu_j$ on $S\M$ such that
\begin{equation}\label{eq:positive}
\left|\Wc_j(f)-\int_{S\M}f\,d\nu_j\right|\le Ch_j^{\frac12}\left(1+\left|\log h_j\right|\right)\|f\|_{C^2(S\M)}\quad\text{for }f\in C^\infty(S\M).
\end{equation}
\end{lemma}

\begin{proof}
Since $\{v_{j,m}\}_{m=-\infty}^\infty$ is orthonormal, $\|\Psi_{j,m_0}\|_{L^2(S\M)}=1$. So $\nu_{j,m_0}$ is a probability measure on $S\M$. Let $f\in C^\infty(S\M)$. For $m_1,m_2\in\Z$, write
\[c_{m_1,m_2}=\left\langle fv_{j,m_1},v_{j,m_2}\right\rangle_{L^2(S\M)}.\]
Then $|c_{m_1,m_2}|\le\|f\|_{L^\infty(S\M)}$, and
\[\int_{S\M}f\,d\nu_{j,m_0}=\left\langle f\Psi_{j,m_0},\Psi_{j,m_0}\right\rangle_{L^2(S\M)}=\frac{1}{2m_0+1}\sum_{m_1,m_2=-m_0}^{m_0}c_{m_1,m_2}.\]
We compare this sum with $\Wc_j(f)$ in three steps.

\textit{Step 1}. We show that for $m_1,m_2\in\Z$,
\begin{equation}\label{eq:amn}
\left|c_{m_1+1,m_2+1}-c_{m_1,m_2}\right|\le\frac{Ch_j}{1+|m_1-m_2|}\|f\|_{C^2(S\M)}.
\end{equation}
By \eqref{eq:RLall},
\[E_+\left(fv_{j,m_1}\right)=\left(E_+f\right)v_{j,m_1}+\left(i\lambda_j+m_1+\frac12\right)fv_{j,m_1+1},\]
and
\[E_-v_{j,m_2+1}=\left(i\lambda_j-m_2-\frac12\right)v_{j,m_2}.\]
Since $E_+^\star=-E_-$ by (i) of Proposition \ref{prop:RL}, we have
\[\left\langle E_+\left(fv_{j,m_1}\right),v_{j,m_2+1}\right\rangle_{L^2(S\M)}=-\left\langle fv_{j,m_1},E_-v_{j,m_2+1}\right\rangle_{L^2(S\M)}=\left(i\lambda_j+m_2+\frac12\right)c_{m_1,m_2}.\]
Hence,
\[\left(i\lambda_j+m_1+\frac12\right)c_{m_1+1,m_2+1}=\left(i\lambda_j+m_2+\frac12\right)c_{m_1,m_2}-\left\langle\left(E_+f\right)v_{j,m_1},v_{j,m_2+1}\right\rangle_{L^2(S\M)}.\]
That is,
\[c_{m_1+1,m_2+1}-c_{m_1,m_2}=\frac{\left(m_2-m_1\right)c_{m_1,m_2}-\left\langle\left(E_+f\right)v_{j,m_1},v_{j,m_2+1}\right\rangle_{L^2(S\M)}}{i\lambda_j+m_1+\frac12}.\]
Let $l=m_2-m_1$. Denote by $f_l$ the $K$-type component of $f$ of weight $l$. Weight orthogonality shows that $c_{m_1,m_2}$ depends only on $f_l$, while the last matrix element depends only on $E_+f_l=(E_+f)_{l+1}$. Now
\[\left|\left(m_2-m_1\right)c_{m_1,m_2}\right|\le|l|\|f_l\|_{L^\infty(S\M)}\le\frac{C}{1+|l|}\left\|W^2f\right\|_{L^\infty(S\M)}\le\frac{C}{1+|l|}\|f\|_{C^2(S\M)}.\]
For any $g\in C^\infty(S\M)$ and $k\in\Z$, its $K$-type component $g_k$ satisfies
\[\|g_k\|_{L^\infty(S\M)}\le\frac{C}{1+|k|}\left(\|g\|_{L^\infty(S\M)}+\|Wg\|_{L^\infty(S\M)}\right).\]
Indeed, for $k\ne0$, this follows from $Wg_k=ikg_k$, while for $k=0$, it follows from $\|g_0\|_{L^\infty(S\M)}\le\|g\|_{L^\infty(S\M)}$.

Applying this estimate to $g=E_+f$ and $k=l+1$, and using $1+|l|\le2(1+|l+1|)$, we obtain
\[\left\|\left(E_+f\right)_{l+1}\right\|_{L^\infty(S\M)}\le\frac{C}{1+|l|}\left(\|E_+f\|_{L^\infty(S\M)}+\|WE_+f\|_{L^\infty(S\M)}\right)\le\frac{C}{1+|l|}\|f\|_{C^2(S\M)}.\]
Consequently,
\[\left|\left\langle\left(E_+f\right)v_{j,m_1},v_{j,m_2+1}\right\rangle_{L^2(S\M)}\right|
\le\frac{C}{1+|l|}\|f\|_{C^2(S\M)}.\]
Thus, the numerator in the formula for $c_{m_1+1,m_2+1}-c_{m_1,m_2}$ satisfies
\[\left|l\,c_{m_1,m_2}\right|+\left|\left\langle(E_+f)v_{j,m_1},v_{j,m_2+1}\right\rangle_{L^2(S\M)}\right|\le\frac{C}{1+|l|}\|f\|_{C^2(S\M)}.\]
Since
\[\left|i\lambda_j+m_1+\frac12\right|\ge\lambda_j=h_j^{-1},\]
dividing by this coefficient proves \eqref{eq:amn}.

\textit{Step 2}. Let $-m_0\le m_1,m_2\le m_0$. Moving along the diagonal $|m_2|$ times, \eqref{eq:amn} gives
\[\left|c_{m_1,m_2}-c_{m_1-m_2,0}\right|\le\frac{Ch_j|m_2|}{1+|m_1-m_2|}\|f\|_{C^2(S\M)}.\]
Since the number of $(m_1,m_2)$ with $-m_0\le m_1,m_2\le m_0$ and $m_1-m_2=l$ is $2m_0+1-|l|$, we obtain
\begin{align*}
\left|\int_{S\M}f\,d\nu_{j,m_0}-\sum_{l=-2m_0}^{2m_0}\left(1-\frac{|l|}{2m_0+1}\right)c_{l,0}\right|
&\le Ch_jm_0\sum_{|l|\le2m_0}\frac{1}{1+|l|}\|f\|_{C^2(S\M)}\\
&\le Ch_jm_0\log\left(2+m_0\right)\|f\|_{C^2(S\M)}.
\end{align*}
That is, only the coefficients $c_{l,0}$ remain, with the weights $1-\frac{|l|}{2m_0+1}$.

\textit{Step 3}. By Theorem \ref{thm:SDOrep}, $\Wc_j(f)=\sum_{l=-\infty}^\infty c_{l,0}$. As in the proof of Lemma \ref{lemma:apriori}, $|c_{l,0}|\le l^{-2}\|W^2f\|_{L^\infty(S\M)}$ for $l\ne0$. Hence,
\begin{align*}
&\quad\left|\Wc_j(f)-\sum_{l=-2m_0}^{2m_0}\left(1-\frac{|l|}{2m_0+1}\right)c_{l,0}\right|\\
&\le\sum_{1\le|l|\le2m_0}\frac{|l|}{2m_0+1}\cdot\frac{\|W^2f\|_{L^\infty(S\M)}}{l^2}+\sum_{|l|>2m_0}\frac{\|W^2f\|_{L^\infty(S\M)}}{l^2}\\
&\le\left(\frac{1+\log(2m_0)}{m_0}+\frac{1}{m_0}\right)\|W^2f\|_{L^\infty(S\M)}.
\end{align*}
Since $\|W^2f\|_{L^\infty(S\M)}\le\|f\|_{C^2(S\M)}$, Steps 2 and 3 prove the first estimate. Finally, let $m_0=\lfloor h_j^{-\frac12}\rfloor$. Then $m_0\ge\frac12h_j^{-\frac12}$, and
\[\log\left(2+m_0\right)\left(m_0^{-1}+h_jm_0\right)\le Ch_j^{\frac12}\left(1+\left|\log h_j\right|\right).\]
This completes the proof.
\end{proof}

\section{Egorov's theorem forward in time}\label{sec:forward}
In this section, we prove (i) of Theorem \ref{thm:Egorov}. Throughout the rest of the paper, we assume that $\mu_j>\frac54$, so that
\[h_j=\lambda_j^{-1}\in(0,1)\quad\text{and}\quad T_E=\frac12\left|\log h_j\right|.\]
Forward in time, the weight $e^{-2s}$ in \eqref{eq:exact} is bounded by $1$. Thus, for $s\ge0$,
\begin{equation}\label{eq:Lsbound}
\left\|\left(Q^2+e^{-2s}U_+^2\right)f\right\|_{C^n(S\M)}\le2\|f\|_{C^{n+2}(S\M)}\quad\text{for }f\in C^\infty(S\M),n\in\N.
\end{equation}
So we can apply \eqref{eq:exact} repeatedly, losing two derivatives each time. This gives the following proposition.

\begin{prop}\label{prop:iterate}
Let $k\in\N$ with $k\ge1$. Then for $f\in C^\infty(S\M)$ and $t\ge0$,
\[\left|\Wc_j\left(f\circ a_t\right)-\Wc_j(f)\right|\le5\left(\sum_{l=1}^{k-1}\left(h_jt\right)^l+\left(h_jt\right)^ke^{2t}\right)\|f\|_{C^{2k+2}(S\M)}.\]
Here, the sum is $0$ if $k=1$.
\end{prop}

\begin{proof}
We use induction on $k$. Let $f\in C^\infty(S\M)$ and $t\ge0$. For $0\le s\le t$, write
\[u_s=\left(Q^2+e^{-2s}U_+^2\right)f.\]
Then the exact form of Egorov's theorem reads
\begin{equation}\label{eq:exactus}
\Wc_j\left(f\circ a_t\right)-\Wc_j(f)=\frac{ih_j}{2}\int_0^t\Wc_j\left(u_s\circ a_s\right)\,ds,
\end{equation}
where by \eqref{eq:Lsbound},
\[\left\|u_s\right\|_{C^n(S\M)}\le2\|f\|_{C^{n+2}(S\M)}.\]
For $k=1$, \eqref{eq:aprioriflow} gives
\[\left|\Wc_j\left(u_s\circ a_s\right)\right|\le5e^{2s}\left\|u_s\right\|_{C^2(S\M)}\le10e^{2t}\|f\|_{C^4(S\M)}.\]
Then \eqref{eq:exactus} gives
\[\left|\Wc_j\left(f\circ a_t\right)-\Wc_j(f)\right|\le5h_jte^{2t}\|f\|_{C^4(S\M)}.\]
This proves the proposition for $k=1$.

Now assume that the estimate holds for $k$. To prove it for $k+1$, write
\[\Wc_j\left(u_s\circ a_s\right)=\Wc_j\left(u_s\right)+\left[\Wc_j\left(u_s\circ a_s\right)-\Wc_j\left(u_s\right)\right].\]
By Lemma \ref{lemma:apriori},
\[\left|\Wc_j\left(u_s\right)\right|\le5\left\|u_s\right\|_{C^2(S\M)}\le10\|f\|_{C^4(S\M)}.\]
By the proposition for $k$, applied to $u_s$ at time $s\le t$,
\begin{align*}
\left|\Wc_j\left(u_s\circ a_s\right)-\Wc_j\left(u_s\right)\right|&\le5\left(\sum_{l=1}^{k-1}\left(h_js\right)^l+\left(h_js\right)^ke^{2s}\right)\left\|u_s\right\|_{C^{2k+2}(S\M)}\\
&\le10\left(\sum_{l=1}^{k-1}\left(h_jt\right)^l+\left(h_jt\right)^ke^{2t}\right)\|f\|_{C^{2k+4}(S\M)}.
\end{align*}
Putting these two estimates into \eqref{eq:exactus}, we have
\begin{align*}
\left|\Wc_j\left(f\circ a_t\right)-\Wc_j(f)\right|&\le\frac{h_j}{2}\int_0^t\left|\Wc_j\left(u_s\circ a_s\right)\right|\,ds\\
&\le\frac{h_j}{2}\int_0^t\left(\left|\Wc_j\left(u_s\right)\right|+\left|\Wc_j\left(u_s\circ a_s\right)-\Wc_j\left(u_s\right)\right|\right)\,ds\\
&\le5\left(h_jt\right)\|f\|_{C^4(S\M)}+5\|f\|_{C^{2k+4}(S\M)}\left(\sum_{l=2}^k\left(h_jt\right)^l+\left(h_jt\right)^{k+1}e^{2t}\right)\\
&\le5\left(\sum_{l=1}^k\left(h_jt\right)^l+\left(h_jt\right)^{k+1}e^{2t}\right)\|f\|_{C^{2k+4}(S\M)}.
\end{align*}
\end{proof}

We are now ready to prove (i) of Theorem \ref{thm:Egorov}.

\begin{proof}[Proof of (i) of Theorem \ref{thm:Egorov}]
Let $\rho>0$, $k=\lfloor\rho\rfloor+2$, and $0\le t\le\rho T_E$. For $0\le s\le t$, write
\[u_s=\left(Q^2+e^{-2s}U_+^2\right)f.\]
By \eqref{eq:exact},
\begin{align*}
\Wc_j\left(f\circ a_t\right)-\Wc_j(f)
&=\frac{ih_j}{2}\left[t\Wc_j(Q^2f)+\frac{1-e^{-2t}}{2}\Wc_j(U_+^2f)\right]\\
&\quad+\frac{ih_j}{2}\int_0^t\left[\Wc_j\left(u_s\circ a_s\right)-\Wc_j(u_s)\right] \,ds.
\end{align*}
Applying \eqref{eq:transport} to $Qf$, we obtain
\[\Wc_j(Q^2f)=\frac{ih_j}{2}\Wc_j\left(\left(Q^2+U_+^2\right)Qf\right).\]
The expression inside $\Wc_j$ contains three derivatives of $f$. Hence, Lemma \ref{lemma:apriori} gives
\[\left|\Wc_j(Q^2f)\right|\le Ch_j\|f\|_{C^5(S\M)}.\]
Applying the same lemma directly to $U_+^2f$ gives
\[\left|\Wc_j(U_+^2f)\right|\le C\|f\|_{C^4(S\M)}.\]
Moreover, Proposition \ref{prop:iterate} with $k-1$, together with \eqref{eq:Lsbound}, gives
\[\left|\Wc_j\left(u_s\circ a_s\right)-\Wc_j(u_s)\right|\le C\left(\sum_{l=1}^{k-2}(h_js)^l+(h_js)^{k-1}e^{2s}\right)\|f\|_{C^{2k+2}(S\M)},\]
where the sum is empty if $k=2$. Consequently,
\[\left|\Wc_j\left(f\circ a_t\right)-\Wc_j(f)\right|\le C\left[h_j+h_j^2t+\sum_{l=2}^{k-1}(h_jt)^l+(h_jt)^ke^{2t}\right]\|f\|_{C^{2k+2}(S\M)}.\]
Since $0\le t\le\rho T_E$, we have
\[t\le\frac{\rho}{2}\left|\log h_j\right|\quad\text{and}\quad e^{2t}\le h_j^{-\rho}.\]
In particular, for $2\le l\le k-1$,
\[\left(h_jt\right)^l\le C(\rho)h_j^l\left|\log h_j\right|^l\le C(\rho)h_j.\]
The last term satisfies
\[\left(h_jt\right)^ke^{2t}\le C(\rho)h_j^{k-\rho}\left|\log h_j\right|^k\le C(\rho)h_j,\]
since $k-\rho>1$. Also, $h_j^2t\le C(\rho)h_j$. Therefore, every term in the preceding bound is $O_\rho(h_j)$, finishing the proof.
\end{proof}

\begin{rmk}
The case $k=1$ of Proposition \ref{prop:iterate} is the direct estimate of the integral in \eqref{eq:exact} by the uniform bound in Lemma \ref{lemma:apriori}. It only gives Egorov's theorem up to $(1-\ve)T_E$, as in the standard quantization; see Section \ref{sec:Egorov}. Since each application of \eqref{eq:exact} loses two derivatives of $f$, a direct Gr\"onwall argument in a fixed $C^n$ norm does not close. Instead, we apply \eqref{eq:exact} $k$ times. Each application gains a factor $h_jt$ at the cost of two derivatives, and the factor $e^{2t}$ from \eqref{eq:aprioriflow} only appears at the last step. This gives Egorov's theorem forward in time to all logarithmic times. In Section \ref{sec:backward}, we apply the exact identity twice after absorbing the weight $e^{-2s}$ in \eqref{eq:exact} by the Schr\"odinger group along the stable horocycles.
\end{rmk}

\section{Egorov's theorem backward in time}\label{sec:backward}
In this section, we prove (ii) of Theorem \ref{thm:Egorov}. Backward in time, the weight $e^{-2s}$ in \eqref{eq:exact} grows and reaches $h_j^{-1}$ at $s=-T_E$. We establish the properties of the Schr\"odinger group along the stable horocycles defined below and derive \eqref{eq:exactSchr}, which absorbs this weight. Applying this identity twice, we prove the renormalized Egorov estimate in Theorem \ref{thm:renorm}, from which (ii) of Theorem \ref{thm:Egorov} follows.

\subsection{The Schr\"odinger group along the stable horocycles}\label{sec:Schrgroup}
Recall the right translations $R_\gamma$ for $\gamma\in\PSU(1,1)$ in Section \ref{sec:rep}. Then $R_{\gamma_1}R_{\gamma_2}=R_{\gamma_1\gamma_2}$, $R_{a_t}f=f\circ a_t$, and $R_{n_s}f=f\circ n_s$. Moreover,
\begin{equation}\label{eq:RAd}
YR_\gamma=R_\gamma\left(\Ad_{\gamma^{-1}}Y\right)\quad\text{on }C^\infty(S\M)\text{ for }Y\in\suf(1,1).
\end{equation}
By \eqref{eq:WUQ}, we have
\begin{equation}\label{eq:Adn}
\Ad_{n_{-s}}Q=Q+sU_+,\quad\Ad_{n_{-s}}U_+=U_+,\quad\Ad_{n_{-s}}U_-=U_--2sQ-s^2U_+.
\end{equation}
Since the stable horocycle flow preserves $d\Vol_L$, $\{R_{n_s}\}_{s\in\R}$ is a strongly continuous unitary group on $L^2(S\M)$. By Stone's theorem, its infinitesimal generator is skew-adjoint. On $C^\infty(S\M)$, this generator agrees with the vector field $U_+$. Moreover, $C^\infty(S\M)$ is an invariant core for the generator, so the generator is the closure of $U_+$ on $C^\infty(S\M)$; see Reed-Simon \cite[Section VIII.4]{ReSi}. We still denote the closure by $U_+$. Then $iU_+$ is self-adjoint, $R_{n_s}=e^{-is(iU_+)}$, and for a bounded Borel function $F$ on $\R$, the operator $F(iU_+)$ on $L^2(S\M)$ is defined by the spectral theorem \cite[Section VIII.3]{ReSi}.

\begin{defn}[The Schr\"odinger group along the stable horocycles]
For $\sigma\in\R$, the Schr\"odinger group along the stable horocycles is defined as
\[e^{-\frac{i\sigma}{4}U_+^2}=F_\sigma\left(iU_+\right),\quad\text{where }F_\sigma(\eta)=e^{\frac{i\sigma}{4}\eta^2}.\]
\end{defn}

By (iii) of Lemma \ref{lemma:Schr} below, for $f\in C^\infty(S\M)$, $u(\sigma)=e^{-\frac{i\sigma}{4}U_+^2}f$ satisfies $\partial_\sigma u(\sigma)=-\frac i4U_+^2u(\sigma)$ and $u(0)=f$, which is the free Schr\"odinger equation along each stable horocycle. This makes the formal definition \eqref{eq:Schr} precise.

If $F$ is in the Schwartz space $\Sc(\R)$, then by the Fourier inversion formula,
\begin{equation}\label{eq:FiU}
F\left(iU_+\right)f=\int_\R\check F(s)\,f\circ n_s\,ds\quad\text{for }f\in L^2(S\M),\quad\text{where }\check F(s)=\frac{1}{2\pi}\int_\R F(\eta)e^{is\eta}\,d\eta.
\end{equation}
\begin{lemma}\label{lemma:Schr}
The operators $e^{-\frac{i\sigma}{4}U_+^2}$, $\sigma\in\R$, form a unitary group on $L^2(S\M)$, which commutes with $R_{n_s}$ for $s\in\R$. Moreover, $e^{-\frac{i\sigma}{4}U_+^2}1=1$, and hence
\[\int_{S\M}e^{-\frac{i\sigma}{4}U_+^2}f\,d\Vol_L=\int_{S\M}f\,d\Vol_L\quad\text{for }f\in L^2(S\M).\]
Let $f\in C^\infty(S\M)$ and $\sigma\in\R$. Then the following hold.
\begin{enumerate}[(i)]
\item $e^{-\frac{i\sigma}{4}U_+^2}f\in C^\infty(S\M)$, and
\[Xe^{-\frac{i\sigma}{4}U_+^2}f=e^{-\frac{i\sigma}{4}U_+^2}X_\sigma f\quad\text{for }X\in\left\{Q,U_+,U_-\right\},\]
where
\[Q_\sigma=Q-\frac{i\sigma}{2}U_+^2,\quad U_{+,\sigma}=U_+,\quad U_{-,\sigma}=U_-+i\sigma U_+Q+\frac{i\sigma}{2}U_++\frac{\sigma^2}{4}U_+^3.\]
\item For $n\in\N$, there is $C=C(n)>0$ such that
\[\left\|e^{-\frac{i\sigma}{4}U_+^2}f\right\|_{H^n(S\M)}\le C\left(1+|\sigma|\right)^{2n}\|f\|_{H^{3n}(S\M)}.\]
\item For $n\in\N$, the map $\sigma\mapsto e^{-\frac{i\sigma}{4}U_+^2}f$ is $C^1$ from $\R$ to $H^n(S\M)$, and
\[\partial_\sigma e^{-\frac{i\sigma}{4}U_+^2}f=-\frac i4U_+^2e^{-\frac{i\sigma}{4}U_+^2}f=-\frac i4e^{-\frac{i\sigma}{4}U_+^2}U_+^2f.\]
\item For $s\in\R$,
\[\left(e^{-\frac{i\sigma}{4}U_+^2}f\right)\circ a_s=e^{-\frac{i\sigma e^{2s}}{4}U_+^2}\left(f\circ a_s\right).\]
\end{enumerate}
\end{lemma}

\begin{proof}
Since $|F_\sigma|=1$ and $F_{\sigma_1}F_{\sigma_2}=F_{\sigma_1+\sigma_2}$, the operators $e^{-\frac{i\sigma}{4}U_+^2}$ form a unitary group, which commutes with $R_{n_s}=e^{-is(iU_+)}$. Since $U_+1=0$ and $F_\sigma(0)=1$, we have $e^{-\frac{i\sigma}{4}U_+^2}1=1$. Hence,
\[\int_{S\M}e^{-\frac{i\sigma}{4}U_+^2}f\,d\Vol_L=\left\langle e^{-\frac{i\sigma}{4}U_+^2}f,1\right\rangle_{L^2(S\M)}=\left\langle f,e^{\frac{i\sigma}{4}U_+^2}1\right\rangle_{L^2(S\M)}=\int_{S\M}f\,d\Vol_L.\]
(i) We first show that for $F\in\Sc(\R)$, $F(iU_+)f\in C^\infty(S\M)$ and
\begin{align}
QF\left(iU_+\right)f&=F\left(iU_+\right)Qf+iF'\left(iU_+\right)U_+f,\nonumber\\
U_+F\left(iU_+\right)f&=F\left(iU_+\right)U_+f,\nonumber\\
U_-F\left(iU_+\right)f&=F\left(iU_+\right)U_-f-2iF'\left(iU_+\right)Qf+F''\left(iU_+\right)U_+f.\label{eq:multcomm}
\end{align}
Let $X\in\{Q,U_+,U_-\}$. By \eqref{eq:RAd},
\[X\left(f\circ n_s\right)=\left(\left(\Ad_{n_{-s}}X\right)f\right)\circ n_s.\]
By \eqref{eq:Adn}, $\Ad_{n_{-s}}X$ is a polynomial of degree at most $2$ in $s$ with coefficients in $\{Q,U_+,U_-\}$. Since $\check F$ decays rapidly, \eqref{eq:FiU} shows that $F(iU_+)f\in C^\infty(S\M)$ and
\[XF\left(iU_+\right)f=\int_\R\check F(s)\left(\left(\Ad_{n_{-s}}X\right)f\right)\circ n_s\,ds.\]
Notice that
\[\int_\R s\check F(s)e^{-is\eta}\,ds=iF'(\eta)\quad\text{and}\quad\int_\R s^2\check F(s)e^{-is\eta}\,ds=-F''(\eta).\]
So by \eqref{eq:FiU}, the terms of degree $1$ and $2$ in $s$ give $iF'(iU_+)$ and $-F''(iU_+)$, respectively. This proves \eqref{eq:multcomm}.

Now let $c>0$ and $F_{\sigma,c}(\eta)=e^{(\frac{i\sigma}{4}-c)\eta^2}$. Then $F_{\sigma,c}\in\Sc(\R)$,
\[F_{\sigma,c}'(\eta)=\left(\frac{i\sigma}{2}-2c\right)\eta F_{\sigma,c}(\eta),\quad\text{and}\quad F_{\sigma,c}''(\eta)=\left[\left(\frac{i\sigma}{2}-2c\right)+\left(\frac{i\sigma}{2}-2c\right)^2\eta^2\right]F_{\sigma,c}(\eta).\]
Since $(iU_+)u=iU_+u$ and $(iU_+)^2u=-U_+^2u$ for $u\in C^\infty(S\M)$, \eqref{eq:multcomm} becomes
\begin{align*}
QF_{\sigma,c}\left(iU_+\right)f&=F_{\sigma,c}\left(iU_+\right)\left[Q-\left(\frac{i\sigma}{2}-2c\right)U_+^2\right]f,\\
U_-F_{\sigma,c}\left(iU_+\right)f&=F_{\sigma,c}\left(iU_+\right)\left[U_-+(i\sigma-4c)U_+Q+\left(\frac{i\sigma}{2}-2c\right)U_+-\left(\frac{i\sigma}{2}-2c\right)^2U_+^3\right]f.
\end{align*}
As $c\to0$, the spectral theorem gives
\[F_{\sigma,c}\left(iU_+\right)f\to e^{-\frac{i\sigma}{4}U_+^2}f
\quad\text{in }L^2(S\M).\]
For each $X\in\{Q,U_+,U_-\}$, the formulas above show that $XF_{\sigma,c}(iU_+)f$ also converges in $L^2(S\M)$ as $c\to0$. We regard $X$ as the closure of its action on $C^\infty(S\M)$, which is the skew-adjoint generator of the corresponding one-parameter unitary group. By the closedness of $X$,
\[Xe^{-\frac{i\sigma}{4}U_+^2}f=e^{-\frac{i\sigma}{4}U_+^2}X_\sigma f.\]
Applying this identity successively to the smooth functions $X_{1,\sigma}\cdots X_{l,\sigma}f$, we obtain
\begin{equation}\label{eq:wordconj}
X_1\cdots X_le^{-\frac{i\sigma}{4}U_+^2}f=e^{-\frac{i\sigma}{4}U_+^2}X_{1,\sigma}\cdots X_{l,\sigma}f\quad\text{for }X_1,\dots,X_l\in\left\{Q,U_+,U_-\right\}.
\end{equation}
In particular, $e^{-\frac{i\sigma}{4}U_+^2}f\in H^l(S\M)$ for every $l\in\N$, so $e^{-\frac{i\sigma}{4}U_+^2}f\in C^\infty(S\M)$.

(ii) Each $X_\sigma$ is a sum of at most four products of at most three vector fields in $\{Q,U_+,U_-\}$, with coefficients bounded by $(1+|\sigma|)^2$. Hence, by \eqref{eq:wordconj} and the unitarity of $e^{-\frac{i\sigma}{4}U_+^2}$,
\[\left\|X_1\cdots X_le^{-\frac{i\sigma}{4}U_+^2}f\right\|_{L^2(S\M)}\le4^l\left(1+|\sigma|\right)^{2l}\|f\|_{H^{3l}(S\M)}.\]
Summing over $X_1,\dots,X_l\in\{Q,U_+,U_-\}$ and $0\le l\le n$ gives (ii).

(iii) Since $f\in C^\infty(S\M)$ is in the domain of $(iU_+)^2$, the spectral theorem gives
\[\partial_\sigma e^{-\frac{i\sigma}{4}U_+^2}f=\frac i4\left(iU_+\right)^2e^{-\frac{i\sigma}{4}U_+^2}f=-\frac i4e^{-\frac{i\sigma}{4}U_+^2}U_+^2f\quad\text{in }L^2(S\M).\]
By \eqref{eq:wordconj}, $X_1\cdots X_l(e^{-\frac{i\sigma}{4}U_+^2}f)$ is a finite sum of terms
\[p\left(\sigma\right)e^{-\frac{i\sigma}{4}U_+^2}g,\]
where $p$ is a polynomial and $g\in C^\infty(S\M)$. Each such term is $C^1$ in $\sigma$ with values in $L^2(S\M)$, by the spectral theorem. Hence, $\sigma\mapsto e^{-\frac{i\sigma}{4}U_+^2}f$ is $C^1$ from $\R$ to $H^n(S\M)$. Finally, by (i),
\[U_+^2e^{-\frac{i\sigma}{4}U_+^2}f=e^{-\frac{i\sigma}{4}U_+^2}U_+^2f.\]
This proves (iii).

(iv) By $n_{s_1}a_s=a_sn_{s_1e^{-s}}$ in Section \ref{sec:Lie}, we have
\[R_{a_s}R_{n_{s_1}}R_{a_s}^{-1}=R_{n_{s_1e^s}}\quad\text{for }s_1\in\R.\]
So
\[R_{a_s}U_+R_{a_s}^{-1}=e^sU_+,\quad\text{and}\quad R_{a_s}F(iU_+)R_{a_s}^{-1}=F(ie^sU_+)\]
for bounded Borel functions $F$ on $\R$. Taking $F=F_\sigma$, we obtain
\[R_{a_s}e^{-\frac{i\sigma}{4}U_+^2}=e^{-\frac{i\sigma e^{2s}}{4}U_+^2}R_{a_s},\]
which is (iv).
\end{proof}

\subsection{Egorov's theorem with the Schr\"odinger group}\label{sec:exactSchr}
For $s\le0$, recall that
\[\sigma(s)=h_j\left(e^{2|s|}-1\right).\]
Then $\sigma(s)\ge0$ is decreasing in $s$,
\[\sigma'(s)=-2h_je^{-2s},\quad\text{and}\quad e^{2|s|}=1+\frac{\sigma(s)}{h_j}.\]
For each fixed $T\in\R$ and sufficiently large $j$, we have
\[\sigma(-T_E-T)=e^{2T}-h_j.\]
Thus, $\sigma(t)=O(1)$ when $t=-T_E+O(1)$.

\begin{thm}[Egorov's theorem with the Schr\"odinger group]\label{thm:exactSchr}
Let $f\in C^\infty(S\M)$ and $t\le0$. Then \eqref{eq:exactSchr} holds, that is,
\[\Wc_j\left(f\circ a_t\right)-\Wc_j\left(e^{-\frac{i\sigma(t)}{4}U_+^2}f\right)=\frac{ih_j}{2}\int_0^t\Wc_j\left(\left(Q^2e^{-\frac{i(\sigma(t)-\sigma(s))}{4}U_+^2}f\right)\circ a_s\right)\,ds.\]
\end{thm}

\begin{proof}
For $t\le s\le0$, let
\[u_s=e^{-\frac{i(\sigma(t)-\sigma(s))}{4}U_+^2}f.\]
Then
\[u_t=f\quad\text{and}\quad u_0=e^{-\frac{i\sigma(t)}{4}U_+^2}f.\]
Since $\sigma'(s)=-2h_je^{-2s}$, (iii) of Lemma \ref{lemma:Schr} gives
\[\partial_su_s=-\frac{ih_j}{2}e^{-2s}U_+^2u_s.\]
By (iii) of Lemma \ref{lemma:Schr} and \eqref{eq:Sobolev}, the map $s\mapsto u_s\circ a_s$ is $C^1$ with values in $C^2(S\M)$. Using \eqref{eq:transport} and Theorem \ref{thm:QU}, we obtain
\begin{align*}
&\quad\frac{d}{ds}\Wc_j\left(u_s\circ a_s\right)\\
&=\Wc_j\left(\left(\partial_su_s\right)\circ a_s\right)+\Wc_j\left(Q\left(u_s\circ a_s\right)\right)\\
&=-\frac{ih_j}{2}e^{-2s}\Wc_j\left(\left(U_+^2u_s\right)\circ a_s\right)+\frac{ih_j}{2}\Wc_j\left(\left(Q^2u_s+e^{-2s}U_+^2u_s\right)\circ a_s\right)\\
&=\frac{ih_j}{2}\Wc_j\left(\left(Q^2u_s\right)\circ a_s\right).
\end{align*}
Integrating from $t$ to $0$ gives
\[\Wc_j(u_0)-\Wc_j\left(u_t\circ a_t\right)=\frac{ih_j}{2}\int_t^0\Wc_j\left(\left(Q^2u_s\right)\circ a_s\right)\,ds.\]
Since $u_t=f$ and $u_0=e^{-\frac{i\sigma(t)}{4}U_+^2}f$, reversing the integration limits gives \eqref{eq:exactSchr}.
\end{proof}

We next bound the integral in \eqref{eq:exactSchr} by applying Theorem \ref{thm:exactSchr} twice. Recall from (i) of Lemma \ref{lemma:Schr} that
\[Q^2e^{-\frac{i\sigma}{4}U_+^2}=e^{-\frac{i\sigma}{4}U_+^2}Q_\sigma^2\]
on $C^\infty(S\M)$, where
\[Q_\sigma^2=Q^2-\frac{i\sigma}{2}\left(QU_+^2+U_+^2Q\right)-\frac{\sigma^2}{4}U_+^4.\]
If $\sigma\ge0$, then $Q_\sigma^2$ is a sum of four products of at most four vector fields in $\{Q,U_+\}$, with coefficients bounded by $(1+\sigma)^2$. So for $f\in C^\infty(S\M)$ and $n\in\N$,
\begin{equation}\label{eq:Qzeta}
\left\|Q_\sigma^2f\right\|_{H^n(S\M)}\le4(1+\sigma)^2\|f\|_{H^{n+4}(S\M)}.
\end{equation}
Moreover, by Lemma \ref{lemma:apriori} and \eqref{eq:Sobolev}, for $f\in C^\infty(S\M)$ and $\sigma\ge0$,
\begin{align}
\left|\Wc_j\left(e^{-\frac{i\sigma}{4}U_+^2}f\right)\right|&\le5\left\|e^{-\frac{i\sigma}{4}U_+^2}f\right\|_{C^2(S\M)}\nonumber\\
&\le C\left\|e^{-\frac{i\sigma}{4}U_+^2}f\right\|_{H^4(S\M)}\nonumber\\
&\le C(1+\sigma)^8\|f\|_{H^{12}(S\M)},\label{eq:aprioriSchr}
\end{align}
where we used (ii) of Lemma \ref{lemma:Schr} with $n=4$. Here, $C=C(\M)>0$.

\begin{thm}[Renormalized Egorov estimate]\label{thm:renorm}
There is $C=C(\M)>0$ such that for $f\in C^\infty(S\M)$ and $t\le0$,
\[\left|\Wc_j\left(f\circ a_t\right)-\Wc_j\left(e^{-\frac{i\sigma(t)}{4}U_+^2}f\right)\right|\le Ch_j\left(1+\sigma(t)\right)^{13}\|f\|_{H^{20}(S\M)}.\]
In particular, if $0\le\ve<\frac{1}{13}$, then for $-(1+\ve)T_E\le t\le0$,
\[\left|\Wc_j\left(f\circ a_t\right)-\Wc_j\left(e^{-\frac{i\sigma(t)}{4}U_+^2}f\right)\right|\le2^{13}Ch_j^{1-13\ve}\|f\|_{H^{20}(S\M)}.\]
\end{thm}

So the remainder is $O_f(h_j)$ as long as $\sigma(t)$ is bounded, i.e., for $t\ge-T_E-T$ with $T$ fixed.

\begin{proof}
Let $f\in C^\infty(S\M)$ and $t\le0$. Write $\sigma=\sigma(t)$. For $t\le s\le0$, let
\[u_s=e^{-\frac{i(\sigma-\sigma(s))}{4}U_+^2}f.\]
By Theorem \ref{thm:exactSchr},
\[\Wc_j\left(f\circ a_t\right)-\Wc_j\left(e^{-\frac{i\sigma}{4}U_+^2}f\right)=\frac{ih_j}{2}\int_0^t\Wc_j\left(\left(Q^2u_s\right)\circ a_s\right)\,ds.\]
We apply Theorem \ref{thm:exactSchr} again, to $Q^2u_s$ at time $s$:
\[\Wc_j\left(\left(Q^2u_s\right)\circ a_s\right)-\Wc_j\left(e^{-\frac{i\sigma(s)}{4}U_+^2}Q^2u_s\right)=\frac{ih_j}{2}\int_0^s\Wc_j\left(\left(Q^2e^{-\frac{i(\sigma(s)-\sigma(v))}{4}U_+^2}Q^2u_s\right)\circ a_v\right)\,dv.\]
Hence, the difference $\Wc_j(f\circ a_t)-\Wc_j(e^{-\frac{i\sigma}{4}U_+^2}f)$ is $I_1+I_2$, where
\begin{align*}
I_1&=\frac{ih_j}{2}\int_0^t\Wc_j\left(e^{-\frac{i\sigma(s)}{4}U_+^2}Q^2u_s\right)\,ds,\\
I_2&=-\frac{h_j^2}{4}\int_0^t\int_0^s\Wc_j\left(\left(Q^2e^{-\frac{i(\sigma(s)-\sigma(v))}{4}U_+^2}Q^2u_s\right)\circ a_v\right)\,dv\,ds.
\end{align*}
We bound $I_1$ and $I_2$ separately.

We first bound $I_2$. By (i) of Lemma \ref{lemma:Schr},
\[Q^2e^{-\frac{i(\sigma(s)-\sigma(v))}{4}U_+^2}Q^2u_s=e^{-\frac{i(\sigma-\sigma(v))}{4}U_+^2}Q_{\sigma-\sigma(v)}^2Q_{\sigma-\sigma(s)}^2f.\]
For $t\le s\le v\le0$, we have $0\le\sigma-\sigma(v),\sigma-\sigma(s)\le\sigma$. So the Sobolev embedding, (ii) of Lemma \ref{lemma:Schr}, and \eqref{eq:Qzeta} give
\begin{align*}
\left\|e^{-\frac{i(\sigma-\sigma(v))}{4}U_+^2}Q_{\sigma-\sigma(v)}^2Q_{\sigma-\sigma(s)}^2f\right\|_{C^2(S\M)}&\le C(1+\sigma)^8\left\|Q_{\sigma-\sigma(v)}^2Q_{\sigma-\sigma(s)}^2f\right\|_{H^{12}(S\M)}\\
&\le C(1+\sigma)^{12}\|f\|_{H^{20}(S\M)}.
\end{align*}
Thus, by \eqref{eq:aprioriflow}, the integrand in $I_2$ is bounded by $Ce^{2|v|}(1+\sigma)^{12}\|f\|_{H^{20}(S\M)}$. Since
\[\int_s^0e^{2|v|}\,dv\le\frac12e^{2|s|}\quad\text{and}\quad\int_t^0e^{2|s|}\,ds\le\frac12e^{2|t|},\]
we obtain
\[\left|I_2\right|\le Ch_j^2e^{2|t|}(1+\sigma)^{12}\|f\|_{H^{20}(S\M)}\le Ch_j(1+\sigma)^{13}\|f\|_{H^{20}(S\M)},\]
where the last inequality follows from $h_je^{2|t|}=h_j+\sigma\le1+\sigma$.

We next bound $I_1$. By (i) of Lemma \ref{lemma:Schr} and $[Q,U_+^2]=2U_+^2$,
\begin{align*}
e^{-\frac{i\sigma(s)}{4}U_+^2}Q^2u_s&=\left(Q+\frac{i\sigma(s)}{2}U_+^2\right)^2e^{-\frac{i\sigma}{4}U_+^2}f\\
&=\left(Q^2+i\sigma(s)QU_+^2-i\sigma(s)U_+^2-\frac{\sigma(s)^2}{4}U_+^4\right)e^{-\frac{i\sigma}{4}U_+^2}f.
\end{align*}
Write $g=e^{-\frac{i\sigma}{4}U_+^2}f$. Applying \eqref{eq:transport} to $Qg$ and $U_+^2g$, respectively, gives
\[\Wc_j\left(Q^2g\right)=\frac{ih_j}{2}\Wc_j\left(\left(Q^2+U_+^2\right)Qg\right)\]
and
\[\Wc_j\left(QU_+^2g\right)=\frac{ih_j}{2}\Wc_j\left(\left(Q^2+U_+^2\right)U_+^2g\right).\]
We commute the derivatives through the Schr\"odinger group by (i) of Lemma \ref{lemma:Schr}. Each resulting factor $Q_\sigma$ costs a factor $1+\sigma$ and at most two derivatives, as in \eqref{eq:Qzeta}. The first identity involves at most three factors $Q_\sigma$, and the second at most two. Then \eqref{eq:aprioriSchr} gives
\[\left|\Wc_j\left(Q^2g\right)\right|\le Ch_j(1+\sigma)^{11}\|f\|_{H^{18}(S\M)}\quad\text{and}\quad\left|\Wc_j\left(QU_+^2g\right)\right|\le Ch_j(1+\sigma)^{10}\|f\|_{H^{18}(S\M)}.\]
The remaining terms are bounded directly by \eqref{eq:aprioriSchr}, since $U_+$ commutes with the Schr\"odinger group:
\[\left|\Wc_j\left(U_+^2g\right)\right|+\left|\Wc_j\left(U_+^4g\right)\right|\le C(1+\sigma)^8\|f\|_{H^{18}(S\M)}.\]
Since $\sigma(s)\le\sigma$, we obtain
\[\left|\Wc_j\left(e^{-\frac{i\sigma(s)}{4}U_+^2}Q^2u_s\right)\right|\le C\left[h_j(1+\sigma)^{11}+\left(\sigma(s)+\sigma(s)^2\right)(1+\sigma)^8\right]\|f\|_{H^{18}(S\M)}.\]
Moreover, $\sigma'(s)=-2(h_j+\sigma(s))$ and $\sigma(s)\ge0$. So
\[\int_t^0\sigma(s)\,ds\le\frac{\sigma}{2}\quad\text{and}\quad\int_t^0\sigma(s)^2\,ds\le\frac{\sigma^2}{4}.\]
Since $h_j|t|\le\frac12h_je^{2|t|}\le1+\sigma$, we obtain
\[\left|I_1\right|\le C\left[h_j^2|t|(1+\sigma)^{11}+h_j\left(\sigma+\sigma^2\right)(1+\sigma)^8\right]\|f\|_{H^{18}(S\M)}\le Ch_j(1+\sigma)^{12}\|f\|_{H^{18}(S\M)}.\]
The bounds of $I_1$ and $I_2$ prove the first estimate. Finally, if $-(1+\ve)T_E\le t\le0$, then $1+\sigma(t)\le2h_j^{-\ve}$, which gives the second estimate.
\end{proof}

We are now ready to prove (ii) of Theorem \ref{thm:Egorov}.

\begin{proof}[Proof of (ii) of Theorem \ref{thm:Egorov}]
Let $f\in C^\infty(S\M)$ and $-T_E\le t\le0$. Then
\[\sigma(t)=h_j\left(e^{2|t|}-1\right)\le1.\]
By (iii) of Lemma \ref{lemma:Schr},
\begin{align*}
\Wc_j\left(e^{-\frac{i\sigma(t)}{4}U_+^2}f\right)-\Wc_j(f)&=\int_0^{\sigma(t)}\Wc_j\left(\partial_\sigma e^{-\frac{i\sigma}{4}U_+^2}f\right)\,d\sigma\\
&=\int_0^{\sigma(t)}\Wc_j\left(-\frac i4e^{-\frac{i\sigma}{4}U_+^2}U_+^2f\right)\,d\sigma.
\end{align*}
Hence, by \eqref{eq:aprioriSchr}, since $\sigma(t)\le1$,
\begin{align*}
\left|\Wc_j\left(e^{-\frac{i\sigma(t)}{4}U_+^2}f\right)-\Wc_j(f)\right|&\le\frac14\int_0^{\sigma(t)}\left|\Wc_j\left(e^{-\frac{i\sigma}{4}U_+^2}U_+^2f\right)\right|\,d\sigma\\
&\le C\int_0^{\sigma(t)}(1+\sigma)^8\left\|U_+^2f\right\|_{H^{12}(S\M)}\,d\sigma\\
&\le C\sigma(t)\|f\|_{H^{14}(S\M)}.
\end{align*}
Since $\sigma(t)\le1$, Theorem \ref{thm:renorm} gives
\[\left|\Wc_j\left(f\circ a_t\right)-\Wc_j\left(e^{-\frac{i\sigma(t)}{4}U_+^2}f\right)\right|\le Ch_j\|f\|_{H^{20}(S\M)}.\]
Therefore,
\begin{align*}
&\quad\left|\Wc_j\left(f\circ a_t\right)-\Wc_j(f)\right|\\
&\le\left|\Wc_j\left(f\circ a_t\right)-\Wc_j\left(e^{-\frac{i\sigma(t)}{4}U_+^2}f\right)\right|+\left|\Wc_j\left(e^{-\frac{i\sigma(t)}{4}U_+^2}f\right)-\Wc_j(f)\right|\\
&\le Ch_j\|f\|_{H^{20}(S\M)}+C\sigma(t)\|f\|_{H^{14}(S\M)}\\
&\le C\left(h_j+\sigma(t)\right)\|f\|_{H^{20}(S\M)}\\
&=Ch_je^{2|t|}\|f\|_{H^{20}(S\M)}.
\end{align*}
This proves (ii) of Theorem \ref{thm:Egorov}.
\end{proof}

\section{Egorov's theorem at and beyond the Ehrenfest time}\label{sec:beyond}
In this section, we prove Theorem \ref{thm:equiv}. Let $J\subset\N$ be an infinite set such that $\mu_j>\frac54$ for $j\in J$.

\subsection{Semiclassical measures}\label{sec:scm}
By Lemma \ref{lemma:positive}, for each $j$ with $\mu_j>\frac54$, there is a probability measure $\nu_j$ on $S\M$ such that
\[\left|\Wc_j(f)-\int_{S\M}f\,d\nu_j\right|\le Ch_j^{\frac12}\left(1+|\log h_j|\right)\|f\|_{C^2(S\M)}.\]
Since $S\M$ is compact, every sequence $\{\nu_j\}$ with $j\to\infty$ has a weakly convergent subsequence. The preceding estimate shows that the corresponding subsequence of $\Wc_j$ converges to the same limit on smooth functions. Thus, every sequence of eigenfunctions with $j\to\infty$ has a subsequence which induces a semiclassical measure, and this measure is a probability measure on $S\M$.

Suppose that $\{\vp_j\}_{j\in J}$ induces a semiclassical measure $\nu$. By \eqref{eq:transport} and Lemma \ref{lemma:apriori},
\[\int_{S\M}Qf\,d\nu=\lim_{j\in J,j\to\infty}\frac{ih_j}{2}\Wc_j\left(\left(Q^2+U_+^2\right)f\right)=0\quad\text{for }f\in C^\infty(S\M).\]
So $\nu$ is invariant under the geodesic flow $a_t$.

We next show that the only finite signed measures on $S\M$ which are annihilated by $U_+^2$ are the multiples of the Liouville measure. This uses Furstenberg's unique ergodicity of the horocycle flow \cite{Fu}.

\begin{lemma}\label{lemma:horo}
Let $\nu$ be a finite signed measure on $S\M$. Suppose that
\[\int_{S\M}U_+^2f\,d\nu=0\quad\text{for }f\in C^\infty(S\M).\]
Then $\nu=\nu(S\M)\,d\Vol_L$.
\end{lemma}

\begin{proof}
Let $f\in C^\infty(S\M)$. For $s\in\R$, let
\[F(s)=\int_{S\M}f\circ n_s\,d\nu,\quad\text{where }n_s=e^{sU_+}.\]
Since $U_+$ commutes with the right translation by $n_s$, we have
\[F''(s)=\int_{S\M}U_+^2\left(f\circ n_s\right)\,d\nu=0.\]
So $F$ is affine in $s\in\R$. Denote by $|\nu|$ the total variation of $\nu$. Then
\[|F(s)|\le|\nu|(S\M)\|f\|_{L^\infty(S\M)}\quad\text{for }s\in\R.\]
Hence, $F$ is constant, i.e., $\nu$ is invariant under the stable horocycle flow. Since the flow consists of homeomorphisms of $S\M$, $|\nu|$ is also invariant, and so are the positive and negative parts $\frac12(|\nu|\pm\nu)$ of $\nu$. By Furstenberg \cite{Fu}, the stable horocycle flow on the compact quotient $S\M=\Gamma\bs\PSU(1,1)$ is uniquely ergodic. So $\frac12(|\nu|\pm\nu)$ are multiples of $d\Vol_L$. Therefore, $\nu=\nu(S\M)\,d\Vol_L$.
\end{proof}

\subsection{Proof of Theorem \ref{thm:equiv}}
We first compute the limits in Egorov's theorem at the times $-T_E+O(1)$.

\begin{lemma}\label{lemma:blowup}
Suppose that $\{\vp_j\}_{j\in J}$ induces a semiclassical measure $\nu$. Let $t_j\le0$, $j\in J$, satisfy
\[\lim_{j\in J,j\to\infty}h_j\left(e^{2|t_j|}-1\right)=\sigma_0\in[0,\infty).\]
Then for $f\in C^\infty(S\M)$,
\[\lim_{j\in J,j\to\infty}\Wc_j\left(f\circ a_{t_j}\right)=\int_{S\M}e^{-\frac{i\sigma_0}{4}U_+^2}f\,d\nu.\]
\end{lemma}

\begin{proof}
Write
\[\sigma_j=h_j\left(e^{2|t_j|}-1\right).\]
Since $\sigma_j$ converges, there is $C_0>0$ such that $\sigma_j\le C_0$ for $j\in J$. So by Theorem \ref{thm:renorm},
\[\left|\Wc_j\left(f\circ a_{t_j}\right)-\Wc_j\left(e^{-\frac{i\sigma_j}{4}U_+^2}f\right)\right|\le C\left(1+C_0\right)^{13}h_j\|f\|_{H^{20}(S\M)}\to0.\]
Moreover, by Lemma \ref{lemma:apriori} and \eqref{eq:Sobolev},
\begin{align*}
\left|\Wc_j\left(e^{-\frac{i\sigma_j}{4}U_+^2}f\right)-\Wc_j\left(e^{-\frac{i\sigma_0}{4}U_+^2}f\right)\right|&\le C\left\|e^{-\frac{i\sigma_j}{4}U_+^2}f-e^{-\frac{i\sigma_0}{4}U_+^2}f\right\|_{C^2(S\M)}\\
&\le C\left\|e^{-\frac{i\sigma_j}{4}U_+^2}f-e^{-\frac{i\sigma_0}{4}U_+^2}f\right\|_{H^4(S\M)}.
\end{align*}
We use (iii) of Lemma \ref{lemma:Schr} to conclude that
\[\left\|e^{-\frac{i\sigma_j}{4}U_+^2}f-e^{-\frac{i\sigma_0}{4}U_+^2}f\right\|_{H^4(S\M)}\to0,\]
since $\sigma_j\to\sigma_0$.

Finally, $e^{-\frac{i\sigma_0}{4}U_+^2}f\in C^\infty(S\M)$ by (i) of Lemma \ref{lemma:Schr}. Therefore,
\[\lim_{j\in J,j\to\infty}\Wc_j\left(e^{-\frac{i\sigma_0}{4}U_+^2}f\right)=\int_{S\M}e^{-\frac{i\sigma_0}{4}U_+^2}f\,d\nu.\]
This completes the proof.
\end{proof}

\begin{lemma}\label{lemma:profile}
Let $\nu$ be a finite signed measure on $S\M$ which is invariant under the geodesic flow. Suppose that there is $\sigma_0>0$ such that
\[\int_{S\M}e^{-\frac{i\sigma_0}{4}U_+^2}f\,d\nu=\int_{S\M}f\,d\nu\quad\text{for }f\in C^\infty(S\M).\]
Then $\nu=\nu(S\M)\,d\Vol_L$.
\end{lemma}

\begin{proof}
Let $f\in C^\infty(S\M)$. For $\sigma\in\R$, let
\[F(\sigma)=\int_{S\M}e^{-\frac{i\sigma}{4}U_+^2}f\,d\nu.\]
Since $\nu$ is finite, by \eqref{eq:Sobolev},
\[\left|\int_{S\M}u\,d\nu\right|\le|\nu|(S\M)\|u\|_{L^\infty(S\M)}\le C\|u\|_{H^4(S\M)}\quad\text{for }u\in C^\infty(S\M).\]
So by (iii) of Lemma \ref{lemma:Schr}, $F$ is $C^1$ on $\R$, and
\[F'(0)=-\frac i4\int_{S\M}U_+^2f\,d\nu.\]
Let $s\in\R$. Since $\nu$ is invariant under $a_s$, (iv) of Lemma \ref{lemma:Schr} and the assumption give
\begin{align*}
F\left(\sigma_0e^{-2s}\right)&=\int_{S\M}\left(e^{-\frac{i\sigma_0e^{-2s}}{4}U_+^2}f\right)\circ a_s\,d\nu\\
&=\int_{S\M}e^{-\frac{i\sigma_0}{4}U_+^2}\left(f\circ a_s\right)\,d\nu\\
&=\int_{S\M}f\circ a_s\,d\nu\\
&=\int_{S\M}f\,d\nu.
\end{align*}
Hence, $F$ is constant on $(0,\infty)$. Since $F$ is $C^1$ on $\R$, we have $F'(0)=0$, i.e.,
\[\int_{S\M}U_+^2f\,d\nu=0.\]
By Lemma \ref{lemma:horo}, $\nu=\nu(S\M)\,d\Vol_L$.
\end{proof}

We are now ready to prove Theorem \ref{thm:equiv}.

\begin{proof}[Proof of Theorem \ref{thm:equiv}]
Let $f\in C^\infty(S\M)$. We first prove \eqref{eq:defect}. Let $T\in\R$ and $t_j=-T_E-T$. Then $t_j\le0$ for large $j$, and
\[h_j\left(e^{2|t_j|}-1\right)=e^{2T}-h_j\to e^{2T}.\]
Since $\Wc_j(f)\to\int_{S\M}f\,d\nu$ as $j\in J$, $j\to\infty$, Lemma \ref{lemma:blowup} with $\sigma_0=e^{2T}$ gives \eqref{eq:defect}.

Next, we show that the limits in \eqref{eq:defect} determine $\nu$. Recall from Section \ref{sec:scm} that $\nu$ is a probability measure on $S\M$ which is invariant under the geodesic flow. Let $T\in\R$. Suppose that $\nu'$ is another such measure with
\[\int_{S\M}\left(e^{-\frac{ie^{2T}}{4}U_+^2}f-f\right)d\nu'=\int_{S\M}\left(e^{-\frac{ie^{2T}}{4}U_+^2}f-f\right)d\nu\quad\text{for }f\in C^\infty(S\M).\]
Then $\nu-\nu'$ is a finite signed measure on $S\M$ which is invariant under the geodesic flow, and it satisfies the assumption of Lemma \ref{lemma:profile} with $\sigma_0=e^{2T}$. Since $(\nu-\nu')(S\M)=0$, Lemma \ref{lemma:profile} gives $\nu'=\nu$.

(i)$\Rightarrow$(iii). We first show that for each $T\ge0$,
\begin{equation}\label{eq:uniformT}
\lim_{j\in J,j\to\infty}\sup_{-T_E-T\le t\le0}\left|\Wc_j\left(f\circ a_t\right)-\Wc_j(f)\right|=0.
\end{equation}
Let $T\ge0$ be fixed. For $-T_E-T\le t\le0$, we have $0\le\sigma(t)\le e^{2T}$. So by Theorem \ref{thm:renorm},
\[\sup_{-T_E-T\le t\le0}\left|\Wc_j\left(f\circ a_t\right)-\Wc_j\left(e^{-\frac{i\sigma(t)}{4}U_+^2}f\right)\right|\le C\left(1+e^{2T}\right)^{13}h_j\|f\|_{H^{20}(S\M)}\to0.\]
It remains to show that
\[\lim_{j\in J,j\to\infty}\sup_{0\le\sigma\le e^{2T}}\left|\Wc_j\left(e^{-\frac{i\sigma}{4}U_+^2}f\right)-\Wc_j(f)\right|=0.\]
By (i), $\Wc_j(g)\to\int_{S\M}g\,d\Vol_L$ as $j\in J$, $j\to\infty$, for $g\in C^\infty(S\M)$. Since the Schr\"odinger group preserves $d\Vol_L$, the difference in the supremum tends to $0$ for each $\sigma$. Moreover, by (iii) of Lemma \ref{lemma:Schr} and \eqref{eq:aprioriSchr}, for $0\le\sigma_1\le\sigma_2\le e^{2T}$,
\begin{align*}
\left|\Wc_j\left(e^{-\frac{i\sigma_2}{4}U_+^2}f\right)-\Wc_j\left(e^{-\frac{i\sigma_1}{4}U_+^2}f\right)\right|&\le\frac14\int_{\sigma_1}^{\sigma_2}\left|\Wc_j\left(e^{-\frac{i\sigma}{4}U_+^2}U_+^2f\right)\right|\,d\sigma\\
&\le C\left(1+e^{2T}\right)^8\left(\sigma_2-\sigma_1\right)\|f\|_{H^{14}(S\M)}.
\end{align*}
That is, the functions $\sigma\mapsto\Wc_j(e^{-\frac{i\sigma}{4}U_+^2}f)$ are Lipschitz on $[0,e^{2T}]$, uniformly in $j$. So their convergence is uniform on $[0,e^{2T}]$. This proves \eqref{eq:uniformT}.

By \eqref{eq:uniformT}, there are $j_1<j_2<\cdots$ in $J$ such that
\[\sup_{-T_E-l\le t\le0}\left|\Wc_j\left(f\circ a_t\right)-\Wc_j(f)\right|\le\frac1l\quad\text{for }j\in J\text{ with }j\ge j_l.\]
Let $T_j=l$ for $j\in J$ with $j_l\le j<j_{l+1}$, and $T_j=0$ for $j\in J$ with $j<j_1$. Then $T_j\to\infty$, and (iii) holds.

(iii)$\Rightarrow$(ii). Since $T_j\ge0$, we have $-T_E-T_j\le-T_E\le0$. So (iii) gives (ii).

(ii)$\Rightarrow$(i). By (ii), the limits in \eqref{eq:defect} with $T=0$ vanish. They also vanish if $\nu$ is replaced by $d\Vol_L$, since the Schr\"odinger group preserves $d\Vol_L$ by Lemma \ref{lemma:Schr}. Since these limits determine $\nu$, we obtain $\nu=d\Vol_L$.
\end{proof}

\begin{rmk}
We justify the consequences of Theorem \ref{thm:equiv} stated after it in Section \ref{sec:intro}, where $\{\vp_j\}_{j\in J}$ need not induce a semiclassical measure. Since every subsequence of $\{\vp_j\}_{j\in J}$ has a further subsequence which induces a semiclassical measure, $\{\vp_j\}_{j\in J}$ tends to be equidistributed in the phase space if and only if
\[\lim_{j\in J,j\to\infty}\Wc_j(g)=\int_{S\M}g\,d\Vol_L\quad\text{for }g\in C^\infty(S\M).\]
The proof of (i)$\Rightarrow$(iii) only uses this property. So if $\{\vp_j\}_{j\in J}$ tends to be equidistributed in the phase space, then (iii) holds. Conversely, suppose that Egorov's theorem holds at $-T_E-T$ along $J$ for some fixed $T\in\R$. Then by the proof of (ii)$\Rightarrow$(i) with $T$ in place of $0$, every semiclassical measure induced by the subsequences of $\{\vp_j\}_{j\in J}$ is $d\Vol_L$. That is, $\{\vp_j\}_{j\in J}$ tends to be equidistributed in the phase space.
\end{rmk}

\section{Egorov's theorem with a rate and QUE with a rate}\label{sec:rate}
In this section, we prove Theorem \ref{thm:rate}. We use the renormalized Egorov estimate in Theorem \ref{thm:renorm} and the effective equidistribution of the stable horocycle flow. Throughout this section, let $J\subset\N$ be an infinite set such that $\mu_j>\frac54$ for $j\in J$.

\subsection{Effective equidistribution of the horocycle flow}
For $f\in C^\infty(S\M)$ and $s\in\R$, recall that $f\circ n_s=R_{n_s}f$, where $n_s=e^{sU_+}$. By \eqref{eq:RAd},
\[X\left(f\circ n_s\right)=\left(\left(\Ad_{n_{-s}}X\right)f\right)\circ n_s\quad\text{for }X\in\left\{Q,U_+,U_-\right\}.\]
By \eqref{eq:Adn}, each $\Ad_{n_{-s}}X$, $X\in\{Q,U_+,U_-\}$, is a linear combination of $Q,U_+,U_-$ whose coefficients have absolute values summing to at most $(1+|s|)^2$. Applying this identity successively to products of at most $n$ vector fields, and using the invariance of $d\Vol_L$ under the stable horocycle flow, we obtain
\begin{equation}\label{eq:normns}
\left\|f\circ n_s\right\|_{H^n(S\M)}\le C\left(1+|s|\right)^{2n}\|f\|_{H^n(S\M)},
\end{equation}
where $C=C(n)>0$.

Recall from the proof of Lemma \ref{lemma:horo} that the stable horocycle flow on $S\M$ is uniquely ergodic. That is, for $f\in C^\infty(S\M)$, the averages
\[\frac{1}{2R}\int_{-R}^Rf\circ n_s\,ds\]
along the stable horocycles converge to $\int_{S\M}f\,d\Vol_L$ uniformly on $S\M$ as $R\to\infty$. Burger \cite{Bur} proved that the convergence has a polynomial rate in $R$, and Flaminio-Forni \cite{FF} found the optimal rate, which is determined by the first non-zero eigenvalue $\mu_1$ of $\Delta_\Gamma$. We use their results in the following form.

\begin{thm}[Effective equidistribution of the horocycle flow]\label{thm:Burger}
Let
\[\beta=\frac12-\Rer\frac{\sqrt{1-4\mu_1}}{2}\quad\text{if }\mu_1\ne\frac14,\]
and let $0<\beta<\frac12$ if $\mu_1=\frac14$. Then there is $C=C(\M,\beta)>0$ such that for $f\in C^\infty(S\M)$ and $R\ge1$,
\[\left\|\frac{1}{2R}\int_{-R}^Rf\circ n_s\,ds-\int_{S\M}f\,d\Vol_L\right\|_{L^\infty(S\M)}\le CR^{-\beta}\|f\|_{H^6(S\M)}.\]
\end{thm}

Notice that the average of $f$ at $x\in S\M$ in Theorem \ref{thm:Burger} is along the stable horocycle segment of length $2R$ starting at $xn_{-R}$, and the estimate is uniform in $x$. If $\mu_1<\frac14$, then the rate $R^{-\beta}$ comes from the complementary series representation of $\PSU(1,1)$ generated by the eigenfunctions with eigenvalue $\mu_1$. If $\mu_1>\frac14$, the principal series gives the rate $R^{-\frac12}$, while at the threshold $\mu_1=\frac14$ logarithmic factors may occur. The Sobolev order $6$ in Theorem \ref{thm:Burger} may not be optimal but is sufficient for our purpose; see Flaminio-Forni \cite[Theorem 1.5]{FF} and Ravotti \cite[Theorem 2]{Rav} for more discussion.

We now prove a quantitative version of Lemma \ref{lemma:horo}: a probability measure on $S\M$ which is almost annihilated by $U_+^2$ is close to $d\Vol_L$.

\begin{lemma}\label{lemma:horoRate}
Let $n\in\N$. Then there is $C=C(\M,\beta,n)>0$ such that the following holds. Let $\nu$ be a probability measure on $S\M$ and $0<\kappa\le1$. Suppose that
\[\left|\int_{S\M}U_+^2f\,d\nu\right|\le\kappa\|f\|_{H^n(S\M)}\quad\text{for }f\in C^\infty(S\M).\]
Then
\[\left|\int_{S\M}f\,d\nu-\int_{S\M}f\,d\Vol_L\right|\le C\kappa^{\frac{\beta}{2n+2+\beta}}\|f\|_{H^{\max\{n,6\}}(S\M)}\quad\text{for }f\in C^\infty(S\M).\]
\end{lemma}

\begin{proof}
Let $f\in C^\infty(S\M)$. As in the proof of Lemma \ref{lemma:horo}, write
\[F(s)=\int_{S\M}f\circ n_s\,d\nu\quad\text{for }s\in\R.\]
Then
\[F''(s)=\int_{S\M}U_+^2\left(f\circ n_s\right)\,d\nu.\]
So by the assumption and \eqref{eq:normns},
\[\left|F''(s)\right|\le\kappa\left\|f\circ n_s\right\|_{H^n(S\M)}\le C\kappa\left(1+|s|\right)^{2n}\|f\|_{H^n(S\M)}.\]
Let $R\ge1$. By Taylor's formula, for $0\le s\le R$,
\[\left|F(s)+F(-s)-2F(0)\right|
\le s^2\sup_{|u|\le s}|F''(u)|
\le C\kappa s^2(1+R)^{2n}\|f\|_{H^n(S\M)}.\]
Notice that
\begin{align*}
\int_0^R\left(F(s)+F(-s)\right)\,ds&=\int_0^R\left(\int_{S\M}f\circ n_s\,d\nu\right)\,ds+\int_0^R\left(\int_{S\M}f\circ n_{-s}\,d\nu\right)\,ds\\
&=\int_0^R\left(\int_{S\M}f\circ n_s\,d\nu\right)\,ds+\int_{-R}^0\left(\int_{S\M}f\circ n_s\,d\nu\right)\,ds\\
&=\int_{S\M}\left(\int_{-R}^Rf\circ n_s\,ds\right)\,d\nu.
\end{align*}
Since $F(0)=\int_{S\M}f\,d\nu$,
\begin{align*}
\left|\int_{S\M}f\,d\nu-\int_{S\M}\left(\frac{1}{2R}\int_{-R}^Rf\circ n_s\,ds\right)\,d\nu\right|&=\left|F(0)-\frac{1}{2R}\int_0^R\left(F(s)+F(-s)\right)\,ds\right|\\
&\le\frac{1}{2R}\int_0^R\left|2F(0)-F(s)-F(-s)\right|\,ds\\
&\le C\kappa(1+R)^{2n+2}\|f\|_{H^n(S\M)}.
\end{align*}
On the other hand, since $\nu$ is a probability measure, Theorem \ref{thm:Burger} gives
\[\left|\int_{S\M}\left(\frac{1}{2R}\int_{-R}^Rf\circ n_s\,ds\right)\,d\nu-\int_{S\M}f\,d\Vol_L\right|\le CR^{-\beta}\|f\|_{H^6(S\M)}.\]
Now we choose
\[R=\kappa^{-\frac{1}{2n+2+\beta}}\ge1.\]
Then
\[\kappa(1+R)^{2n+2}\le2^{2n+2}\kappa^{\frac{\beta}{2n+2+\beta}}\quad\text{and}\quad R^{-\beta}=\kappa^{\frac{\beta}{2n+2+\beta}}.\]
Combining the above two estimates completes the proof.
\end{proof}

\subsection{Egorov's theorem with a rate implies QUE with a rate}
We first show that the Egorov estimate at $-T_E$ with a rate in \eqref{eq:Egorovrate} implies that $\Wc_j$ is approximately annihilated by $U_+^2$.

\begin{prop}\label{prop:horoRate}
Suppose that \eqref{eq:Egorovrate} holds with $\delta_2>0$ and $n_2\in\N$. Let $N=\max\{n_2,20\}$. Then there is $C>0$ such that for $j\in J$ and $f\in C^\infty(S\M)$,
\[\left|\Wc_j\left(U_+^2f\right)\right|\le C\left(h_j^{\delta_2}+h_j\right)\|f\|_{H^{N+2}(S\M)}.\]
\end{prop}

\begin{proof}
Since $\sigma(-T_E)=1-h_j$, Theorem \ref{thm:renorm} at $t=-T_E$, together with \eqref{eq:Egorovrate}, gives for $g\in C^\infty(S\M)$,
\begin{align*}
\left|\Wc_j\left(e^{-\frac{i(1-h_j)}{4}U_+^2}g\right)-\Wc_j(g)\right|&\le\left|\Wc_j\left(e^{-\frac{i(1-h_j)}{4}U_+^2}g\right)-\Wc_j\left(g\circ a_{-T_E}\right)\right|+\left|\Wc_j\left(g\circ a_{-T_E}\right)-\Wc_j(g)\right|\\
&\le C\left(h_j+h_j^{\delta_2}\right)\|g\|_{H^N(S\M)}.
\end{align*}
By (iii) of Lemma \ref{lemma:Schr} and \eqref{eq:aprioriSchr},
\[\left|\Wc_j\left(e^{-\frac i4U_+^2}g-e^{-\frac{i(1-h_j)}{4}U_+^2}g\right)\right|\le\frac14\int_{1-h_j}^1\left|\Wc_j\left(e^{-\frac{i\sigma}{4}U_+^2}U_+^2g\right)\right|\,d\sigma\le Ch_j\|g\|_{H^{14}(S\M)}.\]
Since $N\ge20$, combining these estimates gives
\begin{equation}\label{eq:singleSchrRate}
\left|\Wc_j\left(e^{-\frac i4U_+^2}g-g\right)\right|\le C\left(h_j^{\delta_2}+h_j\right)\|g\|_{H^N(S\M)}.
\end{equation}
By $[Q,U_+^2]=2U_+^2$ and (i) of Lemma \ref{lemma:Schr},
\[\left[Q,e^{-\frac i4U_+^2}\right]=-\frac i2U_+^2e^{-\frac i4U_+^2}.\]
Since $U_+^2$ commutes with $e^{-\frac i4U_+^2}$, we obtain
\begin{align*}
U_+^2&=2i\left[Q,e^{-\frac i4U_+^2}\right]-U_+^2e^{-\frac i4U_+^2}+U_+^2\\
&=2i\left(Qe^{-\frac i4U_+^2}-Q\right)-2i\left(e^{-\frac i4U_+^2}Q-Q\right)-\left(e^{-\frac i4U_+^2}U_+^2-U_+^2\right).
\end{align*}
It then follows that
\begin{equation}\label{eq:WjU}
\Wc_j\left(U_+^2f\right)=2i\Wc_j\left(Qe^{-\frac i4U_+^2}f-Qf\right)-2i\Wc_j\left(e^{-\frac i4U_+^2}Qf-Qf\right)-\Wc_j\left(e^{-\frac i4U_+^2}U_+^2f-U_+^2f\right).
\end{equation}
We estimate the three terms on the right-hand side of \eqref{eq:WjU}. For the second term, \eqref{eq:singleSchrRate} with $g=Qf$ gives
\[\left|\Wc_j\left(e^{-\frac i4U_+^2}Qf-Qf\right)\right|\le C\left(h_j^{\delta_2}+h_j\right)\|Qf\|_{H^N(S\M)}\le C\left(h_j^{\delta_2}+h_j\right)\|f\|_{H^{N+1}(S\M)}.\]
For the third term, \eqref{eq:singleSchrRate} with $g=U_+^2f$ gives
\[\left|\Wc_j\left(e^{-\frac i4U_+^2}U_+^2f-U_+^2f\right)\right|\le C\left(h_j^{\delta_2}+h_j\right)\left\|U_+^2f\right\|_{H^N(S\M)}\le C\left(h_j^{\delta_2}+h_j\right)\|f\|_{H^{N+2}(S\M)}.\]
For the first term, we apply Egorov's theorem \eqref{eq:transport} to $e^{-\frac i4U_+^2}f$, and obtain
\[\Wc_j\left(Qe^{-\frac i4U_+^2}f\right)=\frac{ih_j}{2}\Wc_j\left(\left(Q^2+U_+^2\right)e^{-\frac i4U_+^2}f\right).\]
By (i) of Lemma \ref{lemma:Schr},
\[\left(Q^2+U_+^2\right)e^{-\frac i4U_+^2}f=e^{-\frac i4U_+^2}\left(Q_1^2+U_+^2\right)f,\]
where $Q_1=Q-\frac i2U_+^2$. Hence, by \eqref{eq:aprioriSchr},
\begin{align*}
\left|\Wc_j\left(Qe^{-\frac i4U_+^2}f\right)\right|&\le Ch_j\left|\Wc_j\left(e^{-\frac i4U_+^2}\left(Q_1^2+U_+^2\right)f\right)\right|\\
&\le Ch_j\left\|\left(Q_1^2+U_+^2\right)f\right\|_{H^{12}(S\M)}\\
&\le Ch_j\|f\|_{H^{16}(S\M)}.
\end{align*}
Here, the last inequality follows from \eqref{eq:Qzeta} with $\sigma=1$.
Similarly, applying \eqref{eq:transport} directly to $f$, we obtain
\begin{align*}
\left|\Wc_j(Qf)\right|&=\left|\frac{ih_j}{2}\Wc_j\left(\left(Q^2+U_+^2\right)f\right)\right|\\
&\le Ch_j\left\|\left(Q^2+U_+^2\right)f\right\|_{C^2(S\M)}\\
&\le Ch_j\left\|\left(Q^2+U_+^2\right)f\right\|_{H^4(S\M)}\\
&\le Ch_j\|f\|_{H^{6}(S\M)},
\end{align*}
where we used Lemma \ref{lemma:apriori} and the Sobolev embedding.

Therefore, putting the estimates of the right-hand side of \eqref{eq:WjU} together, we have that
\[\left|\Wc_j\left(U_+^2f\right)\right|\le C\left(h_j^{\delta_2}+h_j\right)\|f\|_{H^{N+2}(S\M)},\]
as required.
\end{proof}

\begin{proof}[Proof that (ii) implies (i) of Theorem \ref{thm:rate}]
Suppose that (ii) holds with $\delta_2>0$ and $n_2\in\N$. Let
\[n_1=\max\{n_2,20\}+2,\]
and let $\delta_1$ satisfy
\[0<\delta_1<\frac{\beta}{2n_1+2+\beta}\min\left\{\delta_2,\frac12\right\}.\]
By Proposition \ref{prop:horoRate},
\[\left|\Wc_j\left(U_+^2f\right)\right|\le C\left(h_j^{\delta_2}+h_j\right)\|f\|_{H^{n_1}(S\M)}.\]
Let $\nu_j$ be the probability measure in \eqref{eq:positive}. By Lemma \ref{lemma:positive} and \eqref{eq:Sobolev},
\[\left|\int_{S\M}U_+^2f\,d\nu_j\right|\le C\left(h_j^{\delta_2}+h_j^{\frac12}\left(1+\left|\log h_j\right|\right)\right)\|f\|_{H^{n_1}(S\M)}.\]
For sufficiently large $j$, put
\[\kappa_j=C\left(h_j^{\delta_2}+h_j^{\frac12}\left(1+|\log h_j|\right)\right)\le1.\]
By Lemma \ref{lemma:horoRate} with $n=n_1$,
\[\left|\int_{S\M}f\,d\nu_j-\int_{S\M}f\,d\Vol_L\right|\le C\kappa_j^{\frac{\beta}{2n_1+2+\beta}}\|f\|_{H^{n_1}(S\M)}.\]
Let
\[\alpha=\frac{\beta}{2n_1+2+\beta}\quad\text{and}\quad\theta=\min\left\{\delta_2,\frac12\right\}.\]
For every $\eta>0$, the inequality
\[1+\left|\log h_j\right|\le C h_j^{-\eta}\quad\text{with }C=C(\eta)>0\]
gives
\[\kappa_j\le C_\eta h_j^{\theta-\eta}.\]
Choose $\eta>0$ sufficiently small that
\[\delta_1<\alpha(\theta-\eta).\]
Then
\[\kappa_j^\alpha\le Ch_j^{\delta_1}.\]
Since $\delta_1<\frac12$, we have $h_j^{\frac12}(1+|\log h_j|)\le Ch_j^{\delta_1}$. So \eqref{eq:positive} and the Sobolev embedding give
\[\left|\Wc_j(f)-\int_{S\M}f\,d\nu_j\right|\le Ch_j^{\delta_1}\|f\|_{H^{n_1}(S\M)}.\]
Combining these estimates proves
\[\left|\Wc_j(f)-\int_{S\M}f\,d\Vol_L\right|\le Ch_j^{\delta_1}\|f\|_{H^{n_1}(S\M)}\]
for sufficiently large $j$. The finitely many remaining indices are covered by enlarging $C$. This proves (i) with this $\delta_1$.
\end{proof}

\subsection{QUE with a rate implies Egorov's theorem with a rate}
We show that (i) implies (ii) of Theorem \ref{thm:rate}, and in fact \eqref{eq:Egorovbeyond}, i.e., Egorov's theorem with a polynomial remainder slightly beyond $-T_E$.

\begin{prop}\label{prop:converse}
Suppose that \eqref{eq:QErate} holds with $\delta_1>0$, $n_1\in\N$ with $n_1\ge1$, and $C=C_0$. Let
\[0\le\ve<\min\left\{\frac{1}{13},\frac{\delta_1}{2n_1}\right\}.\]
Then there is $C>0$ such that for $j\in J$ and $f\in C^\infty(S\M)$,
\[\sup_{-(1+\ve)T_E\le t\le0}\left|\Wc_j\left(f\circ a_t\right)-\Wc_j(f)\right|\le C\left(h_j^{1-13\ve}+h_j^{\delta_1-2n_1\ve}\right)\|f\|_{H^{\max\{20,3n_1\}}(S\M)}.\]
\end{prop}

\begin{proof}
Let $j\in J$, $f\in C^\infty(S\M)$, and $-(1+\ve)T_E\le t\le0$. Then
\[0\le\sigma(t)\le h_j^{-\ve}-h_j\quad\text{and}\quad1+\sigma(t)\le2h_j^{-\ve}.\]
By Lemma \ref{lemma:Schr},
\[\int_{S\M}e^{-\frac{i\sigma(t)}{4}U_+^2}f\,d\Vol_L=\int_{S\M}f\,d\Vol_L.\]
Hence,
\begin{align*}
\Wc_j\left(f\circ a_t\right)-\Wc_j(f)
&=\left[\Wc_j\left(f\circ a_t\right)-\Wc_j\left(e^{-\frac{i\sigma(t)}{4}U_+^2}f\right)\right]\\
&\quad+\left[\Wc_j\left(e^{-\frac{i\sigma(t)}{4}U_+^2}f\right)-\int_{S\M}e^{-\frac{i\sigma(t)}{4}U_+^2}f\,d\Vol_L\right]\\
&\quad-\left[\Wc_j(f)-\int_{S\M}f\,d\Vol_L\right].
\end{align*}
By Theorem \ref{thm:renorm}, the first term is bounded by
\[Ch_j^{1-13\ve}\|f\|_{H^{20}(S\M)}.\]
For the second term, we apply \eqref{eq:QErate} to $e^{-\frac{i\sigma(t)}{4}U_+^2}f$, which is in $C^\infty(S\M)$ by (i) of Lemma \ref{lemma:Schr}. By (ii) of Lemma \ref{lemma:Schr}, the second term is bounded by
\[Ch_j^{\delta_1}\left\|e^{-\frac{i\sigma(t)}{4}U_+^2}f\right\|_{H^{n_1}(S\M)}\le Ch_j^{\delta_1}\left(1+\sigma(t)\right)^{2n_1}\|f\|_{H^{3n_1}(S\M)}\le Ch_j^{\delta_1-2n_1\ve}\|f\|_{H^{3n_1}(S\M)},\]
where $C=C(n_1)>0$.

By \eqref{eq:QErate}, the third term is bounded by
\[Ch_j^{\delta_1}\|f\|_{H^{n_1}(S\M)}\le Ch_j^{\delta_1-2n_1\ve}\|f\|_{H^{n_1}(S\M)}.\]
This completes the proof.
\end{proof}

\begin{proof}[Proof that (i) implies (ii) of Theorem \ref{thm:rate}]
Suppose that (i) holds with $\delta_1$ and $n_1$. Since \eqref{eq:QErate} remains valid if $\delta_1$ is replaced by a smaller positive number and $n_1$ by a larger integer, we may assume that $0<\delta_1\le1$ and $n_1\ge1$. Proposition \ref{prop:converse} with $\ve=0$ gives (ii) with $\delta_2=\delta_1$ and $n_2=\max\{20,3n_1\}$. Moreover, Proposition \ref{prop:converse} with
\[0<\ve<\min\left\{\frac{1}{13},\frac{\delta_1}{2n_1}\right\}\]
gives \eqref{eq:Egorovbeyond}. This proves the result.
\end{proof}

\appendix

\section{Models that fail QUE}\label{sec:non}
In this appendix, we compare our results with two models in which QUE fails: the quantum cat maps and the Bunimovich stadium. For each model, we clarify in which quantization and up to which time Egorov's theorem holds, and how it is related to QUE. Quantum ergodicity holds in both models; see Degli Esposti-Graffi-Isola \cite{DEGI}, Bouzouina-De Bi\`evre \cite{BDB}, and Zelditch \cite{Ze5} for the quantum cat maps, and G\'erard-Leichtnam \cite{GL} and Zelditch-Zworski \cite{ZZ} for the Bunimovich stadium.

On compact hyperbolic surfaces, Theorem \ref{thm:equiv} relies on the following three ingredients.
\begin{enumerate}[(1)]
\item \textit{A quantization adapted to the weak stable foliation.} In Zelditch's quantization, the remainder in Egorov's theorem for eigenfunctions is given explicitly in \eqref{eq:exact}. By Theorem \ref{thm:commuD}, in the commutator of $\Op_h(f)$ with $h^2(\Delta-\frac14)$, the phase of the hyperbolic waves contributes the transport term $-2ihrQf$, and the symbol contributes the term $-h^2(Q^2+U_+^2)f$, which gives the remainder.
\item \textit{Uniform hyperbolicity.} Backward in time,
\[h_j^2U_+^2\left(f\circ a_t\right)=h_j^2e^{2|t|}\left(U_+^2f\right)\circ a_t\quad\text{for }t\le0.\]
This is much smaller than $h_jQ(f\circ a_t)$ before the Ehrenfest time, and comparable to it at $-T_E$, where $f\circ a_t$ oscillates at the scale $h_j^{\frac12}$ along the stable horocycles. So the remainder in Egorov's theorem is small before $-T_E$ for every eigenbasis. At $-T_E$, it is of order $1$, and Egorov's theorem requires a cancellation in the matrix elements of the eigenfunctions.
\item \textit{Unique ergodicity of the stable horocycle flow.} At $-T_E$, the limit of the remainder in Egorov's theorem is given by the Schr\"odinger group along the stable horocycles in \eqref{eq:defect}. It vanishes for all symbols if and only if the semiclassical measure is invariant under the stable horocycle flow. By the unique ergodicity of this flow, this holds if and only if the semiclassical measure is $d\Vol_L$.
\end{enumerate}
For the quantum cat maps, (2) and (3) hold, while (1) holds in a quantization adapted to the stable foliation of the cat map, but not in the Weyl quantization. For the Bunimovich stadium, (2) fails, and there are no analogues of (1) and (3).

\subsection{Quantum cat maps}
Let $\T^2=\R^2/\Z^2$, and denote by $dx$ the Lebesgue measure on $\T^2$. Suppose that
\[M=\begin{pmatrix}
a & b\\
c & d
\end{pmatrix}\in\SL(2,\Z)\quad\text{with }\mathrm{tr}\,M>2.\]
Then $M$ has eigenvalues $e^{\pm\chi}$ with the Lyapunov exponent $\chi>0$. It induces a hyperbolic automorphism of $\T^2$, also denoted by $M$, i.e., a cat map; see Arnold-Avez \cite{AA} and Katok-Hasselblatt \cite{KH}. The stable and unstable directions of $M$ are the eigendirections of $e^{-\chi}$ and $e^\chi$, respectively. Under the parity condition
\[ab=cd=0\pmod2,\]
Hannay-Berry \cite{HB} quantized $M$ as a unitary operator $\hat M$ on an $N$-dimensional Hilbert space $\Hc_N$ for each $N\in\N$ with $N\ge1$, with the semiclassical parameter $h=\frac{1}{2\pi N}$; see also Kurlberg-Rudnick \cite{KR1}. The Weyl quantization of $f\in C^\infty(\T^2)$ is
\[\Op_N^\Wr(f)=\sum_{k\in\Z^2}f_k\hat T_N(k),\quad\text{where }f(x)=\sum_{k\in\Z^2}f_ke^{2\pi ik\cdot x}.\]
Here, $\hat T_N(k)$ is the quantization of the character $e^{2\pi ik\cdot x}$, which acts on $\Hc_N$ as a quantum translation in the phase space $\T^2$. QUE for $\hat M$ means that for every sequence of normalized eigenstates $\vp_j\in\Hc_{N_j}$ of $\hat M$ with $N_j\to\infty$,
\[\lim_{j\to\infty}\left\langle\Op_{N_j}^\Wr(f)\vp_j,\vp_j\right\rangle_{\Hc_{N_j}}=\int_{\T^2}f\,dx\quad\text{for each }f\in C^\infty(\T^2).\]
Equivalently, every semiclassical measure of the eigenstates of $\hat M$ is $dx$.

\subsubsection{The Weyl quantization}
Since $M$ is linear, Egorov's theorem in the Weyl quantization is exact for all times \cite{HB}, that is,
\[\hat M^{-t}\Op_N^\Wr(f)\hat M^t=\Op_N^\Wr\left(f\circ M^t\right)\quad\text{for }t\in\Z\text{ and }f\in C^\infty(\T^2).\]
So if $\vp\in\Hc_N$ is an eigenstate of $\hat M$, then
\[\left\langle\Op_N^\Wr\left(f\circ M^t\right)\vp,\vp\right\rangle_{\Hc_N}=\left\langle\Op_N^\Wr(f)\vp,\vp\right\rangle_{\Hc_N}\quad\text{for }t\in\Z.\]
That is, Egorov's theorem holds for every eigenstate at all times, including the Ehrenfest time.

Nevertheless, QUE fails. Faure-Nonnenmacher-De Bi\`evre \cite{FNDB} constructed a sequence $N_j\to\infty$ and eigenstates $\vp_j\in\Hc_{N_j}$ of $\hat M$ such that
\begin{equation}\label{eq:FNDB}
\lim_{j\to\infty}\left\langle\Op_{N_j}^\Wr(f)\vp_j,\vp_j\right\rangle_{\Hc_{N_j}}=\frac12f(0)+\frac12\int_{\T^2}f\,dx\quad\text{for each }f\in C^\infty(\T^2),
\end{equation}
that is, these eigenstates are half-scarred at the fixed point $0$. Here, the quantum periods of $\hat M$ on $\Hc_{N_j}$ are of order $\log N_j$, which results in large multiplicities of the eigenvalues. The weight $\frac12$ is maximal by Faure-Nonnenmacher \cite{FN}. On the other hand, Kurlberg-Rudnick proved that QUE holds for the joint eigenstates of $\hat M$ and its Hecke operators \cite{KR1}, and along a density one sequence of $N$ \cite{KR2}; see also Bourgain \cite{Bo}.

Therefore, Egorov's theorem in the Weyl quantization does not distinguish these cases, i.e., it does not detect QUE. It is the analogue of (i) of Theorem \ref{thm:Egorov}, which holds for every eigenbasis, rather than of Theorem \ref{thm:equiv}. In terms of the three ingredients above, (2) and (3) hold for the cat maps: $M$ is uniformly hyperbolic, and the translation flow along the stable direction is a linear flow on $\T^2$ with an irrational slope, which is uniquely ergodic. However, (1) fails: the Weyl quantization is not adapted to any Lagrangian foliation, and it intertwines $M$ exactly. So there is no remainder in Egorov's theorem whose cancellation could detect QUE.

\subsubsection{A quantization adapted to the stable foliation}
In \cite{Han2}, we introduce a quantization $\Op_N^+$ for the quantum cat maps, which is the analogue of Zelditch's quantization. It is adapted to the stable foliation of $M$: on the universal cover $\R^2$, the operators act on the Lagrangian states associated with the lifted stable leaves by multiplication by the restriction of the symbols to the leaves. For each $f\in C^\infty(\T^2)$ independent of $h$, the operators $\Op_N^+(f)$ and $\Op_N^\Wr(f)$ differ by $O_f(h)$ in the operator norm. So QUE does not depend on the choice between the two quantizations. However, in $\Op_N^+$, Egorov's theorem has a remainder, which is given exactly by the Schr\"odinger group along the stable foliation: for $f\in C^\infty(\T^2)$ and $t\in\Z$,
\begin{equation}\label{eq:Ecatplus}
\Op_N^+\left(f\circ M^t\right)=\hat M^{-t}\Op_N^+\left(e^{-\frac{i\sigma(t)}{4}U_+^2}f\right)\hat M^t,\quad\text{where }\sigma(t)=\frac{b}{\sinh\chi}h\left(e^{-2\chi t}-1\right).
\end{equation}
Here, $U_+$ is the derivative along the stable eigenvector of $M$ whose first coordinate is $1$, and $b\ne0$ is the upper right entry of $M$. This is the analogue of \eqref{eq:exactSchr}, without the integral, since $M$ is linear. Forward in time, $|\sigma(t)|<|b|h/\sinh\chi$. Backward in time, $|\sigma(t)|$ grows, and it is of order $1$ at the Ehrenfest time $-T_E$, where
\[T_E=\left\lfloor\frac{|\log h|}{2\chi}\right\rfloor.\]
So Egorov's theorem holds forward in time for all times, and backward in time before the Ehrenfest time.

The analogue of Theorem \ref{thm:equiv} holds at $-T_E$ \cite{Han2}. Suppose that normalized eigenstates $\vp_j\in\Hc_{N_j}$ of $\hat M$, $j\in J$, induce a semiclassical measure $\nu$ on $\T^2$. As in Section \ref{sec:intro}, write
\[\Wc_j(f)=\left\langle\Op_{N_j}^+(f)\vp_j,\vp_j\right\rangle_{\Hc_{N_j}}\quad\text{for }f\in C^\infty(\T^2).\]
Assume that $\sigma(-T_E)\to\sigma_0$ as $j\in J$, $j\to\infty$, where $h=\frac{1}{2\pi N_j}$ in $\sigma$ and $T_E$. Then for $f\in C^\infty(\T^2)$,
\[\lim_{j\in J,j\to\infty}\left(\Wc_j\left(f\circ M^{-T_E}\right)-\Wc_j(f)\right)=\int_{\T^2}\left(e^{-\frac{i\sigma_0}{4}U_+^2}f-f\right)d\nu,\]
and these limits determine $\nu$. In particular, Egorov's theorem holds at $-T_E$ along $J$ if and only if $\nu=dx$, in which case it holds for all times along $J$. So in the quantization $\Op_N^+$, Egorov's theorem at the Ehrenfest time detects whether QUE holds or fails.

For the half-scarred eigenstates in \eqref{eq:FNDB}, i.e., for $\nu=\frac12\delta_0+\frac12dx$, the limit of the remainder is
\[\frac12\left(\left(e^{-\frac{i\sigma_0}{4}U_+^2}f\right)(0)-f(0)\right),\]
which does not vanish for some $f\in C^\infty(\T^2)$. In \cite{Han2}, we compute this limit directly from the construction of the eigenstates in \cite{FNDB}, and then recover \eqref{eq:FNDB}, since the limits determine $\nu$. On the other hand, Egorov's theorem holds for all times for the Hecke eigenstates in \cite{KR1}, and for all eigenstates when $N$ is prime and the period of $M$ modulo $N$ is much larger than $N^{\frac12}$, as in \cite{KR2}.

Since $M$ is linear, several steps of this paper simplify for the quantum cat maps. The remainder in Egorov's theorem is given exactly by \eqref{eq:Ecatplus}, and forward in time, Egorov's theorem holds for all times. Moreover, the only distributions on $\T^2$ which are invariant under the linear flow along the stable direction are the multiples of $dx$. So the positivity of the matrix elements in Lemma \ref{lemma:positive} is not needed, and QUE with a polynomial rate is equivalent to Egorov's theorem at $-T_E$ with a polynomial remainder with the same exponent \cite{Han2}.

\subsection{The Bunimovich stadium}
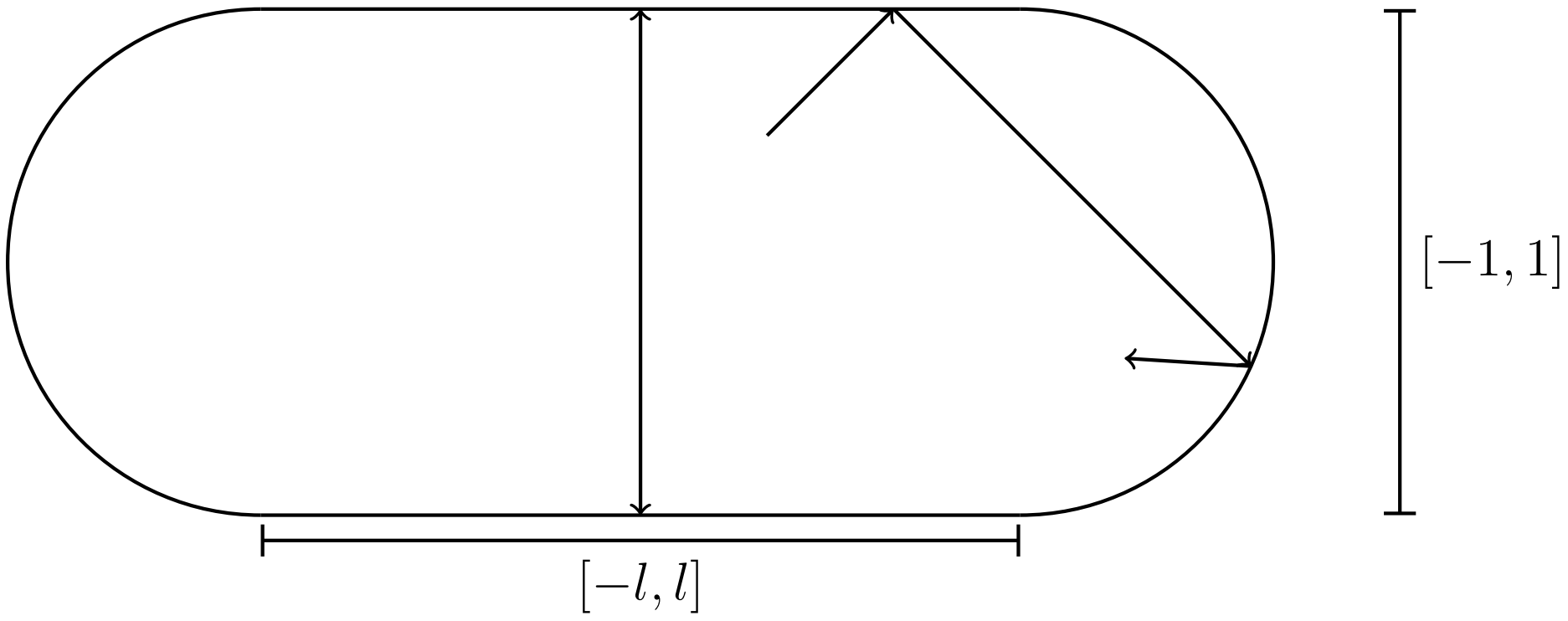
\begin{figure}[htb]
\centering
\begin{tikzpicture}
\draw[-,thick] (-3,2) -- (3,2);
\draw[-,thick] (-3,-2) -- (3,-2);
\draw[thick] (3,-2) arc (270:450:2cm);
\draw[thick] (-3,2) arc (90:270:2cm);
\draw[|-|,thick] (-3,-2.2) -- (3,-2.2) node[midway,below] {$[-l,l]$};
\draw[|-|,thick] (6,-2) -- (6,2) node[midway,right] {$[-1,1]$};

\draw[->,thick] (1,1) -- (2,2);
\draw[->,thick] (2,2) -- (4.823,-0.823);
\draw[->,thick] (4.823,-0.823) -- (3.825,-0.76);
\draw[<->,thick] (0,-2) -- (0,2);
\end{tikzpicture}
\caption{The billiard flow on the Bunimovich stadium $\Omega_l$ and a bouncing ball trajectory (vertical)}\label{fig:stadium}
\end{figure}

Let $l>0$. Denote by $\Omega_l\subset\R^2$ the union of the rectangle $R_l=[-l,l]\times[-1,1]$ and the two half disks of radius $1$ attached to its vertical sides, i.e., the Bunimovich stadium; see Figure \ref{fig:stadium}. Denote by $a_t$ the billiard flow on $S\Omega_l=\Omega_l\times\Sb$, i.e., the straight motion at unit speed in $\Omega_l$ with reflection at the boundary. Bunimovich \cite{Bun} proved that $a_t$ is ergodic with respect to the (normalized) Liouville measure $d\Vol_L$ on $S\Omega_l$; see Chernov-Markarian \cite[Chapter 8]{CM} for more discussion.

Let $\{\vp_j\}_{j=1}^\infty$ be an orthonormal basis of eigenfunctions of the non-negative Laplacian on $\Omega_l$ with the Dirichlet boundary condition, with eigenvalues $\mu_j\to\infty$. Write
\[h_j=\mu_j^{-\frac12}.\]
QUE on $\Omega_l$ means that for every such eigenbasis,
\[\lim_{j\to\infty}\left\langle\Op_{h_j}(f)\vp_j,\vp_j\right\rangle_{L^2(\Omega_l)}=\int_{S\Omega_l}f\,d\Vol_L\]
for each symbol $f$ independent of $h$ and compactly supported over the interior of $\Omega_l$, where $\Op_h$ is a standard semiclassical quantization on $\R^2$.

\subsubsection{Egorov's theorem}
In a standard semiclassical quantization, Egorov's theorem holds for $|t|\le T$ with $T$ independent of $h$, away from the trajectories tangent to the boundary; see Farris \cite{Fa}, and also G\'erard-Leichtnam \cite{GL} and Zelditch-Zworski \cite{ZZ}, where such a version is used to prove quantum ergodicity. To the author's knowledge, it is unknown whether Egorov's theorem holds at logarithmic times $|t|=O(|\log h|)$ on the stadium.

\subsubsection{Relation to QUE}
Hassell \cite{Has} proved that QUE fails for almost every $l\in[1,2]$: there is a sequence of eigenfunctions whose semiclassical measure has positive mass on the bouncing ball trajectories, i.e., the vertical trajectories in $R_l$; see Figure \ref{fig:stadium}. The set of these trajectories is invariant under $a_t$ and has zero Liouville measure. It is believed that QUE fails on every stadium, but this remains open. Since Egorov's theorem at logarithmic times is unknown, it is unclear whether an analogue of Theorem \ref{thm:equiv} holds on the stadium.

The decisive feature here is the lack of uniform hyperbolicity. The billiard flow $a_t$ is ergodic but not uniformly hyperbolic: the bouncing ball trajectories form a family of periodic trajectories, and the nearby trajectories in $R_l$ separate from them only linearly in time, rather than exponentially. The non-QUE eigenfunctions of \cite{Has} concentrate exactly on this non-hyperbolic family, and the result extends to the partially rectangular billiards, which all contain such a family; see the appendix of \cite{Has}. Accordingly, for the billiard map, i.e., the return map of $a_t$ to the boundary, the correlations of H\"older functions decay only at the rate $O(n^{-1})$, and the billiard flow is not expected to mix faster; see B\'alint-Melbourne \cite[Section 3]{BM}. This is in sharp contrast with the exponential decay of correlations for the geodesic flow on compact hyperbolic surfaces; see Ratner \cite{Rat}.

In our approach, the uniform hyperbolicity of the geodesic flow makes $f\circ a_t$ oscillate at the scale $e^{-|t|}$ along the stable horocycles, so the cancellation is required at the logarithmic time $-T_E$. Heuristically, near the bouncing ball trajectories, $f\circ a_t$ only develops oscillations at the scale $|t|^{-1}$, which reaches the scale $h^{\frac12}$ at the polynomial time $|t|\sim h^{-\frac12}$. So no cancellation is required near the bouncing ball trajectories at logarithmic times, which is consistent with the concentration of the eigenfunctions there. Thus, unlike for the quantum cat maps, the obstruction to an analogue of Theorem \ref{thm:equiv} on the stadium is dynamical, rather than due to the choice of quantization.

\section*{Acknowledgements}
I have benefited greatly from discussions with Suresh Eswarathasan, Andrew Hassell, Bingrong Huang, Long Jin, St\'ephane Nonnenmacher, Ze\'ev Rudnick, Shu Shen, Melissa Tacy, Steve Zelditch, and Zhiyuan Zhang. I am grateful to them for their valuable insights and continued support throughout this research project.

Part of this work was completed during my time as an invited researcher at the Institut des Hautes \'Etudes Scientifiques (IHES) in the spring of 2023. I thank IHES for its hospitality, the Chern Foundation for supporting my stay, and Maxim Kontsevich at IHES for his interest in this research.


\end{document}